%% file: agct11.tex
\documentclass{conm-p-l}
\usepackage{amssymb}
\usepackage{amsmath}
\usepackage{amsthm}
\usepackage{booktabs}
\usepackage{setspace}
\usepackage{array}
\usepackage{url}
\usepackage{color}
\usepackage{epic}
\usepackage[colorlinks=true,bookmarks=false,citecolor=webbrown,linkcolor=webbrown,filecolor=webcyan,urlcolor=webcyan]{hyperref}
\definecolor{webbrown}{rgb}{0,0,0}
\definecolor{webcyan}{rgb}{0,0,0}

\newcommand{\Fp}{\mathbb{F}_p}

\newcommand{\Fpk}{\mathbb{F}_{p^k}}

\newcommand{\ZetaCp}{Z(C/\mathbb{F}_p;T)}

\newcommand{\JCp}{J(C/\mathbb{F}_p)}

\newcommand{\GalQ}{{\rm Gal}(\overline{\mathbb{Q}}/\mathbb{Q})}
\newcommand{\barLp}{\bar{L}_p(T)}
\newcommand{\Aut}{{\rm Aut}}
\newcommand{\EndC}{{\rm End}_\mathbb{C}}
\newcommand{\Exp}{{\bf E}}
\newcommand{\mean}{{\rm m}}
\newcommand{\Dop}{{\bf D}}
\newcommand{\sgn}{{\rm sgn}}
\newcommand{\tr}{{\rm tr}}
\newcommand{\Aseq}[1]{\htmladdnormallink{A#1} {http://www.research.att.com/projects/OEIS?Anum=#1}}
\newcommand{\z}{\drawline}
\newtheorem{conj}{Conjecture}
\newtheorem{proposition}{Proposition}
\newtheorem{theorem}{Theorem}
\newtheorem{corollary}{Corollary}
\newtheorem{lemma}{Lemma}
\begin{document}

\title{Hyperelliptic Curves, $L$-polynomials, and Random Matrices}

\author{Kiran S. Kedlaya}
\thanks{Kedlaya was supported by 
NSF CAREER grant DMS-0545904 and a Sloan Research Fellowship.}
\author{Andrew V. Sutherland}
\address{Department of Mathematics \\ Massachusetts Institute of Technology \\
77 Massachusetts Avenue \\ Cambridge, MA 02139}
\email{(kedlaya|drew)@math.mit.edu}

\begin{abstract}
We analyze the distribution of unitarized $L$-polynomials $\barLp$
(as $p$ varies) obtained from a hyperelliptic curve of genus $g\le 3$
defined over $\mathbb{Q}$. In the generic case, we find experimental agreement
with a predicted correspondence (based on the Katz-Sarnak random matrix
model) between the distributions of $\barLp$ and of
characteristic polynomials of random matrices in the compact Lie
group $USp(2g)$. We then formulate an analogue of the Sato-Tate
conjecture for curves of genus 2, in which the generic distribution
is augmented by 22 exceptional distributions, each corresponding
to a compact subgroup of $USp(4)$. In every case,
we exhibit a curve closely matching the proposed
distribution, and can find no curves unaccounted for by our classification.
\end{abstract}
\keywords{Sato-Tate conjecture, trace of Frobenius, zeta function, Haar measure, moment sequence}
\subjclass[2000]{Primary 11G40; Secondary 11G30, 05E15}
\maketitle


\section{Introduction}

For $C$ a smooth projective curve of genus $g$ defined over $\mathbb{Q}$ and each prime $p$ where $C$ has good reduction, we consider the polynomial $L_p(T)$, the numerator of the zeta function $\ZetaCp$.  This polynomial is intimately related to many arithmetic properties of the curve, appearing in the Euler product of the $L$-series
$$L(C,s) = \prod_p L_p(p^{-s})^{-1},$$
the characteristic polynomial of the Frobenius endomorphism, 
$$\chi_p(T)=T^{2g}L_p(T^{-1}),$$
and the order of the group of $\Fp$-rational points on the Jacobian of $C$,
$$\#\JCp=L_p(1).$$

In genus 1, $C$ is an elliptic curve, and $L_p(T) = pT^2 - a_pT + 1$ is determined by the trace of Frobenius, $a_p$.  The distribution of $a_p$ as $p$ varies has been and remains a subject of considerable interest, forming the basis of several well known conjectures, including those of Lang-Trotter \cite{Lang:FrobeniusDistributions} and Sato-Tate \cite{Tate:SatoTate}.  Considerable progress has been made on these questions, particularly the latter, much of it quite recently \cite{Baier:SatoTateLangTrotter,Baier:SatoTateSmallAngles,Banks:SatoTateOnAverage,Harris:SatoTateProof,Mazur:SatoTate}.

In genus 2, $C$ is a hyperelliptic curve, and the study of $L_p(T)$ may be viewed as a natural generalization of these questions (we also consider hyperelliptic curves of genus 3).  Our first objective is to understand the shape of the distribution of $L_p(T)$, which leads us to focus primarily on Sato-Tate type questions.

The random matrix model developed by Katz and Sarnak provides a ready generalization of the Sato-Tate conjecture to higher genera.  They show (see Theorems 10.1.18.3 and 10.8.2 in \cite{Katz:RandomMatrices}) that over a universal family of hyperelliptic curves of genus $g$, the distribution of {\em unitarized} $L$-polynomials, $\barLp = L_p(p^{-1/2}T)$, corresponds to the distribution of characteristic polynomials $\chi(T)$ in the compact Lie group $USp(2g)$ (the group of $2g\times 2g$ complex matrices that are both unitary and symplectic).  By also considering infinite compact subgroups of $USp(2g)$, we are able to frame a generalization of the Sato-Tate conjecture applicable to any smooth curve defined over $\mathbb{Q}$.\footnote{With the generous assistance of Nicholas Katz.}

To test this conjecture, we rely on a collection of highly efficient algorithms to compute $L_p(T)$, described by the authors in \cite{KedlayaSutherland:HyperellipticLSeries}.  The performance of these algorithms has improved dramatically in recent years (due largely to an interest in cryptographic applications \cite{Cohen:HECHECC}). For a hyperelliptic curve, it is now entirely practical to compute $L_p(T)$ for all $p\le N$ (where the curve has good reduction), with $N$ on the order of $10^8$ in genus 2 and more than $10^7$ in genus 3.  Alternatively, for much smaller $N$ (less than $10^4$) one can perform similar computations with $10^{10}$ curves or more.

We characterize the otherwise overwhelming abundance of data with {\em moment statistics}.  If $\{x_p\}$ is a set of unitarized values derived from $\barLp$, say $x_p=a_p/\sqrt{p}$, we compute the first several terms of the sequence
\begin{equation}\label{expr:MeanSequence}
\mean(x_p),\enspace \mean(x_p^2),\enspace \mean(x_p^3),\enspace\ldots,
\end{equation}
where $\mean(x_p^k)$ denotes the mean of $x_p^k$ over $p$.  Under the conjecture, the moment statistics converge, term by term, to the {\em moment sequence}
\begin{equation}\label{expr:MomentSequence}
\Exp[X],\enspace \Exp[X^2],\enspace \Exp[X^3],\enspace \ldots,
\end{equation}
where the corresponding random variable $X$ is derived from the characteristic polynomial $\chi(T)$ of a random matrix $A$, say $X=\tr(A)$.  In all the cases of interest to us, the moments $\Exp[X^n]$ exist and determine the distribution of $X$.  Furthermore, they are integers.  The distributions we encounter can typically be distinguished by the first eight terms of ($\ref{expr:MomentSequence}$).  It will be convenient to begin our sequences with $\Exp[X^0]=1$.

To apply this approach we must determine these moment sequences explicitly.  This is an interesting problem in its own right, with applications to representation theory and combinatorics.  The particular cases we consider include the elementary and Newton symmetric functions (power sums) of the eigenvalues of a random matrix in $USp(2g)$.  We derive explicit formulae for the corresponding moment generating functions.  Our results intersect other work in this area, most notably that of Grabiner and Magyar \cite{Grabiner:WeylChamber}, and also Eric Rains \cite{Rains:EigenvalueDistribution}.  We take a somewhat different approach, relying on the Haar measure to encode the combinatorial structure of the group.

With moment sequences for $USp(2g)$ in hand, we then survey the distributions of $\barLp$.  The case $g=1$ is easily described, and we do so here.  The polynomial $\barLp=p+a_1T+1$ is determined by the coefficient $a_1=-a_p/\sqrt{p}$, and there are two distributions of $a_1$ that arise.

For elliptic curves without complex multiplication, the moment statistics of $a_1$ converge to the corresponding moment sequence in $USp(2)$:
\begin{align*}
1,\enspace 0,\enspace 1,\enspace 0, \enspace 2, \enspace 0, \enspace 5, \enspace 0, \enspace 14, \enspace 0, \enspace 42,\enspace \ldots,
\end{align*}
whose $(2n)$th term is the $n$th Catalan number.  Convergence follows from the Sato-Tate conjecture, which for curves with multiplicative reduction at some prime (almost all curves) is now proven, thanks to the work of Clozel, Harris, Shepherd-Barron, and Taylor \cite{Clozel:SatoTate,Harris:SatoTateProof,Taylor:SatoTate} (see Mazur \cite{Mazur:SatoTate} for an overview).  Testing at least one example of every curve with conductor less than $10^7$ (including all the curves in Cremona's tables \cite{Cremona:Database,SteinWatkins:Database}) revealed no apparent exceptions among curves with purely additive reduction.

For elliptic curves with complex multiplication, the moment statistics of $a_1$ converge instead to the sequence:
\begin{align*}
1,\enspace 0,\enspace 1,\enspace 0, \enspace 3, \enspace 0, \enspace 10, \enspace 0, \enspace 35, \enspace 0, \enspace 126,\enspace \ldots,
\end{align*}
whose $(2n)$th term is $\binom{2n}{n}/2$ for $n>0$.  This is the moment sequence of a compact subgroup of $USp(2)$, specifically, the normalizer of $SO(2)$ in $SU(2)=USp(2)$.  The elements with nonzero traces have uniformly distributed eigenvalue angles, and for elliptic curves with complex multiplication, convergence follows from a theorem of Deuring
\cite{Deuring:ComplexMultiplication} and known equidistribution results for Hecke characters \cite[Ch. XV]{Lang:AlgebraicNumberTheory}.  The only other infinite compact subgroup of $USp(2)$ (up to conjugacy) is $SO(2)$, and its moment sequence does not appear to correspond to the moment statistics of any elliptic curve (such a curve would necessarily contradict the Sato-Tate conjecture).

The main purpose of the present work is to undertake a similar study in genus 2.  We also lay some groundwork for genus 3, but consider only the case of a typical hyperelliptic curve.  Already in genus 2 we find a much richer set of possible distributions.

For the typical hyperelliptic curve in genus 2 (resp. 3), when computed for suitably large $N$, the moment statistics closely match the corresponding moment sequences in $USp(4)$ (resp. $USp(6)$), as predicted.  In genus 2 we can make a stronger statement.  In a family of one million curves with randomly chosen coefficients, every single one appeared to have the $\barLp$ distribution of characteristic polynomials in $USp(4)$.  Under reasonable assumptions, we can reject alternative distributions with a high level of statistical confidence.

To find exceptional distributions in genus 2 we must cast our net wider, searching very large families of curves with constrained coefficient values, as well as examples taken from the literature.  Such a search yielded over 30,000 nonisomorphic exceptional curves, but among these we find only 22 clearly distinct distributions, each with integer moments (Table \ref{table:g2a1dist}, p. \pageref{table:g2a1dist}).  For each of these distributions we are able to identify a specific compact subgroup $H$ of $USp(4)$ with a matching distribution of characteristic polynomials (Table \ref{table:USp4groups}, p. \pageref{table:USp4groups}).  The method we use to construct these subgroups is quite explicit, and a converse statement is very nearly true.  Of the subgroups we can construct, only two do not correspond to a distribution we have found; we can rule out one of these
and suspect the other can be ruled out also (\S~\ref{section:nonexistence}).  We believe we have accounted for all possible $L$-polynomial distributions of a genus 2 curve.

\section{Some Motivating Examples}\label{section:preliminaries}

We work throughout with smooth, projective, geometrically irreducible algebraic curves defined over $\mathbb{Q}$.
Recall that for a curve $C$ with good reduction at $p$, the zeta function $\ZetaCp$ is defined by the formal power series
\begin{equation*}
\ZetaCp = \exp\left(\sum_{k=1}^{\infty}N_kT^k/k\right),
\end{equation*}
where $N_k$ counts the (projective) points on $C$ over $\Fpk$.  From the seminal work of Emil Artin \cite{Artin:ZetaFunction}, we know that $\ZetaCp$ is a rational function of the form
\begin{equation*}
\ZetaCp = \frac{L_p(T)}{(1-T)(1-pT)},
\end{equation*}
where the monic polynomial $L_p(T)\in\mathbb{Z}[T]$ has degree $2g$ ($g$ is the genus of $C$) and constant coefficient 1. 

By the Riemann hypothesis for curves (proven by Weil \cite{Weil:ZetaFunction}), the roots of $L_p(T)$ lie on a circle of radius $p^{-1/2}$ about the origin of the complex plane.  To study the distribution of $L_p(T)$ as $p$ varies, we use the unitarized polynomial $$\barLp = L_p(p^{-1/2}T),$$
which has roots on the unit circle.  As $\barLp$ is a real polynomial of even degree with $\bar{L}_p(0)=1$, these roots occur in conjugate pairs.  We may write
$$\barLp = T^{2g} + a_1T^{2g-1} + a_2T^{2g-2} + \cdots + a_2T^2 + a_1T + 1.$$
Since $\barLp$ has unitary roots, we know that
\begin{equation}\label{equation:WeilInterval}
|a_i|\le\binom{2g}{i},
\end{equation}
and ask how $a_i$ is distributed within this interval as $p$ varies.

The next three pages show the distribution of $a_1$ for arbitrarily chosen curves of genus 1, 2, and 3 and various values of $N$.  The coefficient $a_1$ is the negative sum of the roots of $\barLp$, and may be written as $a_1= -a_p/\sqrt{p}$, where $a_p$ is the trace of Frobenius.

Each graph represents a histogram of nearly $\pi(N)$ samples (one for each prime where $C$ has good reduction) placed into approximately $\sqrt{\pi(N)}$ buckets which partition the interval $[-2g,2g]$ determined by (\ref{equation:WeilInterval}).  The horizontal axis spans this interval, and the vertical axis has been suitably scaled, with the height of the uniform distribution, $1/(4g)$, indicated by a dotted line.

\clearpage
\include{converge}

The familiar semicircular shape in genus 1 is the Sato-Tate distribution.  The examples in higher genera also appear to converge to distinct distributions.  Provided the curve is ``typical" (a notion we will define momentarily) this distribution is the same for every curve of a given genus.  Even in atypical cases, there is (empirically) a small distinct set of distributions that arise for a given genus.  In genus 1, only one exceptional distribution for $a_1$ is known, exhibited by all curves with complex multiplication:

\vspace{12pt}
\input{g1_excep}
\begin{center}
\small
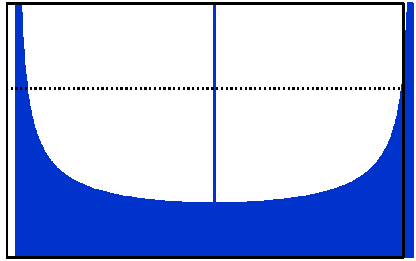
Fig. 1: Distribution of $a_1$ for $y^2 = x^3 - 15x + 22$.
\normalsize
\end{center}
\vspace{8pt}

The central spike has area 1/2 (asymptotically), arising from the fact that a curve with CM-field $\mathbb{Q}[\sqrt{D}]$ has $a_p=0$ precisely when $D$ is a not a quadratic residue in $\Fp$ \cite{Deuring:ComplexMultiplication}.

In higher genera, a richer set of exceptional distributions arises.  Below is an example for a genus 2 curve whose Jacobian splits as the product of two elliptic curves (one with complex multiplication).  Histograms for several other exceptional genus 2 distributions are provided in Appendix II.

\vspace{12pt}
\input{g2_excep}
\begin{center}
\small
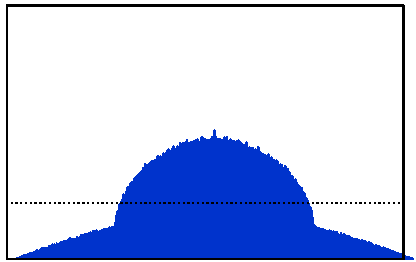
Fig. 2: Distribution of $a_1$ for $y^2 = x^5 + 20x^4 - 26x^3 + 20x^2 + x$.
\normalsize
\end{center}

\section{A Generalized Sato-Tate Conjecture}\label{section:SatoTate}

We wish to give a conjectural basis for these distributions, both in the typical and atypical cases.  The formulation presented here follows the model developed by Katz and Sarnak \cite{Katz:RandomMatrices} and relies heavily on additional detail provided by Nicholas Katz \cite{Katz:PersonalCommunication}, whom we gratefully acknowledge.  Most of the statements below readily generalize to abelian varieties, but we restrict ourselves to curves defined over $\mathbb{Q}$.

We begin with the Sato-Tate conjecture, which may be stated as follows:

\begin{conj}[Sato-Tate]\label{conjecture:SatoTate}
For an elliptic curve without complex multiplication, the distribution of the roots $e^{i\theta}$ and $e^{-i\theta}$ of $\barLp$ for $p\le N$ converges (as $N\to\infty$) to the distribution given by the measure $\mu = \frac{2}{\pi}\sin^2\theta d\theta$ over $\theta\in[0,\pi]$.
\end{conj}

As noted earlier, this has been proven for elliptic curves with multiplicative reduction at some prime \cite{Harris:SatoTateProof}.  Using $a_1=-2\cos\theta$, one finds that
$$\Pr[a_1\le x] = \frac{1}{2\pi}\int_{-2}^x\sqrt{4-t^2}dt,$$
giving the familiar semicircular distribution.

An equivalent formulation of Conjecture \ref{conjecture:SatoTate} is that the distribution of $\barLp$ corresponds to the distribution of the characteristic polynomial of a random matrix in $USp(2)$.  More generally, the Haar measure on $USp(2g)$ provides a natural distribution of unitary symplectic polynomials of degree $2g$: the eigenvalues of a unitary matrix lie on the unit circle and the symplectic condition ensures that they occur in conjugate pairs (giving a polynomial with real coefficients).  Conversely, each unitary symplectic polynomial corresponds to a conjugacy class of matrices in $USp(2g)$.

Let the eigenvalues of a random matrix in $USp(2g)$ be $e^{\pm i\theta_1}$, \ldots, $e^{\pm i\theta_g}$, with $\theta_j \in [0,\pi]$.  The joint probability density function on $(\theta_1,\ldots,\theta_g)$ given by the Haar measure on $USp(2g)$ is 
\begin{equation}\label{equation:EigenvalueDistribution}
\mu(USp(2g)) = \frac{1}{g!}\Bigl(\prod_{j<k}(2\cos\theta_j - 2\cos\theta_k)\Bigr)^2\prod_j\left(\frac{2}{\pi}\sin^2\theta_j d\theta_j\right),
\end{equation}
as shown by Weyl \cite[Thm. 7.8B, p. 218]{Weyl:ClassicalGroups} (also see \cite[5.0.4, p. 107]{Katz:RandomMatrices}).  For $g=1$, we obtain the Sato-Tate distribution above.

In view of the atypical examples in the previous section, we cannot expect every curve to achieve the distribution given by $\mu(USp(2g))$.  One might impose a restriction comparable to that of Sato-Tate by requiring that $\EndC(J(C))=\mathbb{Z}$, i.e. that the Jacobian have minimal endomorphism ring.  While necessary, it is not clear that this restriction is sufficient in general. A stronger condition uses the $\ell$-adic representation of $\GalQ$ induced by the Galois action on the Tate module $T_\ell(C)$ (the inverse limit of the $\ell^n$-torsion subgroups of $J(C)$).  Specifically, we require that the image of the representation
\begin{equation*}
\rho_\ell:\GalQ\rightarrow \Aut(T_\ell(C))\cong GL(2g,\mathbb{Z}_\ell)
\end{equation*}
be Zariski dense in $GSp(2g,\mathbb{Z}_\ell)\subset GL(2g,\mathbb{Z}_\ell)$.  We know from results of Serre (\cite[Sec. 7, Thm 3]{Serre:GSp26odd} and \cite[p. 104]{Serre:CourseNotes1986}) that if this is achieved for any $\ell$, then it holds for all $\ell$, and we say such a curve has {\em large Galois image}.\\
\\
\noindent Each of the conditions below suffices for $C$ to have large Galois image.
\begin{enumerate}
\item
$C$ is a genus $g$ curve with $g$ odd, 2, or 6, and $\EndC(J(C))=\mathbb{Z}$ (Serre \cite{Serre:GSp26odd,Serre:MTGSp26odd}). This does not hold in genus 4 (Mumford \cite{Mumford:Genus4Example}).
\item
$C$ is a hyperelliptic curve $y^2=f(x)$ with $f(x)$ of degree $n\ge 5$ and the Galois group of $f(x)$ is isomorphic to $S_n$ or $A_n$ (Zarhin \cite{Zarhin:HyperellipticNoCM}).
\item
$C$ is a genus $g$ curve with good reduction outside of a set of primes $S$ and 
the mod $\ell$ reduction of the image of $\rho_\ell$ is equal to $GSp(2g,\mathbb{Z}/\ell\mathbb{Z})$ for some $\ell \ge c(g,S)$ (Faltings \cite{Faltings:Theorem}, Serre \cite{Serre:GSp26odd,Serre:CourseNotes1986}).
\end{enumerate}

Condition 3 was suggested to us by Katz.  The constant $c(g,S)$ depends only on $g$ and $S$.  From Faltings \cite[Thm. 5]{Faltings:Theorem}, we know that for a given $g$ and $S$, there are only finitely many nonisomorphic Jacobians of genus $g$ curves with good reduction outside of $S$.  Each such curve has large Galois image if and only if for all $\ell>c$ the mod $\ell$ image of $\rho_\ell$ is $GSp(2g,\mathbb{Z}/\ell\mathbb{Z})$, for some constant $c$.  By the results of Serre cited above, it suffices to find one such $\ell$, and taking the maximum of $c$ over the finite set of Jacobians gives the desired $c(g,S)$.  At present, effective bounds on $c(g,S)$ are known only in genus 1 \cite{Cojocaru:SurjectiveGalois,Kraus:EffectiveBound}.  Even without effective bounds, Condition 3 gives an easily computable heuristic that is useful in practice.\footnote{A fourth condition has recently been proven by Hall \cite{Hall:OpenImageTheorem}.  In the same paper, Kowalski proves that almost all hyperelliptic curves have large Galois image.}

We now consider the situation for a curve which does not have large Galois image.  The simplest case is an elliptic curve with complex multiplication.  For such a curve the distribution of $a_1$ clearly does not match the distribution of traces in $USp(2)$ (see Fig. 1 above).  In particular, the density of primes for which $a_1=-a_p/\sqrt{p}=0$ is one half.  There is, however, a compact subgroup of $USp(2)$ which gives the correct distribution.  Consider the subgroup
\begin{equation*}
H=
\left\{\left(\begin{matrix}
\cos{\theta} &\sin{\theta}\\
-\sin{\theta} &\cos{\theta}
\end{matrix}
\right),\left(
\begin{matrix}
i\cos{\theta} & i\sin{\theta}\\
i\sin{\theta} & -i\cos{\theta}\\
\end{matrix}
\right): \theta\in[0,2\pi]\right\}.
\end{equation*}
$H$ contains $SO(2)$ as a subgroup of index 2 (it is in fact the normalizer of $SO(2)$ in $USp(2)$).  The elements of $H$ not in $SO(2)$ all have zero traces, giving the desired density of 1/2.  The Haar measure on $SO(2)$ gives uniformly distributed eigenvalue angles, matching the distribution of nonzero traces for elliptic curves with complex multiplication.  Note that $H$ is disconnected.  This is a common (but not universal) feature among the subgroups we wish to consider.

For an elliptic curve with complex multiplication, the mod $\ell$ Galois image lies in the normalizer of a Cartan subgroup (see Lang \cite[Thm 3.2]{Lang:ModularForms}). For sufficiently large $\ell$, one finds that  is in fact equal to the normalizer of a Cartan subgroup.  There is then a correspondence with the subgroup $H$ above, which is the normalizer of the Cartan subgroup $SO(2)$ in $SU(2)=USp(2)$.

In general, we anticipate a relationship between the $\ell$-adic Galois image of a curve $C$, call it $G_\ell$ (a subgroup of $GSp(2g,\mathbb{Z}_\ell)$), and a compact subgroup $H$ of $USp(2g)$ whose distribution of characteristic polynomials matches the distribution of unitarized $L$-polynomials $\bar{L}_p(T)$ of $C$.  One can (conjecturally) describe this relationship, albeit in a  nonexplicit manner.  Briefly, one takes the Zariski closure of $G_\ell$ through a series of embeddings:
$$GSp(2g,\mathbb{Z}_\ell)\rightarrow GSp(2g,\mathbb{Q}_\ell)\rightarrow GSp(2g,\overline{\mathbb{Q}}_\ell)\rightarrow GSp(2g,\mathbb{C}).$$
The last step is justified by the fact that $\overline{\mathbb{Q}}_\ell$ and $\mathbb{C}$ are algebraically closed fields containing $\mathbb{Q}$ of equal cardinality, hence isomorphic, and we choose a particular embedding of $\overline{\mathbb{Q}}_\ell$ in $\mathbb{C}$. One then takes the intersection with $Sp(2g)$, obtaining a reductive group over $\mathbb{C}$.  After dividing by $\sqrt{p}$, the image of each Frobenius element lies in this intersection, as a diagonalizable matrix with unitary eigenvalues.  We now consider a maximal compact subgroup of the intersection, $H$, lying in $USp(2g)$.  Each unitarized Frobenius element is conjugate to some element of $H$ (unique up to conjugacy in $H$) whose characteristic polynomial is $\bar{L}_p(T)$.  When $G_\ell$ is Zariski dense in $GSp(2g,\mathbb{Z}_\ell)$, we obtain $H=USp(2g)$, but in general, $H$ is some compact subgroup of $USp(2g)$.

There is an analogous construction involving the Mumford-Tate group ${\rm MT}(A)$ of an abelian variety $A$, which contains $G_\ell$.  The result is the Hodge group ${\rm Hg}(A)$, corresponding to our subgroup $H$ above.  For simple abelian varieties of low dimension (up to genus 5) the possibilities for ${\rm Hg}(A)$ have been classified by Moonen and Zarhin \cite{Moonen:HodgeClasses}.\footnote{We thank David Zywina for bringing this to our attention.}  For genus 2 curves we  consider the classification of $H$ in Section \ref{section:g2groups}, and give what we believe to be a complete list of the possibilities (most of these do not correspond to simple Jacobians).  We note here that curves with nonisomorphic Hodge groups may have identical $L$-polynomial distributions, so the notions are not equivalent.

We can now state the conjecture.

\begin{conj}[Generalized Sato-Tate]\label{conjecture:GeneralSatoTate}
For a curve $C$ of genus $g$, the distribution of $\barLp$ converges to the distribution of characteristic polynomials $\chi(T)$ in an infinite compact subgroup $H\subseteq USp(2g)$.  Equality holds if and only if $C$ has large Galois image.
\end{conj}

We say that the subgroup $H$ {\em represents} the $L$-polynomial distribution of $C$.

\section{Moment Sequences Attached to $L$-polynomials}

Let $P_C(N)$ denote the set of primes $p\le N$ for which the curve $C$ has good reduction.  We may compute $\barLp$ for all $p\in P_C(N)$, and if $x$ is a quantity derived from $\barLp$ (e.g., the coefficients $a_k$ or some function of the $a_k$), we consider the mean value $m(x^n)$ over $p\in P_C(N)$ as an approximation of the $n$th moment $E[X^n]$ of a corresponding random variable $X$.  Under Conjecture \ref{conjecture:GeneralSatoTate}, we assume there is a compact subgroup $H\subseteq USp(2g)$ which represents the $L$-polynomial distribution of $C$.  If $X$ is a real random variable defined as a polynomial of the eigenvalues of an element of $H$, it clearly has bounded support under the Haar measure on $H$ (the eigenvalues lie on the unit circle).  Therefore its moments all exist.  Further, one may determine absolute bounds on $X$ depending only on $g$, which we regard as fixed.  Carleman's condition \cite[p. 126]{Koosis:LogarithmicIntegral} then implies that the moment sequence for $X$ uniquely determines its distribution.  We summarize this argument with the following proposition.

\begin{proposition}
Under Conjecture \ref{conjecture:GeneralSatoTate}, let $H$ be a compact subgroup of $USp(2g)$ which represents the $L$-polynomial distribution of a curve $C$.  Let $x_p$ be a real-valued polynomial of the roots of $\barLp$, and let the random variable $X$ be the corresponding polynomial of the eigenvalues of a random matrix in $H$.  Then the moments of $X$  exist and determine the distribution of $X$.  For all nonnegative integers $n$, the mean value of $x_p^n$ over $p\in P_C(N)$ converges to $E[X^n]$ as $N\to\infty$.
\end{proposition}

We now consider random variables $X$ that are symmetric functions (polynomials) of the eigenvalues, with integer coefficients.  In this case, it is easy to see that the moments of $X$ must be integers.

\begin{proposition}\label{prop:IntegerMoments}
Let $V$ be a vector space of finite dimension over $\mathbb{C}$.  Let the random variable $X$ be a symmetric polynomial with integer coefficients over the eigenvalues of a random matrix in a compact group $G\subseteq GL(V)$, distributed with Haar measure.  Then $E[X^n]\in\mathbb{Z}$ for all nonnegative integers $n$.
\end{proposition}
\begin{proof}
The sum of the eigenvalues, $e_1$, is the trace of the standard representation $V$ of $G$, and $E[e_1]=\int_G\tr(A)dA$ counts the multiplicity of the trivial representation in $V$, an integer.  Similarly, the $k$th symmetric function $e_k$ is the trace of the $k$th exterior power $\Lambda^kV$, hence $E[e_k]$ must be an integer.  The product of the traces of two representations is the trace of their tensor product, hence $E[e_k^n]$ is an integer, as are the moments of any product of elementary symmetric functions.  Linearity of expectation implies that $\mathbb{Z}$-linear combinations of products of elementary symmetric functions also have integer moments, and every symmetric polynomial with integer coefficients may be expressed in this way.
\end{proof}

Proposition \ref{prop:IntegerMoments} is quite useful in practice, particularly when $H$ is unknown.  We can ``determine" $E[X^n]$ for small $n$ by computing the mean of $x_p^n$ up to a suitable bound $N$.  The value obtained is purely heuristic of course, depending both on Conjecture \ref{conjecture:GeneralSatoTate} and an assumption about the rate of convergence.  Nevertheless, the consistency of the results so obtained are quite compelling.  

For most curves (those with large Galois image) we expect $H=USp(2g)$, and we may compute the distribution of $X$ using the measure $\mu$ defined in (\ref{equation:EigenvalueDistribution}).

\subsection{Moment Generating Functions in $USp(2g)$}

Let $\chi(T)=\sum a_kT^k$ be the characteristic polynomial of a random matrix in $USp(2g)$ with eigenvalues $e^{\pm i\theta_1},\ldots,e^{\pm i\theta_g}$.  Up to a sign $(-1)^k$, the $a_k$ are the elementary symmetric functions of the eigenvalues.  We also define
\begin{equation*}
s_k=\sum_j(e^{ki\theta_j}+e^{-ki\theta_j}) = \sum_j 2\cos k\theta_j,
\end{equation*}
the $k$th power sums of the eigenvalues (Newton symmetric functions).

We wish to compute the (integer) sequence $M[X] = (1, E[X], E[X^2], \ldots)$, where $X$ is $a_k$ or $s_k$.
We first consider $s_k$ (including $a_1=-s_1$).\footnote{The other $a_k$ are addressed in Section \ref{section:ExplicitMoments2} (for $g\le 3$).}  We have
\begin{equation}\label{equation:esk1}
M[s_k](n) = E[s_k^n] = \int_V\Bigl(\sum_j 2\cos k\theta_j\Bigr)^n\mu,
\end{equation}
where $V=[0,\pi]^g$ denotes the volume of integration.  If we expand (\ref{equation:esk1}) using the formula for $\mu$  in (\ref{equation:EigenvalueDistribution}), we need only consider univariate integrals of the form
\begin{equation}\label{equation:BaseIntegral}
C_k^m(n) = \frac{1}{\pi}\int_0^\pi(2\cos k\theta)^n(2\cos\theta)^m(2\sin^2\theta)d\theta.
\end{equation}

We regard $C_k^m$ as a sequence indexed by $n\in \mathbb{Z}^+$ (in fact, an integer sequence), and define the exponential generating function (egf) of $C_k^m$ by
\begin{equation*}
\mathcal{C}_k^m(z)=\sum_{n=0}^\infty C_k^m(n)\frac{z^n}{n!}.
\end{equation*}
Similarly, $\mathcal{M}[s_k]$, the egf of $M[s_k]$, is the moment generating function of $s_k$.
We will compute $C_k^m$ in terms of sequences $B_\nu$, defined by\footnote{$B_\nu(n)$ counts paths of length $n$ from 0 to $\nu$ on the real line using step set $\{\pm 1\}$.} 
\begin{equation}\label{definition:Bnu}
B_\nu(n) = \binom{n}{\frac{n+\nu}{2}},
\end{equation}
where the binomial coefficient is zero when $(n+\nu)/2$ is not an integer in the interval $[0,n]$.  With this understanding, we allow $\nu$ to take arbitrary values, with $B_\nu$ identically zero for $\nu\notin\mathbb{Z}$.
Letting $\mathcal{B}_\nu(z) =\sum_{n=0}^\infty B_\nu(n)\frac{z^n}{n!}$, we find that
\begin{equation}\label{equation:BesselDef}
\mathcal{B}_\nu(z) = \mathcal{I}_\nu(2z) = \sum_{n=0}^\infty \frac{z^{2n+\nu}}{n!(n+\nu)!},\qquad\qquad{\rm for}\enspace\nu\in\mathbb{Z}.
\end{equation}

The function $\mathcal{I}_\nu(z)$ is a hyperbolic Bessel function (of the first kind, of order $\nu$) \cite[Ch. 49-50]{Spanier:AtlasOfFunctions}.  Hyperbolic Bessel functions are defined for nonintegral values of $\nu$ (and are not identically zero), so the condition $\nu\in\mathbb{Z}$ in (\ref{equation:BesselDef}) should be duly noted.

We can now state a concise formula for $\mathcal{M}[s_k]$.
\begin{theorem}\label{theorem:msk}
Let $s_k$ denote the $k$th power sum of the eigenvalues of a random element of $USp(2g)$.  The moment generating function of $s_k$ is
\begin{equation}\label{equation:msk1}
\mathcal{M}[s_k] = \det_{g\times g}\left(\mathcal{C}_k^{i+j-2}\right),
\end{equation}
where $\mathcal{C}_k^m(z)$ is given by
\begin{equation}\label{equation:msk2}
\mathcal{C}_k^m = \sum_j\binom{m}{j}\bigl(\mathcal{B}_{(2j-m)/k}-\mathcal{B}_{(2j-m+2)/k}\bigr),
\end{equation}
with $\mathcal{B}_\nu$ defined as above.
\end{theorem}
We postpone the proof until Section \ref{subsection:Proof}.

The determinantal formula in Theorem \ref{theorem:msk} contains some redundancy (e.g., when $g=3$ the term $C_k^1C_k^2C_k^3$ appears twice), but we find the simple form of Theorem \ref{theorem:msk} well suited to both hand and machine computation.
For example, when $g=2$ and $k=1$ we have
$$\mathcal{M}[s_1] = \mathcal{C}_1^0\mathcal{C}_1^2-\mathcal{C}_1^1\mathcal{C}_1^1 = (\mathcal{B}_0-\mathcal{B}_2)(\mathcal{B}_0-\mathcal{B}_4) - (\mathcal{B}_1-\mathcal{B}_3)^2.$$
{}From (\ref{equation:Ms1genus2}) of Section \ref{section:ExplicitMoments1}, we obtain the identity\footnote{The similarity of $\mathcal{C}_1^0\mathcal{C}_1^2-\mathcal{C}_1^1\mathcal{C}_1^1$ and $c(n)c(n+2)-c(n+1)^2$ is somewhat misleading, both expressions involve Catalan numbers, but the terms do not correspond.}
\begin{equation}\label{equation:g2s1moments}
M[s_1](2n) = c(n)c(n+2)-c(n+1)^2,
\end{equation}
where $c(n)$ is the $n$th Catalan number.  The odd moments are zero and the even moments form sequence \Aseq{005700} in the On-line Encyclopedia of Integer Sequences (OEIS) \cite{Sloane:OEIS}.  A complete list of the sequences $M[s_k]$ and $g\le 3$ can be found in Section \ref{section:ExplicitMoments1}.

The sequence \Aseq{005700} = (1, 1, 3, 14, 84, 594, \ldots) is well known.  It counts lattice paths of length $2n$ in $\mathbb{Z}^2$ with step set $\{(\pm 1,0), (0,\pm 1)\}$ that return to the origin and are constrained by $x_1\ge x_2\ge 0$.  In general, the sequence $M[s_1]$ counts returning lattice paths in $\mathbb{Z}^g$ which remain in the region $x_1\ge\ldots\ge x_g\ge 0$.  This follows from a general result of Grabiner and Magyar\footnote{As explained to us by Arun Ram, this equality may be interpreted in terms of crystal bases, which leads directly to analogues for groups other than $USp(2g)$.} which relates the decomposition of tensor powers of certain Lie group representations to lattice paths in a chamber of the associated Weyl group \cite{Grabiner:WeylChamber}: as in the proof of Proposition (\ref{prop:IntegerMoments}), interpret $E[s_1^n]$ as the multiplicity of the trivial representation in $V^{\otimes n}$, where $V$ is the standard representation of $USp(2g)$, then apply Theorem 2 of \cite{Grabiner:WeylChamber}.

For the group $USp(2g)$ and $k=1$, our results intersect those of \cite{Grabiner:WeylChamber}, where an equivalent determinantal formula is obtained by counting lattice paths in the Weyl chamber of the corresponding Lie algebra.  By contrast, we compute $\mathcal{M}[s_k]$ directly from the measure $\mu(USp(2g))$, using only elementary methods.  The Haar measure, via the Weyl integration formula, effectively encodes the relevant combinatorial content.  Particularly when $k>1$, this is  simpler than a combinatorial or representation-theoretic approach.  More generally, Haar measure, and moment sequences in particular, can provide convenient access to the combinatorial structure of compact groups.

Determinantal formulas of the type above arise in many combinatorial questions related to lattice paths and Young tableaux (see \cite{Gessel:SymmetricFunctions,Chen:CrossingsNestings} for examples).  One might ask what the moment sequences for $s_k$ count when $k>1$.  For $g=2$, some answers may be found in the OEIS (see links in Table \ref{table:MomentSequences1}).

Before proving Theorem \ref{theorem:msk}, we note the following corollary.

\begin{corollary}\label{corollary:randomidd}
For all $k> 2g$, the random variables $s_k$ are identically distributed with moment generating function
$\mathcal{M}[s_k] = (\mathcal{B}_0)^g$.
\end{corollary}
\begin{proof}
From (\ref{equation:msk1}) of Theorem \ref{theorem:msk}, $\mathcal{M}[s_k]$ is a polynomial expression in $\mathcal{C}_k^m$ with $m\le 2g-2$.  From (\ref{equation:msk2}), it follows that $\mathcal{C}_k^m = \binom{m}{m/2}\mathcal{B}_0-\binom{m}{m/2-1}\mathcal{B}_0$ for all $k>m+2$, since $B_\nu$ is zero for nonintegral $\nu$.  Thus $\mathcal{M}[s_k]$ is an integer multiple of $\mathcal{B}_0^g$, and $\mathcal{M}[s_k]$ must be equal to $\mathcal{B}_0^g$, since $M[s_k](0) = 1 = B_0(0)$.\footnote{In fact, for $k>2g$, factoring out $\mathcal{B}_0$ from $\mathcal{M}[s_k]$ leaves the Hankel determinant of the sequence $c(0)$, 0, $c(1)$, 0, $c(2)$, \ldots, which is 1 (see \cite{Aigner:Catalan}).}
\end{proof}

The distribution of $s_k$ given by Corollary \ref{corollary:randomidd} corresponds to the trace of a random matrix under a uniform distribution of $\theta_1,\ldots,\theta_g$.
This is a special case of a general phenomenon first noticed by Eric Rains, who has proven similar results for all the compact classical groups \cite{Rains:Powers,Rains:EigenvalueDistribution}.  The sudden transition to a fixed distribution is quite startling when first encountered; one might na\"{i}vely expect the distribution of $s_k$ to gradually converge as $k\to\infty$.
There is, however, an elementary explanation (see the proof of Lemma \ref{lemma:Ckm} below).

\subsection{Proof of Theorem \ref{theorem:msk}}\label{subsection:Proof}

We rearrange the integral for $M[s_k](n)=E[s_k^n]$ to obtain a determinantal expression in $C_k^m$.  Lemma \ref{lemma:Ckm} then evaluates $C_k^m(n)$.

\begin{proof}[Proof of Theorem $\ref{theorem:msk}$]
Let $w_j=2\cos k\theta_j$ for $1\le j\le g$.  Then $s_k = \sum w_j$ and the integral for $M[s_k](n)$ in (\ref{equation:esk1}) becomes
\begin{equation}\label{equation:pmsk1}
M[s_k](n) = \int_V\Bigl(\sum_j w_j\Bigr)^n\mu = \sum_vn_v\int_V\prod_jw_j^{v_j}\mu
\end{equation}
where $V=[0,\pi]^g$ is the volume of integration, $\mu$ is the Haar measure on $USp(2g)$ defined in (\ref{equation:EigenvalueDistribution}), $v$ ranges over vectors of $g$ nonnegative integers, and $n_v=\binom{n}{v_1,\ldots,v_g}$.  Now let $x_j=2\cos\theta_j$ and $y_j=2\sin^2\theta_j$ so that (\ref{equation:EigenvalueDistribution}) becomes
$$\mu = \frac{1}{g!\pi^g}\prod_{i<j}(x_i-x_j)^2\prod_j(y_jd\theta_j).$$
By Lemma \ref{lemma:Discriminant}, we may write this as
\begin{equation}\label{equation:pmsk2}
\mu = \frac{1}{g!\pi^g}\sum_\sigma\det_{g\times g}\left(x_{\sigma(j)}^{i+j-2}\right)\prod_j (y_jd\theta_j),
\end{equation}
where $\sigma$ ranges over permutations of $\{1,\ldots,g\}$.  From the definition of $C_k^m(n)$ in (\ref{equation:BaseIntegral}) we have  $C_k^m(n)=\frac{1}{\pi}\int_0^\pi w_j^nx_j^my_jd\theta_j$, for any $j$. Combining (\ref{equation:pmsk1}) and (\ref{equation:pmsk2}),
\begin{align*}
M[s_k](n) = \sum_v n_v\frac{1}{g!}\sum_\sigma \det_{g\times g}\left(C_k^{i+j-2}(v_{j})\right)
          = \sum_v n_v\det_{g\times g}\left(C_k^{i+j-2}(v_j)\right).
\end{align*}
In terms of egfs (by Lemma \ref{lemma:egfDeterminant}), we then have the desired expression
$$\mathcal{M}[s_k] = \det_{g\times g}\left(\mathcal{C}_k^{i+j-2}\right).$$
Applying Lemma \ref{lemma:Ckm} to evaluate the integral $C_k^m(n)$ completes the proof.
\end{proof}

Lemmas (\ref{lemma:Discriminant}) and (\ref{lemma:egfDeterminant}) are elementary facts; lacking a ready reference, we provide short proofs.

\begin{lemma}\label{lemma:Discriminant}
If $r_1$, \ldots, $r_g$ are indeterminates in a commutative ring, then
$$\prod_{j< k}(r_j-r_k)^2 = \sum_\sigma \det_{g\times g}\left(r_{\sigma(j)}^{i+j-2}\right),$$
where $\sigma$ ranges over permutations of $\{1,\ldots,g\}$.
\end{lemma}
\begin{proof}
Recall that $\prod_{j< k}(r_j-r_k)$ is the Vandermonde determinant \cite[p. 71]{Knapp:BasicAlgebra}.  Define $r^{\alpha}=\prod_jr_j^{\alpha(j)}$ for any function $\alpha:\{1,\ldots,g\}\to\mathbb{Z}^+$, and for a permutation $\sigma$ let $r_\sigma=(r_{\sigma(1)},\ldots,r_{\sigma(g)})$.  Then $\prod_{j< k}(r_j-r_k)^2$ is given by
\begin{align*}
\left(\det\left(r_i^{j-1}\right)\right)^2 &= \left(\sum_\pi\sgn(\pi) r^{\pi-1}\right)^2 = \sum_{\pi,\phi}\sgn(\pi)\sgn(\phi) r^{\pi+\phi-2}\\
		&= \sum_{\pi,\phi}\sgn(\pi\phi)r_{\phi}^{\pi(\phi^{-1})+{\rm id}-2} = \sum_{\pi,\phi}\sgn(\pi\phi^{-1})r_{\phi}^{\pi(\phi^{-1})+{\rm id}-2} 
\\
	   &=\sum_{\sigma}\sum_{\phi}\sgn(\phi) r_\sigma^{\phi+{\rm id}-2} = \sum_{\sigma}\det\left(r_{\sigma(j)}^{i+j-2}\right),
\end{align*}
and the lemma is proven.
\end{proof}

\begin{lemma}\label{lemma:egfDeterminant}
For $1\le i,j\le m$, let $\mathcal{A}_{i,j}\in\mathbb{C}[[z]]$. Then
$$[n]\det_{m\times m}(\mathcal{A}_{i,j}) = \sum_{n_1,\ldots,n_m}\binom{n}{n_1,\ldots,n_m}\det_{m\times m}\bigl([n_j]\mathcal{A}_{i,j}\bigr),$$
where $[n]\mathcal{F}$ denotes the coefficient of $z^n/n!$ in $\mathcal{F}\in\mathbb{C}[[z]]$.  
\end{lemma}
\begin{proof}
For any $\mathcal{F}_1, \ldots, \mathcal{F}_m\in\mathbb{C}[[z]]$ we have
$$[n]\prod_j\mathcal{F}_j = \sum_{n_1,\ldots,n_m}\binom{n}{n_1,\ldots,n_m}\prod_j[n_j]\mathcal{F}_j.$$
It follows that
\begin{align*}
[n]\det_{m\times m}(\mathcal{A}_{i,j}) &= \sum_\sigma\sgn(\sigma)\sum_i[n]\prod_j\mathcal{A}_{\sigma(i),j}\\
&= \sum_\sigma\sgn(\sigma)\sum_i\sum_{n_1,\ldots,n_m}\binom{n}{n_1,\ldots,n_m}\prod_j[n_j]\mathcal{A}_{\sigma(i),j}\\
&= \sum_{n_1,\ldots,n_m}\binom{n}{n_1,\ldots,n_m}\det_{m\times m}\bigl([n_j]\mathcal{A}_{i,j}\bigr),
\end{align*}
as desired.
\end{proof}

Lemma \ref{lemma:Ckm} computes the integral $C_k^m$.

\begin{lemma}\label{lemma:Ckm}
Define $C_k^m(n) = \frac{1}{\pi}\int_0^\pi(2\cos k\theta)^n(2\cos\theta)^m(2\sin^2\theta)d\theta$, for positive integers $k$ and nonnegative integers $m$ and $n$, and let $B_\nu(n) = \binom{n}{\frac{n+\nu}{2}}$.  Then
$$C_k^m(n) = \sum_j\binom{m}{j}\bigl(B_{(2j-m)/k}(n)-B_{(2j-m+2)/k}(n)\bigr)$$
for all nonnegative integers $n$.
\end{lemma}
\begin{proof}
We may write $C_k^m(n)$ as
$$\frac{1}{\pi}\int_0^\pi\left(e^{ik\theta}+e^{-ik\theta}\right)^n \left(e^{i\theta}+e^{-i\theta}\right)^m \left(1-\left(e^{2i\theta}+e^{-2i\theta}\right)/2\right) d\theta.$$
In terms of $\delta(t) = \frac{1}{\pi}\int_0^\pi e^{it\theta} d\theta$, one finds that $C_k^m(n)$ is equal to
\begin{multline*}
\sum_r\sum_j\binom{n}{r}\binom{m}{j}\delta(k(2r-n)+2j-m)\\
- \frac{1}{2}\sum_r\sum_j\binom{n}{r}\binom{m}{j}\delta(k(2n-r)+2j-m+2)\\
- \frac{1}{2}\sum_r\sum_j\binom{n}{r}\binom{m}{j}\delta(k(2n-r)+2j-m-2).
\end{multline*}
As $C_k^m(n)$ is a real number, we need only consider the real parts of the sums above.
For real $t$, the real part of $\delta(t)$ is nonzero only when $t=0$, in which case it is 1.
Hence, in Iverson notation, $\Re(\delta(t)) = [t=0]$ holds for all $t\in\mathbb{R}$.\footnote{The function $[P]$ is 1 when the boolean predicate $P$ is true and 0 otherwise, see \cite{Knuth:Notation}.} 

The real parts of the second two sums are equal, since we may replace $j$ with $m-j$ and $r$ with $n-r$ in the last sum and then apply $[t=0]=[-t=0]$.  Interchanging the order of summation we obtain
$$\sum_j\sum_r\binom{m}{j}\binom{n}{r}\Bigl([k(2r-n)+2j-m=0] - [k(2n-r)+2j-m+2=0]\Bigr).$$
We now note that
$$\sum_r\binom{n}{r}[2r-n=\nu] = \binom{n}{\frac{n+\nu}{2}} = B_\nu(n),$$
for all nonnegative integers $n$ and arbitrary $\nu$.  Thus we have
$$C_k^m(n) = \sum_j\binom{m}{j}\bigl(B_{(m-2j)/k}(n) - B_{(m-2j-2)/k}(n)\bigr).$$
Applying the identity $B_\nu = B_{-\nu}$ completes the proof.
\end{proof}

\subsection{Explicit Computation of $M[s_k]$ in $USp(2g)$, for $g\le 3$.}\label{section:ExplicitMoments1}

For $g\le 3$, Theorem \ref{theorem:msk} gives:
\begin{tabbing}
\hspace{30pt}\=$g$=1:\hspace{10pt}\=$M[s_k] = C_k^0$;\\
\>$g$=2:\>$M[s_k] = C_k^0C_k^2-C_k^1C_k^1$;\\
\>$g$=3:\>$M[s_k] = C_k^0C_k^2C_k^4 + 2C_k^1C_k^2C_k^3 - C_k^0C_k^3C_k^3 - C_k^1C_k^1C_k^4 - C_k^2C_k^2C_k^2$.
\end{tabbing}
We will compute these moment sequences explicitly.

\pagebreak
\noindent
For convenience, define
\begin{onehalfspace}
\begin{tabbing}
\hspace{30pt}\=$A = B_0-B_1$\hspace{30pt}\=\Aseq{126930}=(1, -1, 2, -3, 6, -10, 20, -35, 70, -126, \dots),\\
\>$B = B_0$ \>\Aseq{126869}=(1, 0, 2, 0, 6, 0, 20, 0, 70, 0, 252, \dots),\\
\>$C = B_0-B_2$ \>\Aseq{126120}=(1, 0, 1, 0, 2, 0, 5, 0, 14, 0, 42, \dots),\\
\>$D = B_1$ \>\Aseq{138364}=(0, 1, 0, 3, 0, 10, 0, 35, 0, 126, 0, \dots),
\end{tabbing}
\end{onehalfspace}
\noindent and let $\mathcal{A},\mathcal{B},\mathcal{C}$, and $\mathcal{D}$ denote the corresponding egfs.
\begin{lemma}\label{lemma:CkmIdentities}
Let $\Dop$ denote the derivative operator.
\begin{onehalfspace}
\begin{tabbing}
\hspace{30pt}\={\rm 1}. $\mathcal{C}_1^m = \Dop^m\mathcal{C}.$\hspace{100pt}\=\\
\>{\rm 2}. $\mathcal{C}_2^m = (\Dop+2)^{(m-2)/2}\mathcal{C},$\>for even $m>0$.\\
\>{\rm 3}. $\mathcal{C}_k^m=0,$\>for $m$ odd and $k$ even.\\
\>{\rm 4}. $\mathcal{C}_k^m=C(m)\mathcal{B},$\>for $k = m + 1 > 1$, or $k > m + 2$.
\end{tabbing}
\end{onehalfspace}
\end{lemma}
\vspace{-24pt}
\begin{proof}
Recall from (\ref{equation:msk2}) of Theorem \ref{theorem:msk} that
$$\mathcal{C}_k^m = \sum_j\binom{m}{j}\bigl(\mathcal{B}_{(2j-m)/k}-\mathcal{B}_{(2j-m+2)/k}\bigr).$$
By Pascal's identity, we have
\begin{equation}\label{equation:BPascal}
\Dop\mathcal{B}_\nu = \mathcal{B}_{\nu+1} + \mathcal{B}_{\nu-1}.
\end{equation}
The proofs of statements 1 and 2 are then straightforward inductions on $m$.  Statements 3 and 4 follow immediately from Lemma \ref{lemma:Ckm}.
\end{proof}
Applying Lemmas \ref{lemma:Ckm} and \ref{lemma:CkmIdentities}, we compute Table \ref{table:Ckm}. From Table \ref{table:Ckm} and (\ref{equation:msk1}) of Theorem \ref{theorem:msk} we obtain a closed form for $\mathcal{M}[s_k]$ in terms of $\mathcal{A},\mathcal{B},\mathcal{C}$, and $\mathcal{D}$.  To determine $M[s_k](n)$ we must compute linear combinations of multinomial convolutions of various $B_j$. 
This is reasonably efficient for small values of $g$ and $n$, however we can speed up the process significantly with the following lemma.

\setlength\extrarowheight{2pt}
\begin{table}
\begin{center}
\begin{tabular}{ | c | c | c | c | c | c | c | c |}
\hline
$k\backslash m$ & $\quad 0\quad$ & $\quad 1\quad$ & $\quad 2\quad$ & $\quad 3\quad$ & $\quad 4\quad$ & $\quad 5\quad$ & $\quad 6\quad$\\\hline
1 & $\mathcal{C}$ & $\Dop \mathcal{C}$ & $\Dop^2 \mathcal{C}$ & $\Dop^3 \mathcal{C}$ & $\Dop^4 \mathcal{C}$ & $\Dop^5 \mathcal{C}$ & $\Dop^6 \mathcal{C}$\\\hline
2 & $\mathcal{A}$ & 0 & $\mathcal{C}$ & 0 & $(\Dop+2)\mathcal{C}$ & 0 & $(\Dop+2)^2\mathcal{C}$\\\hline
3 & $\mathcal{B}$ & $-\mathcal{D}$ & $\mathcal{B}$ & $-\mathcal{D}$ & $\mathcal{B}+\mathcal{C}$ & $-\mathcal{D}$ & $\mathcal{B}+4\mathcal{C}$\\\hline
4 & $\mathcal{B}$ & 0 & $\mathcal{A}$ & 0 & $2\mathcal{A}$ & 0 & $4\mathcal{A}+\mathcal{C}$\\\hline
5 & $\mathcal{B}$ & 0 & $\mathcal{B}$ & $-\mathcal{D}$ & $2\mathcal{B}$ & $-3\mathcal{D}$ & $5\mathcal{B}$\\\hline
6 & $\mathcal{B}$ & 0 & $\mathcal{B}$ & 0 & $\mathcal{A}+\mathcal{B}$ & 0 & $4\mathcal{A}+\mathcal{B}$\\\hline
7 & $\mathcal{B}$ & 0 & $\mathcal{B}$ & 0 & $2\mathcal{B}$ & $-\mathcal{D}$ & $5\mathcal{B}$\\\hline
8 & $\mathcal{B}$ & 0 & $\mathcal{B}$ & 0 & $2\mathcal{B}$ & 0 & $5\mathcal{B}-\mathcal{D}$\\\hline
9 & $\mathcal{B}$ & 0 & $\mathcal{B}$ & 0 & $2\mathcal{B}$ & 0 & $5\mathcal{B}$\\\hline
\end{tabular}
\vspace{12pt}
\caption{Exponential Generating Functions $\mathcal{C}_k^m$.}\label{table:Ckm}
\end{center}
\end{table}

\begin{lemma}\label{lemma:Z2paths}
Let $\mathcal{B}_\nu(z)=\sum_{n=0}^\infty B_\nu(n)z^n/n!$, where $B_\nu(n) = \binom{n}{\frac{n+\nu}{2}}$.  Then
$$\mathcal{B}_a(z)\mathcal{B}_b(z) = \sum_{n=0}^\infty B_{a+b}(n)B_{a-b}(n)\frac{z^n}{n!}$$
for all $a,b\in\mathbb{Z}$.
\end{lemma}
\begin{proof}
The coefficient of $z^n/n!$ on both sides of the equality count lattice paths from $(0,0)$ to $(a,b)$ in $\mathbb{Z}^2$ with step set $\{(\pm 1,0),(0,\pm 1)\}$.  This is immediate for the LHS.  A simple bijective proof for the RHS appears in \cite{Guy:LatticePaths}.
\end{proof}

In genus 2, Lemma \ref{lemma:Z2paths} gives us a closed form for $\mathcal{M}[s_k](n)$ in terms of binomial coefficients (or Catalan numbers).  For example, from
$$\mathcal{M}[s_1] = \mathcal{C}_1^0\mathcal{C}_1^2-\mathcal{C}_1^1\mathcal{C}_1^1 = (\mathcal{B}_0-\mathcal{B}_2)(\mathcal{B}_0-\mathcal{B}_4)-(\mathcal{B}_1-\mathcal{B}_3)^2,$$
we obtain
\begin{align*}
M[s_1] &= B_0^2 - B_4^2 - B_2^2+B_2B_6-B_0B_2+2B_2B_4-B_0B_6\\
		&= (B_0-B_2)(2(B_0-B_4)+B_2-B_6) - (B_0-B_4)^2.
\end{align*}
This may be expressed more compactly as
\begin{equation}\label{equation:Ms1genus2}
M[s_1](n) = C(n)C(n+2)-C(n+4).
\end{equation}
Here $C(2n)=c(n)$ is the $n$th Catalan number, giving the identity (\ref{equation:g2s1moments}) noted earlier.  Similar formulas for the other $M[s_k]$ in genus 2 are listed in Table \ref{table:MomentSequences1}.  In higher genera we do not obtain a closed form, but computation is considerably faster with Lemma \ref{lemma:Z2paths}; in genus 4 we use $O(n)$ multiplications to compute $M[s_k](n)$, rather than $O(n^3)$.

By Corollary \ref{corollary:randomidd}, the sequences for $k > 2g$ are all the same, so it suffices to consider $k\le 2g+1$.  For even $k\le 2g$ we find that $E[s_k]=-1$, hence we also consider $M[s_k+1]=M[s_k^+]$, the sequence of central moments.  These may be obtained by computing the binomial convolution of $M[s_k]$ with the sequence $(1, 1, 1, \ldots)$.  In genus 1 we obtain

\vspace{-9pt}
\begin{doublespace}
\begin{tabbing}
\hspace{12pt}\=$M[s_2^+]$ = (1, 0, 1, 1, 3, 6, 15, 36, 91, 232, 603, 1585, \ldots)\hspace{24pt}\=\Aseq{005043},\\
and in genus 2 we find\\
\>$M[s_2^+]$ = (1, 0, 2, 1, 11, 16, 95, 232, 1085, 3460, 14820 \ldots)\>\Aseq{138351},\\
\>$M[s_4^+]$ = (1, 0, 3, 1, 21, 26, 215, 493, 2821, 9040, 43695, \ldots)\>\Aseq{138354}.
\end{tabbing}
\end{doublespace}
\setlength\extrarowheight{1pt}
\begin{table}
\begin{center}
\begin{tabular}{ccll}
\hline
$\space\space g\space\space$ & $\space\space k\space\space$ & $M[s_k](n) = E[s_k^n]$& OEIS\\\hline
1 & 1 & $C(n)$&\\
  &   & 1, 0, 1, 0, 2, 0, 5, 0, 14, 0, 42, 0, 132, 0, 429, \ldots&\Aseq{126120}\\
  & 2 & $A(n)$&\\
  &   & 1, -1, 2, -3, 6, -10, 20, -35, 70, -126, 252, -462, 924, \ldots&\Aseq{126930}\\
  & 3 & $B(n)$&\\
  &   & 1, 0, 2, 0, 6, 0, 20, 0, 70, 0, 252, 0, 924, 0, 3432, \ldots&\Aseq{126869}\\\hline 
2 & 1 & $C(n)C(n+4) - C(n+2)^2$&\\
  &   & 1, 0, 1, 0, 3, 0, 14, 0, 84, 0, 594, 0, 4719, 0, 40898, \ldots&\Aseq{138349}\\
  & 2 & $C(n)D(n+1) - D(n)C(n+1)$&\\
  &   & 1, -1, 3, -6, 20, -50, 175, -490, 1764, -5292, 19404, \ldots&\Aseq{138350}\\
  & 3 & $B(n)C(n)$&\\
  &   & 1, 0, 2, 0, 12, 0, 100, 0, 980, 0, 10584, 0, 121968, \ldots&\Aseq{000888}*\\
  & 4 & $B(n)^2-D(n)^2$&\\
  &   & 1, -1, 4, -9, 36, -100, 400, -1225, 4900, -15876, 63504, \ldots&\Aseq{018224}*\\
  & 5 & $B(n)^2$\\
  &   & 1, 0, 4, 0, 36, 0, 400, 0, 4900, 0, 63504, 0, 853776, \ldots&\Aseq{002894}*\\\hline
\end{tabular}
\vspace{12pt}
\caption{Moment Sequences $M[s_k]$ for $g\le 2$.}\label{table:MomentSequences1}
\setlength\extrarowheight{4pt}
\vspace{30pt}
\begin{tabular}{clr}
\hline
$\space\space X\space\space$ & $M[X](n) = E[X^n]$& OEIS\\\hline
$s_1$ & 1, 0, 1, 0, 3, 0, 15, 0, 104, 0, 909, 0, 9449, 0, 112398, \ldots&\Aseq{138540}\\
$s_2$ & 1, -1, 3, -7, 24, -75, 285, -1036, 4242, -16926, 73206, \ldots&\Aseq{138541}\\
$s_2^+$ & 1, 0, 2, 0, 11, 1, 95, 36, 1099, 982, 15792, 25070,\ldots&\Aseq{138542}\\
$s_3$ & 1, 0, 3, 0, 26, 0, 345, 0, 5754, 0, 110586, 0, 2341548 \ldots&\Aseq{138543}\\
$s_4$ & 1, -1, 4, -9, 42, -130, 660, -2415, 12810, -51786, 281736, \ldots&\Aseq{138544}\\
$s_4^+$ & 1, 0, 3, 1, 27, 26, 385, 708, 7231, 20296, 164277, \ldots&\Aseq{138545}\\
$s_5$ & 1, 0, 4, 0, 42, 0, 660, 0, 12810, 0, 281736, 0, 6727644, \ldots&\Aseq{138546}\\
$s_6$ & 1, -1, 6, -15, 90, -310, 1860, -7455, 44730, -195426,  \ldots&\Aseq{138547}\\
$s_6^+$ & 1, 0, 5, 1, 63, 46, 1135, 1800, 25431, 66232, 666387, \ldots&\Aseq{138548}\\
$s_7$ & 1, 0, 6, 0, 90, 0, 1860, 0, 44730, 0, 1172556, \ldots&\Aseq{002896}\\\hline
\end{tabular}
\vspace{12pt}
\caption{Moment Sequences $M[s_k]$ and $M[s_{2k}^+]$ for $g=3$.$^\dagger$}\label{table:MomentSequences2}
\end{center}
\small
\noindent
\raggedright
*The OEIS sequence differs slightly.\\
$^\dagger$ The notation $X^+$ denotes the random variable $X+1$.
\normalsize
\end{table}

\subsection{Explicit Computation of $M[a_k]$ in $USp(2g)$ for $g\le 3$}\label{section:ExplicitMoments2}

To complete our study of moment sequences, we now consider the coefficients $a_k$ of the characteristic polynomial $\chi(T)$ of a random matrix in $USp(2g)$.  We have already addressed $a_1=-s_1$.  For $k>1$, the Newton identities allow us to express $a_k$ in terms of $s_1$, \ldots, $s_k$, however this does not allow us to easily compute $M[a_k]$ from the sequences $M[s_j]$ (the covariance among the $s_j$ is nonzero).  Instead, we note that by writing
\begin{equation}\label{equation:charpoly}
\chi(T) = \prod_j\Bigl((T-e^{i\theta_j})(T-e^{-i\theta_j})\Bigr) = \prod_j(T^2-2\cos\theta_jT+1),
\end{equation}
each $a_k$ may be expressed as a polynomial in $2\cos\theta_1,\ldots,2\cos\theta_g$.  It follows that we may compute
$M[a_k]$ in terms of the sequences $C_1^m$ already considered. The following proposition addresses the cases that arise for $g\le 3$.
\begin{proposition}\label{prop:Mak}
Let $C(n) = B_0(n)-B_2(n)$ as above and let $a_k$ be the coefficient of $T^k$ in $\chi(T)$, the characteristic polynomial of a random matrix in $USp(2g)$.  For $g=2$ we have
\begin{equation}\label{equation:Ma2g2}
E[(a_2-2)^n] = C(n)C(n+2)-C(n+1)^2,
\end{equation}
and if $g=3$ then
\begin{equation}\label{equation:Ma2g3}
E[(a_2-3)^n] = \sum_{n_1,n_2,n_3}\binom{n}{n_1, n_2, n_3}\det_{3\times 3}C(n-n_j+i+j-2).
\end{equation}
Also for $g=3$ we have
\begin{equation}\label{equation:Ma3g3}
E[a_3^n] = \sum_{n_1,n_2,n_3,m}\binom{n}{n_1, n_2, n_3, m}2^{n-m}\det_{3\times 3}C(m+n_j+i+j-2).
\end{equation}
\end{proposition}
\begin{proof}
In genus 2, equation (\ref{equation:charpoly}) gives $a_2-2 = (2\cos\theta_1)(2\cos\theta_2)$ and we have
\begin{equation}\label{equation:Ea2_1}
E[(a_2-2)^n] = \int_0^\pi\int_0^\pi(2\cos\theta_1)^n(2\cos\theta_2)^n\mu.
\end{equation}
By Lemma \ref{lemma:Ckm}, we note that
\begin{equation}\label{equation:Ea2_2}
C(n+m) = C_1^m(n) = \frac{1}{\pi}\int_0^\pi(2\cos\theta)^n(2\cos\theta)^m(2\sin^2\theta)d\theta.
\end{equation}
Expanding (\ref{equation:Ea2_1}) and applying (\ref{equation:Ea2_2}) yields (\ref{equation:Ma2g2}).  In genus 3 we write
$$a_2-3 = (2\cos\theta_1)(2\cos\theta_2)+(2\cos\theta_1)(2\cos\theta_3)+(2\cos\theta_2)(2\cos\theta_3),$$
and apply Lemma \ref{lemma:Discriminant} to write the expanded integral for $E[(a_2-3)^n]$ in determinantal form to obtain (\ref{equation:Ma2g3}) (we omit the details).  For (\ref{equation:Ma3g3}), note that
$$a_3 = (2\cos\theta_1)(2\cos\theta_2)(2\cos\theta_3) + 2(2\cos\theta_1+2\cos\theta_2+2\cos\theta_3),$$
and proceed similarly.
\end{proof}

Taking the binomial convolution of $M[a_2-g]$ with the sequence $(1,g,g^2,\ldots)$ gives the moment sequence for $a_2$ in genus $g$.  One finds that $E[a_2] = 1$, hence we also consider the sequence of central moments, $M[a_2-1]$.
Table \ref{table:MomentSequences3} gives the complete set of moment sequences for $a_k$ in genus $g\le 3$, including $a_1=-s_1$.

\setlength\extrarowheight{2pt}
\begin{table}
\begin{center}
\begin{tabular}{lclr}
\hline
$g$&$\quad X\quad$ & $M[X](n) = E[X^n]$& OEIS\\\hline
1 & $a_1$ & 1, 0, 1, 0, 2, 0, 5, 0, 14, 0, 42, 0, 132, \ldots&\Aseq{126120}\\\hline
2 & $a_1$ & 1, 0, 1, 0, 3, 0, 14, 0, 84, 0, 594, 0, 4719, \ldots&\Aseq{138349}\\
  & $a_2$ & 1, 1, 2, 4, 10, 27, 82, 268, 940, 3476, 13448, \ldots&\Aseq{138356}\\
  & $a_2^-$ & 1, 0, 1, 0, 3, 1, 15, 15, 105, 190, 945, 2410,\ldots&\Aseq{095922}\\\hline
3 & $a_1$ & 1, 0, 1, 0, 3, 0, 15, 0, 104, 0, 909, 0, 9449, \ldots&\Aseq{138540}\\
  & $a_2$ & 1, 1, 2, 5, 16, 62, 282, 1459, 8375, 52323 \ldots&\Aseq{138549}\\
  & $a_2^-$ & 1, 0, 1, 1, 5, 16, 75, 366, 2016, 11936, 75678, \ldots&\Aseq{138550}\\
  & $a_3$ & 1, 0, 2, 0, 23, 0, 684, 0, 34760, 0, 2493096, \ldots&\Aseq{138551}\\\hline
\end{tabular}
\vspace{12pt}
\caption{Moment Sequences $M[a_k]$ and $M[a_{2k}^-]$ for $g\le 3$.$^\dagger$}\label{table:MomentSequences3}
\end{center}
\raggedright
$^\dagger$The notation $X^-$ denotes the random variable $X-1$.
\end{table}

\section{Moment Statistics of Hyperelliptic Curves}

Having computed moment sequences attached to characteristic polynomials of random matrices in $USp(2g)$, we now consider the corresponding moment statistics of a hyperelliptic curve.  Under Conjecture \ref{conjecture:GeneralSatoTate}, the latter should converge to the former, provided the curve has large Galois image.  Tables \ref{table:g1stats}-\ref{table:g3stats} list moment statistics of a hyperelliptic curves known to have large Galois image.  

These tables were constructed by computing $\barLp$ to determine sample values of $a_k$ and $s_k$ for each $p$ (the $a_k$ are the coefficients, the $s_k$ are derived via the Newton identities).  Central moment statistics of $X=a_k-E[a_k]$ and $X=s_k-E[s_k]$ were then computed by averaging $X^n$ over all $p\le N$ where the curve has good reduction.  The values of $E[a_k]$ and $E[s_k]$ are as determined in the previous two sections, 0 for $k$ odd or $k>2g$ and $\pm 1$ otherwise.

Using central moments, we find the moment statistic $M_1$ (the mean of $x$) very close to zero ($|M_1|<0.001$), so we list only $M_2$, \ldots, $M_{10}$ for each $x$.  Beneath each row the corresponding moments for $USp(2g)$ are listed.  Note that the value of $N$ is not the same in each table (we are able to use larger $N$ in lower genus), and the number of sample points is approximately $\pi(N)\approx N/\log{N}$.

Tables \ref{table:g2seqconverge} and \ref{table:g3seqconverge} show the progression of the moment statistics in genus 2 and 3 as $N$ increases, giving a rough indication of the rate of convergence and the degree of uncertainty in the higher moments.

The agreement between the moment statistics listed in Tables \ref{table:g1stats}-\ref{table:g3stats} and the moment sequences computed in Sections \ref{section:ExplicitMoments1}-\ref{section:ExplicitMoments2} is consistent with Conjecture \ref{conjecture:GeneralSatoTate}.  Indeed, on the basis of these results we can quite confidently reject certain alternative hypotheses.

As an example, consider the fourth moment of $a_1$ in genus 2.  The value $M_4=3.004$ represents the average of 3,957,807 data points ($\approx\pi(2^{26})$).  A uniform distribution on $a_1$ would imply the mean value of $a_1^4$ is greater than 50.  The probability of then observing $M_4=3.004$ over a sample of nearly four million data points is astronomically small.  A uniform distribution on the eigenvalue angles gives a mean of 36 yielding a similarly improbable event.  In fact, let us suppose only that $a_1^4$ has a distribution with integer mean not equal to 3.  We can then bound the standard deviation by 256 (since $a_1^8\le 256^2$) and apply a $Z$-test to obtain an event probability less than one in a trillion. 

\pagebreak
\setlength\extrarowheight{0pt}
\begin{table}
\begin{center}
\begin{tabular}{lrrrrrrrrr}
\toprule
$X$ & \hspace{10pt}$M_2$ & \hspace{10pt}$M_3$ & \hspace{10pt}$M_4$ & \hspace{10pt}$M_5$ & \hspace{10pt}$M_6$  & \hspace{10pt}$M_7$  & \hspace{10pt}$M_8$ & \hspace{10pt}$M_9$ & \hspace{10pt}$M_{10}$\\\midrule
$a_1$   & 1.000 & 0.000 & 2.000 & 0.000 &  5.000 &  0.001 &  14.000 &   0.002 &  42.000\\\vspace{8pt}
        & 1 & 0 & 2 & 0 & 5 & 0 & 14 & 0 & 42\\
$s_2^+$ & 1.000 & 1.000 & 3.000 & 6.000 & 15.001 & 36.003 &  91.010 & 232.03 & 603.11\\\vspace{8pt}
        & 1 & 1 & 3 & 6 & 15 & 36 & 91 & 232 & 603\\
$s_3$   & 2.000 & 0.000 & 6.000 & 0.000 & 20.000 &  0.000 &  69.998 &   0.000 & 252.00\\\vspace{8pt}
        & 2 & 0 & 6 & 0 & 20 & 0 & 70 & 0 & 252\\
$s_4$   & 2.000 & 0.000 & 6.000 & 0.000 & 20.000 &  0.000 &  70.001 &  -0.001 & 252.00\\
        & 2 & 0 & 6 & 0 & 20 & 0 & 70 & 0 & 252\\
\bottomrule
\end{tabular}
\vspace{12pt}
\caption{Central moment statistics of $a_k$ and $s_k$ in genus 1, $N=2^{35}$.}\label{table:g1stats}
\vspace{-20pt}
\small
$$y^2 = x^3 + 314159x + 271828$$.\\\vspace{24pt}
\normalsize
\begin{tabular}{lrrrrrrrrr}
\toprule
$X$ & $M_2$ & $M_3$ & $M_4$ & $M_5$ & $M_6$  & $M_7$  & $M_8$ & $M_9$ & $M_{10}$\\\midrule
$a_1$   & 1.001 & -0.001 & 3.004 & -0.006 & 14.014 & -0.031 & 84.041 & -0.178 & 594.02\\\vspace{8pt}
        & 1 & 0 & 3 & 0 & 14 & 0 & 84 & 0 & 594\\
$a_2^-$ &  1.001 &  0.000  & 3.003  & 0.997  & 15.013  & 14.964  & 105.10  & 190.00  & 947.38 \\\vspace{8pt}
        & 1 & 0 & 3 & 1 & 15 & 15 & 105 & 190 & 945\\
$s_2^+$ & 2.001 &  1.002 &  11.014  & 16.044 &  95.247 &  232.90 &  1089.3 &  3476.4 &  14891\\\vspace{8pt}
        & 2 & 1 & 11 & 16 & 95 & 232 & 1085 & 3460 & 14820\\
$s_3$   &  2.002 &  0.001 &  12.014 &  -0.001 &  100.14 &  -0.147 &  981.54 &  -2.850 &  10603\\\vspace{8pt}
        & 2 & 0 & 12 & 0 & 100 & 0 & 980 & 0 & 10584\\
$s_4^+$ &  3.004 &  1.010 &  21.049  & 26.150 &  215.66 &  500.32 &  2830.6 &  9075.6 &  43836\\\vspace{8pt}
        & 3 & 1 & 21 & 26 & 215 & 493 & 2821 & 9040 & 43695\\
$s_5$   &  3.996 &  -0.015 &  35.958  & -0.211 &  399.62 &  -3.152 &  4897.2  & -47.602 & 63492\\\vspace{8pt}
        & 4 & 0 & 36 & 0 & 400 & 0 & 4900 & 0 & 63504\\
$s_6$   & 3.999 & -0.002 & 35.983 & -0.023 & 399.81 &-0.490 & 4898.0 & -9.460 & 63487\\
        & 4 & 0 & 36 & 0 & 400 & 0 & 4900 & 0 & 63504\\
\bottomrule
\end{tabular}
\vspace{12pt}
\caption{Central moment statistics of $a_k$ and $s_k$ in genus 2, $N=2^{26}$.}\label{table:g2stats}
\vspace{-20pt}
\small
$$y^2 = x^5 + 314159x^3 + 271828x^2 + 1644934x + 57721566.$$
\normalsize
\end{center}
\end{table}

\begin{table}
\begin{center}
\begin{tabular}{lrrrrrrrrr}
\toprule
$X$& $M_2$ & $M_3$ & $M_4$ & $M_5$ & $M_6$  & $M_7$  & $M_8$ & $M_9$ & $M_{10}$\\\midrule
$a_1$ & 0.999 &  -0.007 &  2.995 &  -0.046 &  14.968 &  -0.343 &  103.76 &  -2.723 &  906.67\\\vspace{10pt}
&1 & 0 & 3 & 0 & 15 & 0 & 104 & 0 & 909\\
$a_2^-$ & 0.999 &  0.996 &  4.992  & 15.982  & 75.049 &  366.87 &  2023.7 &  11990  & 75992 \\\vspace{10pt}
&1 & 1 & 5 & 16 & 75 & 366 & 2016 & 11936 & 75678\\
$a_3$ & 1.996  & -0.034 &  22.940 &  -0.950 &  684.23  & -22.334  & 34938 &  2360.8  & 2512126   \\\vspace{10pt}
&2 & 0 & 23 & 0 & 684 & 0 & 34760 & 0 & 2493096\\
$s_2^+$ & 2.000 &  0.001  & 10.998  & 0.969 &  94.977 &  35.182 &  1099.5 &  966.05 &  15812 \\\vspace{10pt}
&2 & 0 & 11 & 1 & 95 & 36 & 1099 & 982 & 15792\\
$s_3$& 2.996& 0.008  & 25.953 &  0.129  & 344.64 &  2.935  & 5759.4 &  73.138 &  111003 \\\vspace{10pt}
&3 & 0 & 26 & 0 & 345 & 0 & 5754 & 0 & 110586\\
$s_4^+$ & 3.002 &  0.980 &  27.023 &  25.574 &  384.80 &  697.45  & 7207.4 &  20004 &  163235 \\\vspace{10pt}
&3 & 1 & 27 & 26 & 385 & 708 & 7231 & 20296 & 164277\\
$s_5$ & 3.995 &  -0.036 &  41.906 & -0.719 & 658.28  & -16.625 & 12776 &  -428.23  & 281027 \\\vspace{10pt}
&4 & 0 & 42 & 0 & 660 & 0 & 12810 & 0 & 281736\\
$s_6^+$ & 5.001 &  0.968  & 63.005 &  45.334 &  1134.1  & 1782.0 &  25376  & 65650  & 663829\\\vspace{10pt}
&5 & 1 & 63 & 46 & 1135 & 1800 & 25431 & 66232 & 666387\\
$s_7$ & 6.000 &  0.002 &  90.015 &  0.356 &  1860.6 &  13.010 &  44746  & 380.75 & 1172844   \\
&6 & 0 & 90 & 0 & 1860 & 0 & 44730 & 0 & 1172556\\
\bottomrule
\end{tabular}
\vspace{12pt}
\caption{Central moment statistics of $a_k$ and $s_k$ in genus 3, $N=2^{25}$.}\label{table:g3stats}
\vspace{-20pt}
\small
$$y^2 = x^7 + 314159x^5 + 271828x^4 + 1644934x^3 +57721566x^2 +1618034x + 141021.$$
\normalsize
\end{center}
\end{table}

\begin{table}
\begin{center}
\begin{tabular}{rrrrrrrrrrr}
\toprule
$N$&\hspace{16pt}$M_1$& \hspace{16pt}$M_2$ & \hspace{16pt}$M_3$ & \hspace{16pt}$M_4$ & \hspace{16pt}$M_5$ & \hspace{16pt}$M_6$  & \hspace{16pt}$M_7$  & \hspace{16pt}$M_8$\\\midrule 
$2^{11}$&-0.071  & 1.031 &  -0.276  & 3.167  & -2.295  & 15.250  & -21.145   & 97.499\\
$2^{12}$&-0.036  & 1.112 &  -0.087 &  3.565 &  -0.475  & 17.251 &  -4.539  & 105.082\\
$2^{13}$&-0.067  & 1.085  & -0.249 &  3.407 &  -1.567 &  16.537  & -12.893 &  103.344\\
$2^{14}$&-0.046  & 1.029 &  -0.232 &  3.181 &  -1.529  & 15.795 &  -13.309  & 104.558\\
$2^{15}$&-0.044  & 1.031 &  -0.121&   3.256 &  -0.428   &16.325  & -2.396 &  107.173\\
$2^{16}$&-0.025  & 1.022 &  -0.069 &  3.143 &  -0.251  & 15.251 &  -1.673 &  96.837\\
$2^{17}$&-0.016  & 1.011 &  -0.041 &  3.079 &  -0.204 &  14.594 &  -1.717 &  88.871\\
$2^{18}$&-0.009  & 1.002 &  -0.022 &  3.041  & -0.138 &  14.441  & -1.456 &  88.636\\
$2^{19}$&-0.002  & 1.003 &  -0.013 &  3.031 &  -0.108 &  14.259  & -1.023 &  86.288\\
$2^{20}$& 0.001  & 0.998  & 0.001 &  3.003  & -0.041 &  14.126  & -0.687  & 85.815\\
$2^{21}$&-0.000  & 1.003  & -0.002 &  3.016  & -0.045  & 14.088 &  -0.577  & 84.746\\
$2^{22}$&0.002   &1.002 &  0.009  & 3.013  & 0.037 &  14.058  & 0.101  & 84.166\\
$2^{23}$&0.001   &1.001 &  0.002 &  3.006  & 0.001&   13.999  & -0.103&   83.715\\
$2^{24}$&0.000  & 1.001 &  0.001 &  3.002  & 0.008 &  13.964  & 0.036  & 83.346\\
$2^{25}$&0.000  & 1.000 &  -0.000 &  2.995 &  -0.010 &  13.950  & -0.120  & 83.500\\
$2^{26}$&0.000  & 1.001 &  -0.001  & 3.004 &  -0.006 &  14.014 &  -0.031  & 84.041\\
\bottomrule
\end{tabular}
\vspace{6pt}
\caption{Convergence of moment statistics for $a_1$ in genus 2 as $N$ increases.}\label{table:g2seqconverge}
\vspace{-20pt}
\small
$$y^2 = x^5 + 314159x^3 + 271828x^2 + 1644934x + 57721566.$$\vspace{12pt}
\normalsize
\begin{tabular}{rrrrrrrrrrr}
\toprule
$N$&\hspace{16pt}$M_1$& \hspace{16pt}$M_2$ & \hspace{16pt}$M_3$ & \hspace{16pt}$M_4$ & \hspace{16pt}$M_5$ & \hspace{16pt}$M_6$  & \hspace{16pt}$M_7$  & \hspace{16pt}$M_8$\\\midrule 
$2^{11}$&0.033  & 0.942 &  0.204  & 2.611 &  1.313  & 11.365  & 9.198 &  62.068\\
$2^{12}$&0.009  & 0.921 &  0.092  & 2.447  & 0.649  & 10.149 &  4.617 &  52.838\\
$2^{13}$&0.015  & 0.961 &  0.075 &  2.676 &  0.535  & 11.971 &  4.345 &  69.641\\
$2^{14}$&-0.011  & 0.983  & -0.060 &  2.893 &  -0.245 &  14.316  & 0.704 &  99.690\\
$2^{15}$&-0.005  & 1.011  & -0.018  & 3.134 &  -0.067  & 16.286  & 0.836 &  116.675\\
$2^{16}$&-0.017  & 1.007  & -0.054  & 3.154 &  -0.105  & 16.813  & 2.952  & 127.212\\
$2^{17}$&-0.006  & 0.993  & -0.024  & 3.027  & -0.041  & 15.431  & 1.622 &  109.717\\
$2^{18}$&-0.005  & 0.996  & -0.026 &  3.006 &  -0.110 &  15.196  & 0.239  & 106.901\\
$2^{19}$&-0.001  & 0.999  & -0.013  & 2.985 &  -0.087 &  14.793  & -0.418 &  101.662\\
$2^{20}$&0.000  & 0.989  & -0.007 &  2.934  & -0.072 &  14.440  & -0.759 &  98.109\\
$2^{21}$&0.003  & 0.997  & 0.003 & 2.979  & -0.017  & 14.796  & -0.562 &  101.690\\
$2^{22}$&0.002  & 0.999  & 0.005 &  3.003  & 0.038  & 15.098  & 0.446 &  105.733\\
$2^{24}$&0.000  & 1.001 &  0.001  & 3.015  & -0.005  & 15.138 &  -0.102 &  105.418\\
$2^{24}$&0.000 &  0.999  & -0.004  & 2.990  & -0.043 &  14.916  & -0.397  & 103.271\\
$2^{25}$&0.000  & 0.999 &  -0.007 &  2.995  & -0.046 &  14.968  & -0.343  & 103.755\\
\bottomrule
\end{tabular}
\vspace{6pt}
\caption{Convergence of moment statistics for $a_1$ in genus 3 as $N$ increases.}\label{table:g3seqconverge}
\vspace{-20pt}
\small
$$y^2 = x^7 + 314159x^5 + 271828x^4 + 1644934x^3 +57721566x^2 +1618034x + 141021.$$
\normalsize
\end{center}
\end{table}
\clearpage
After perusing the data in Tables \ref{table:g1stats}-\ref{table:g3seqconverge}, several questions come to mind:
\begin{enumerate}
\item
What is the rate of convergence as $N$ increases?
\item
Are these results representative of typical curves?
\item
Can we distinguish exceptional distributions that may arise?
\end{enumerate}

Question 1 has been considered in genus 1 (but not, to our knowledge, in higher genera).  For an elliptic curve without complex multiplication, the conjectured discrepancy between the observed distribution and the Sato-Tate prediction is $O(N^{-1/2})$ (see Conjecture 1 of \cite{Akiyama:SatoTateConvergence} or Conjecture 2.2 of \cite{Mazur:SatoTate}).  This conjecture implies that the generalized Riemann Hypothesis then holds for the $L$-series of the curve \cite[Theorem 2]{Akiyama:SatoTateConvergence}.  We will not attempt to address this question here, other than noting that the figures listed in Tables \ref{table:g2seqconverge} and \ref{table:g3seqconverge} are not inconsistent with a convergence rate of $O(N^{-1/2})$.

\subsection{Random Families of Genus 2 Curves}

We can say more about Questions 2 and 3, at least in genus 2.  To address Question 2 we tested over a million randomly generated curves of the form
$$y^2=x^5+f_4x^4+f_3x^3+f_2x^2+f_1x+f_0,$$
with the integer coefficients $f_0$, \ldots, $f_5$ obtained from a uniform distribution on the interval $[-2^{63}+1,2^{63}-1]$.

Table \ref{table:RandomCurves} describes the distribution of moment statistics for $a_1$ over three sets of computations:  one million curves with $N=2^{16}$, ten thousand curves with $N=2^{20}$, and one hundred curves with $N=2^{24}$.  The rows list bounds on the deviation from the moment sequence for $a_1$ in $USp(2g)$ that apply to $m\%$ of the curves, with $m$ equal to 50, 90, or 99.  One sees close agreement with the predicted moment statistics, with $\Delta M_n$ decreasing as $N$ increases.  The maximum deviation in $M_4$ observed for any curve was 0.56 with $N=2^{16}$.

\begin{table}
\begin{center}
\begin{tabular}{lrrrrrrrrrr}
\toprule
\# & $N$&\%&\hspace{4pt}$\Delta M_1$& $\Delta M_2$ & $\Delta M_3$ & $\Delta M_4$ & $\Delta M_5$ & $\Delta M_6$  & $\Delta M_7$  & $\Delta M_8$\\\midrule
$10^6$&$2^{16}$&50&0.008&0.012&0.031&0.072&0.204&0.563&1.685&5.104\\
&&90&0.020&0.029&0.076&0.177&0.497&1.371&4.120&12.397\\
&&99&0.032&0.045&0.120&0.277&0.781&2.156&6.512&19.633\\\midrule
$10^4$&$2^{20}$&50&0.002&0.003&0.009&0.020&0.057&0.159&0.470&1.433\\
&&90&0.006&0.008&0.021&0.049&0.138&0.384&1.154&3.485\\
&&99&0.009&0.013&0.033&0.078&0.214&0.604&1.801&5.432\\\midrule
$10^2$&$2^{24}$&50&0.001&0.001&0.002&0.005&0.017&0.044&0.138&0.424\\
&&90&0.002&0.002&0.006&0.013&0.035&0.101&0.277&0.933\\
&&99&0.002&0.003&0.008&0.019&0.054&0.165&0.543&1.519\\
\bottomrule
\end{tabular}
\vspace{6pt}
\caption{Moment deviations for families of random genus 2 curves.}\label{table:RandomCurves}
\end{center}
\end{table}

We also looked for exceptional distributions among the outliers, considering the possibility that one or more curves in our random sample might not have large Galois image.  From our initial family of one million random curves we selected one thousand curves whose moment statistics showed the greatest deviation from the predicted values.  We recomputed the $a_1$ moment statistics of these curves, with the bound $N$ increased from $2^{16}$ to $2^{20}$.  In each and every case, we saw convergence toward the moment sequence for $a_1$ in $USp(4)$,
\begin{align}
M[a_1] = 1,\enspace 0,\enspace 1,\enspace 0,\enspace 3,\enspace 0,\enspace 14,\enspace 0,\enspace 84,\enspace 0,\enspace 594,\enspace \ldots,\tag{\Aseq{138349}}
\end{align}
as predicted by Conjecture \ref{conjecture:GeneralSatoTate} for curves with large Galois image.  An additional test of one hundred of the most deviant curves from within this group with $N=2^{24}$ yielded further convergence, with $\Delta M_n < 1$ for all $n \le 8$.  This is strong evidence that all the curves in our original random sample had large Galois image; every exceptional distribution we have found for $a_1$ in genus 2 has $\Delta M_8 > 200$.

One may ask whether convergence of the $a_1$ moment statistics to $M[a_1]$ is enough to guarantee that the $L$-polynomial distribution of the curve is represented by $USp(2g)$, since we have not examined the distribution of $a_2$.  If we assume the distribution of $\barLp$ is given by some infinite compact subgroup of $USp(4)$ (as in Conjecture \ref{conjecture:GeneralSatoTate}), then it suffices to consider $a_1$.  In fact, under this assumption, much more is true:  if the fourth moment statistic of $a_1$ converges to $3$, then the distribution of $\barLp$ converges to the distribution of $\chi(T)$ in $USp(4)$.  This remarkable phenomenon is a consequence of {\em Larsen's alternative} \cite{Larsen:SatoTate,Katz:LarsensAlternative}.

\subsection{Larsen's Alternative}

To apply Larsen's alternative we need to briefly introduce a representation theoretic definition of ``moment" which will turn out to be equivalent to the usual statistical moment in the case of interest to us.  Here we parallel the presentation in \cite[Section 1.1]{Katz:LarsensAlternative}, but assume $G$ to be compact rather than reductive.  Let $V$ be a complex vector space of dimension at least two and $G\subset GL(V)$ a compact group.  Define
\begin{equation}
M_{a,b}(G,V) = \dim_\mathbb{C}(V^{\otimes a}\otimes \check{V}^{\otimes b})^G,
\end{equation}
and set $M_{2n}(G,V)=M_{n,n}(G,V)$.  Let $\chi(A)=\tr(A)$ denote the character of $V$ as a $G$-module (the standard representation of $G$).  We then have
$$M_{2n}(G,V) = \int_{G}\chi(A)^n\bar\chi(A)^n dA = \int_{G}|\chi(A)|^{2n}dA.$$
We now specialize to the case $V=\mathbb{C}^{2g}$ and suppose $G\subset USp(2g)$.  Then
$$M_{2n}(G,V) = \int_{G}(\tr(A))^{2n}dA = E[(\tr(A))^{2n}] = M[a_1](2n),$$
where $a_1=-\tr(A)$ and $M[a_1](n) = E[a_1^n]$ as usual.  We can now state Larsen's alternative as it applies to our present situation.

\begin{theorem}[Larsen's Alternative]
Let $V$ a complex vector space of even dimension greater than $2$ and suppose $G$ is a compact subgroup of $USp(V)$.
If $M_{4}(G,V) = 3$, then either $G$ is finite or $G=USp(V)$.
\end{theorem}

This is directly analogous to Part 3 of Theorem 1.1.6 in \cite{Katz:LarsensAlternative}, and the proof is the same.

\begin{corollary}
Let $C$ be a curve of genus $g>1$.  Under Conjecture \ref{conjecture:GeneralSatoTate}, the distribution of $\barLp$ converges to the distribution of $\chi(T)$ in $USp(2g)$ if and only if the fourth moment statistic of $a_1$ converges to $3$.
\end{corollary}

The corollary provides a wonderfully effective way for us to distinguish curves with exceptional $\barLp$ distributions.

\section{Exceptional $\barLp$ Distributions in Genus 2}\label{section:g2Exceptions}

While we were unable to find {\em any} exceptional $\barLp$ distributions among random genus 2 curves with large coefficients, if one restricts the size of the coefficients such cases are readily found.  We tested every curve of the form $y^2=f(x)$, with $f(x)$ a monic polynomial of degree 5 with coefficients in the interval $[-64,64]$, more than $2^{35}$ curves.  As not every hyperelliptic curve of genus 2 can be put in this form, we also included curves with $f(x)$ of degree 6 (not necessarily monic) and coefficients in the interval $[-16,16]$.  With such a large set of curves to test, we necessarily used a much smaller value of $N$, approximately $2^{12}$.  We computed the moment statistics of $a_1$ using parallel point-counting techniques described in \cite{KedlayaSutherland:HyperellipticLSeries} to process 32 curves at once.

To identify exceptional curves, we applied a heuristic filter to bound the deviation of the fourth and sixth moment statistics from the $USp(2g)$ values $M[a_1](4)=3$ and $M[a_1](6)=14$.  By Larsen's alternative, it suffices only to consider the fourth moment, however we found the sixth moment also useful as a distinguishing metric: the smallest sixth moment observed among any of the exceptional distributions was 35, compared to 14 in the typical case.  A combination of the two moments proved to be most effective.

Searching for a small subset of exceptional curves in a large family using a statistical test necessarily generates many false positives: nonexceptional curves which happen to deviate significantly from the $USp(2g)$ distribution for $p\le N$.  The filter criteria were tuned to limit this, at the risk of introducing more false negatives (unnoticed exceptional curves).  After filtering the entire family with $N\approx 2^{12}$, the remaining curves were filtered again with $N=2^{16}$ to remove false positives.  Finally, we restricted the resulting list to curves with distinct Igusa invariants \cite{Igusa:Genus2Invariants}, leaving a set of some 30,000 nonisomorphic curves with (apparently) exceptional distributions.

One additional criterion used to distinguish distributions was the ratio $z(C,N)$ of zero traces, that is, the proportion of primes $p$ for which $a_p=0$, among $p\le N$ where $C$ has good reduction.  For a typical curve, $z(C,N)\to 0$ as $N\to\infty$, but for many exceptional distributions, $z(c,N)$ converges to a nonzero rational number.  In most cases where this arises, one can readily compute
\begin{equation}
\lim_{N\to\infty} z(C,N) = z(C),
\end{equation}
using the Hasse-Witt matrix.  This is described in detail in \cite{Sutherland:HasseWitt}, and the following is a typical example.  One can show that the curve $y^2=x^6+2$ has $a_p=0$ unless $p$ is of the form $p=6n+1$, in which case
\begin{equation}
a_p\equiv \binom{3n}{n}2^n(2^n+1) \mod p.
\end{equation}

It follows that for $p=6n+1$ we have $a_p=0$ if and only if $2^n\equiv -1 \bmod p$.  The integer $2^n=2^{(p-1)/6}$ is necessarily a sixth root of unity mod $p$, and exactly one of these is congruent to -1.  By the \v{C}ebotarev density theorem, this occurs for a set of density 1/6 among primes of the form $p=6n+1$.  Combine this with the fact that $a_p=0$
when $p$ is not of this form and we obtain $z(C)=7/12$.

Under Conjecture 2, one can show that for any curve $C$ (of arbitrary genus), $z(C)$ must exist and is a rational number, however we need not assume this here.\footnote{A compact subgroup of $USp(2g)$ has a finite number of connected components.  The density of zero trace elements must be zero or one on each component.}  For each nonzero $z(C)$ in Table \ref{table:g2a1dist}, we have identified a specific curve exhibiting the distribution for which one can show of $\lim_{N\to\infty}z(C,N) = z(C)$.  Typically, the Hasse-Witt matrix gives a lower bound on $\liminf_{N\to\infty}z(C,N)$, and by computing $G[\ell]$, the mod $\ell$ image of the Galois representation in $GSp(2g,Z/\ell z)$, one obtains the density of nonzero traces mod $\ell$, establishing an upper bound on $\limsup_{N\to\infty}z(c,N)$.  It is generally not difficult to find an $\ell$ for which these two bounds are equal (the \v{C}ebotarev density theorem is invoked on both sides of the argument).

After sorting the sequences of moment statistics and considering the values of $z(C)$ among our set of more than 30,000 exceptional curves, we were able to identify only 22 distributions that were clearly distinct (within the precision of our computations).  We also tested a wide range of genus 2 curves taken from the literature \cite{Cassels:Prolegomena,Frey:ModularCurves,Gonzalez:ModularGenus2,Hashimoto:SatoTateGenus2,
Howe:TorsionSubgroupsJacobians,Kuwata:SplitJacobians,Rodriguez-Villegas:SplitCM,Smart:Genus2,
Wamelen:Genus2CM,Wamelen:ProvingCM}, most with coefficient values outside the range of our search family.  In every case the $a_1$ moment statistics appeared to match one of our previously identified distributions.  Conversely, several of the distributions found in our search did not arise among the curves we tested from the literature.

Table \ref{table:g2a1dist} lists the 23 distinct distributions for $a_1$ we found for genus two curves, including the typical case, which is listed first. The value of $z(C)$ and the first six moments of $a_1$ suffice to distinguish every distribution we have found.  We also list the eighth moment statistic, which, while not accurate to the nearest integer, is almost certainly within one percent of the ``true" value.  We list only the moment statistics of $a_1$.  Histograms of the first twelve $a_1$ distributions can be found in Appendix II.  Additional $a_1$ histograms, along with moment statistics and histograms for $a_2$ and $s_k$ are available at \url{http://math.mit.edu/~drew/}.

\setlength\extrarowheight{3pt}
\begin{table}
\begin{center}
\begin{tabular}{rcrrrrcl}
\toprule
\#&$\enspace z(C)\enspace $&$M_2$ & $\enspace M_4$ & $\quad M_6$ & $\quad\enspace M_8$&$\space$& $f(x)$\\\midrule
1&0&1&3&14&84&&$x^5+x+1$\\
2&0&2&10&70&$588^*$&&$x^5-2x^4+x^3+2x-4$\\
3&0&2&11&90&$888^*$&&$x^5+20x^4-26x^3+20x^2+x$\\
4&0&2&12&110&$1203^*$&&$x^5+4x^4+3x^3-x^2-x$\\
5&0&4&32&320&$3581^*$&&$x^5+7x^3+32x^2+45x+50$\\
6&1/6&2&12&100&$979^*$&&$x^5-5x^3-5x^2-x$\\
7&1/4&2&12&100&$1008^*$&&$x^5+2x^4+2x^2-x$\\
8&1/4&2&12&110&$1257^*$&&$x^5-4x^4-2x^3-4x^2+x$\\
9&1/2&1&5&35&$293^*$&&$x^5-2x^4+11x^3+4x^2+4x$\\
10&1/2&1&6&55&$601^*$&&$x^5-2x^4-3x^3+2x^2+8x$\\
11&1/2&2&16&160&$1789^*$&&$x^5+x^3+x$\\
12&1/2&2&18&220&$3005^*$&&$x^5-3x^4+19x^3+4x^2+56x-12$\\
13&1/2&4&48&640&$8949^*$&&$x^6+1$\\
14&7/12&1&6&50&$489^*$&&$x^5-4x^4-3x^3-7x^2-2x-3$\\
15&7/12&2&18&200&$2446^*$&&$x^6+2$\\
16&5/8&1&6&50&$502^*$&&$x^5+x^3+2x$\\
17&5/8&2&18&200&$2515^*$&&$x^5-10x^4+50x^2-25x$\\
18&3/4&1&8&80&$894^*$&&$x^5-2x^3-x$\\
19&3/4&1&9&100&$1222^*$&&$x^5-1$\\
20&3/4&1&9&110&$1501^*$&&$11x^6+11x^3-4$\\
21&3/4&2&24&320&$4474^*$&&$x^5+x$\\
22&13/16&1&9&100&$1254^*$&&$x^5+3x$\\
23&7/8&1&12&160&$2237^*$&&$x^5+2x$\\
\bottomrule
\end{tabular}
\vspace{6pt}
\caption{Moments of $a_1$ for genus 2 curves $y^2=f(x)$ with $N=2^{26}$.}\label{table:g2a1dist}
\end{center}
\begin{minipage}{5in}
\small
Column $z(C)$ is the density of zero traces ($a_p$ values).  The starred values indicate uncertainty in the eighth moment statistic.  In each case, if $T_8=\lim_{N\to\infty}M_8$, we estimate that $-0.005\le 1-M_8/T_8\le 0.01$ with very high probability (the larger uncertainty on the positive side is primarily due to an observed excess of zero traces for small values of $N$, we expect $M_8\le T_8$ in most cases).  See Table \ref{table:USp4groups} for predicted values of $T_8$ and $T_{10}$.
\end{minipage}
\normalsize
\end{table}

The third distribution in Table \ref{table:g2a1dist} went unnoticed in our initial analysis (we later found several examples that had been misclassified) and is not a curve taken from the literature.  We constructed the curve
\begin{align}
\qquad y^2 = x^5 + 20x^4 - 26x^3 + 20x^2 + x\tag{$C$}
\end{align}
to have a split Jacobian, isogenous to the product of the elliptic curve
\begin{align}
y^2 = x^3 - 11x +14,\tag{$E_1$}
\end{align}
which has complex multiplication, and the elliptic curve
\begin{align}
y^2=x^3 + 4x^2 - 4x,\tag{$E_2$}
\end{align}
which does not.  For every $p$ where $C$ has good reduction, the trace of Frobenius is simply the sum of the traces of $E_1$ and $E_2$.  As $E_1$ and $E_2$ are not isogenous (over $\mathbb{C}$), we expect their $a_1$ distributions to be uncorrelated.\footnote{The genus 1 traces may be correlated in a way that does not impact $a_1=-a_p/\sqrt{p}$, e.g., both curves might have the property that $p\bmod 3$ determines $a_p\bmod 3$.}  It follows that the moment sequence of $a_1$ for the curve $C$ is simply the binomial convolution of the moment sequences of $a_1$ for the curves $E_1$ and $E_2$.  Thus we expect the moment statistics of $a_1$ for $C$ to converge to
\begin{align}
1,\enspace 0,\enspace 2, \enspace 0, \enspace 11, \enspace 0, \enspace 90, \enspace 0, \enspace 889, \enspace 0, \enspace 9723, \enspace \ldots,\tag{\Aseq{138552}*}
\end{align}
which is the binomial convolution of the sequences $(1, 0, 1, 0, 3, 0, 10, 0, 35, \ldots)$ and $(1, 0, 1, 0, 2, 0, 5, 0, 14, \ldots)$ mentioned in the introduction.  These are the $a_1$ moment sequences of elliptic curves with and without complex multiplication.  In terms of moment generating functions, we simply have
$$\mathcal{M}_C[a_1] = \mathcal{M}_{E_1}[a_1]\mathcal{M}_{E_2}[a_1].$$

Provided the covariance between the $a_1$ distributions of $E_1$ and $E_2$ is zero, one can prove the $a_1$ moment statistics of $C$ converge to the sequence above using known results for genus 1 curves (note that $E_1$ has complex multiplication and $E_2$ has multiplicative reduction at $p=107$).

Many of the distributions in Table \ref{table:g2a1dist} can be obtained from genus 1 moment sequences by constructing an appropriate genus 2 curve with split Jacobian, as shown in the next section.  It is important to note that a curve whose Jacobian is simple (not split over $\overline{\mathbb{Q}}$) may still have an $L$-polynomial distribution matching that of a split Jacobian.  Distribution \#2, for example, corresponds to a Jacobian which splits as the product of two nonisogenous elliptic curves without complex multiplication.  This distribution also arises for some genus 2 curves with simple Jacobians, including
$$y^2 = x^5 - x^4 + x^3 + x^2 - 2x + 1.$$
This is a modular genus 2 curve which appears as $C_{188,A}$ in \cite{Gonzalez:ModularGenus2}, along with many similar examples.

We speculate that simple Jacobians with Distribution \#2 are all of type I(2) in the classification of Moonen and Zarhin \cite{Moonen:HodgeClasses}, corresponding to Jacobians whose endomorphism ring is isomorphic to the ring of integers in a real quadratic extension of $\mathbb{Q}$.  A similar phenomenon occurs with Distribution \#11, which arises for split Jacobians that are isogenous to the product of an elliptic curve and its twist, but also for simple Jacobians of type II(1) in the Moonen-Zarhin classification (these are QM-curves, see \cite{Hashimoto:SatoTateGenus2} for examples).  The remaining two types of simple genus 2 Jacobians in the Moonen-Zarhin classification are type I(1), which is the typical case (Jacobians with endomorphism ring $\mathbb{Z}$), and type IV(2,1), which occurs for curves with complex multiplication over a quartic CM field (with no imaginary quadratic subfield).  These correspond to Distributions \#1 and \#19 respectively (examples of the latter can be found in \cite{Wamelen:Genus2CM}).  The remaining distributions appear to arise
only for curves with split Jacobians.

\section{Representation of Genus 2 Distributions in $USp(4)$}\label{section:g2groups}

Conjecture \ref{conjecture:GeneralSatoTate} implies that each distribution in Table \ref{table:g2a1dist} is represented by the distribution of characteristic polynomials in some infinite compact subgroup $H$ of $USp(4)$.  In this section we will exhibit such an $H$ for each distribution.  We do not claim that each $H$ we give is the ``correct" subgroup for every curve with the corresponding distribution, in the sense of corresponding to the Galois image $G_\ell$, as discussed in Section \ref{section:SatoTate}.\footnote{Indeed, a single $H$ for each distribution would not suffice.  As noted above, two curves may share the same $\barLp$ distribution for different reasons.}  Rather, for each $H$ we show that its density of zero traces and moment sequence are compatible with the corresponding data in Table \ref{table:g2a1dist}.  In most cases we also have evidence that suggests $H$ is the correct subgroup for the particular corresponding curve listed in Table \ref{table:g2a1dist} (see Section \ref{section:evidence}).

For all but two cases we will construct $H$ using two subgroups of $USp(2)$ that represent distributions of genus 1 curves: $G_1 = USp(2)$ for an elliptic curve without complex multiplication, and $G_2 = N(SO(2))$, the normalizer of $SO(2)$ in $SU(2)$, for an elliptic curve with complex multiplication.  We will construct each $H\subset USp(4)$ from $G_1$ and/or $G_2$ explicitly as a group of matrices, however the motivation behind these constructions are split Jacobians.

The most obvious cases correspond to the product of two nonisogenous elliptic curves, and we have the groups $G_1\times G_1$, $G_1\times G_2$, and $G_2\times G_2$ as groups of block diagonal $2\times 2$ matrices in $USp(4)$.  These correspond to Distributions \#2, \#3, and \#8 in Table \ref{table:g2a1dist}, and their moment sequences are easily computed via binomial convolutions of the appropriate genus 1 moment sequences (the example of the previous section corresponds to $G_1\times G_2$).

To obtain additional distributions, we also consider the product of two isogenous elliptic curves.  We may, for example, pair an elliptic curve with an isomorphic copy of itself, or with one of its twists.\footnote{See \cite[Ch. 14]{Cassels:Prolegomena} and \cite{Kuwata:SplitJacobians} for explicit methods of constructing such curves.}  For two isomorphic curves, the corresponding subgroup $H_i$ contains block diagonal matrices of the form
\begin{equation*}
B = \left(\begin{matrix}
A & 0\\
0 & A\\
\end{matrix}\right)
\end{equation*}
where $A$ is an element of $G_i$ ($i=1$ or 2).  To pair a curve with its twist we also include block diagonal matrices with $A$ and $-A$ on the diagonal to obtain the subgroup $H_i^{-}$.

We now generalize this idea.  Let $G=G_1$ (resp. $G_2$) be a compact subgroup of $USp(2)$, and let $G^*$ be the subgroup of $U(2)$ obtained by extending $G$ by scalars and taking the subgroup of elements whose determinants are $k$th roots of unity (for some positive integer $k$).  For $A\in G^*$, let $\overline{A}$ denote the complex conjugate of $A$ and define the block diagonal matrix
\begin{equation*}
B = \left(\begin{matrix}
A & 0\\
0 & \overline{A}\\
\end{matrix}\right).
\end{equation*}
The matrix $B$ is clearly unitary, and one easily verifies that it is also symplectic, hence $B\in USp(4)$.  The set of all such $B$ forms our subgroup $H$.  As a topological group, $H$ has $k$ (resp. $2k$) connected components, each a closed set consisting of elements with $A$ having a fixed determinant, thus $H$ is compact.  The identity component is isomorphic to $USp(2)$ (resp. $SO(2)$), embedded diagonally.

We may write $A\in G^*$ as $A = \omega^j A_0$ with $A_0\in G\subseteq USp(2)$, $\omega$ a primitive $2k$th root unity, and $1\le j \le k$.  We then have
\begin{equation*}
\tr(B) = \tr(A)+\tr(\overline{A}) = \omega^j\tr(A_0)+\omega^{-j}\tr(\overline{A}_0) = (\omega^j+\omega^{-j})\tr(A_0),
\end{equation*}
since $\tr(\overline{A}_0)=\tr(A_0)$ for any $A_0\in USp(2)$.  It follows that 
\begin{equation}\label{equation:expHk}
\Exp_H[(\tr(B))^n] =\left(\frac{1}{k}\sum_{j=1}^k(\omega^j+\omega^{-j})^n\right)\Exp_G[(\tr(A_0))^n]
\end{equation}
where $\Exp_G[X]$ denotes the expectation of a random variable $X$ over the Haar measure on $G$.

The term $(\omega^j+\omega^{-j})^n\Exp_G[(\tr(A_0))^n]$ corresponds to the $n$th moment of the trace distribution on a component of $H$.  These moments will all be integers precisely when $k\in\{1,2,3,4,6\}$ (they are zero for $n$ odd).  In fact, these are the only values of $k$ for which $H$ plausibly represents the distribution of $\barLp$ for a genus 2 curve defined over $\mathbb{Q}$, as we now argue.

Consider an $L$-polynomial $L_p(T)$ for which $\barLp = L_p(p^{-1/2}T)$ is the characteristic polynomial of some $B\in H$.  We may factor $L_p(T)$ as
\begin{equation}\label{equation:splitLp}
L_p(T) = (\chi(p)pT^2-\alpha T + 1)(\overline{\chi}(p)pT^2 - \overline{\alpha} T +1),
\end{equation}
where $\alpha\in \mathbb{Z}[\zeta]$, with $\zeta$ a primitive $k$th root of unity, and $\chi(p)$ is a $k$th root of unity equal to the determinant of $A$ in the component of $H$ containing $B$.  For the elliptic curve we have in mind as a factor of the Jacobian, $\chi(p)p$ is the determinant of its Frobenius element and $\alpha$ is the trace.  Since $L_p(T)$ has integer coefficients, $\alpha$ must lie in a quadratic extension of $\mathbb{Q}$, giving $k\in\{1,2,3,4,6\}$.

For each of these $k$ we can easily compute closed forms for the parenthesized expression in (\ref{equation:expHk}).  Assuming $n > 0$ is even, we obtain
$$2^n,\quad 2^{n}/2, \quad (2^n+2)/3, \quad (2^n+2^{n/2+1})/4, \quad (2^n+2\cdot 3^{n/2}+2)/6,$$
for $k=1,2,3,4,6$, respectively.  Since $\Exp_G[(\tr(A_0))^n]$ is zero for odd values of $n$, these expressions may be used in (\ref{equation:expHk}) for all positive $n$.  For even values of $k$, the density $z(H)$ of zero traces in $H$ is $1/k$ (resp. $(1+1/k)/2$), and $z(H)$ is zero (resp. 1/2) for $k$ odd.

For any of the subgroups $H$ constructed as above, we may also consider the group $J(H)$ generated by $H$ and the block diagonal matrix
\begin{equation}
J = \left(\begin{matrix}
\enspace 0\enspace & I\\
-I & 0\\
\end{matrix}\right).
\end{equation}
The group $J(H)$ contains $H$ as a subgroup of index 2, with the nonidentity coset having all zero traces.  For $n>0$, the $n$th moment of the trace distribution in $J(H)$ is simply half that of $H$, and $z(J(H)) = (z(H)+1)/2$.

For $i=1$ or 2 and $k\in\{1,2,3,4,6\}$, let $H^k_i$ denote the group $H$ constructed using $G_i$ and $k$.  By also considering $J(H^k_i)$, we can construct a total of 20 groups, 18 of which have distinct eigenvalue distributions.  With the sole exception of $J(H^6_2)$, each of these matches one of the distributions in Table \ref{table:g2a1dist} to a proximity well within the accuracy of our computational methods.  Note that $H^1_i = H_i$, and the group $H^2_i$ has the same eigenvalue distribution as $H_i^{-}$ (but is not conjugate).

\begin{table}
\begin{center}
\begin{tabular}{lccccc}
\toprule
$\quad H\quad$&$\quad k=1\quad$&$\quad k=2\quad$&$\quad k=3\quad$&$\quad k=4\quad$&$\quad k=6\quad$\\
\midrule
$H^k_1$&  5 & 11 & 4 & 7 & 6\\
$J(H^k_1)$& 11 & 18 & 10 & 16 & 14\\
$H^k_2$& 13 & 21 & 12 & 17 & 15\\
$J(H^k_2)$& 21 & 23 & 20 & 22 & *\\
\bottomrule
\end{tabular}
\vspace{6pt}
\caption{Distributions matching $H^k_i$ or $J(H^k_i)$.}
\end{center}
\end{table}

We may also consider $J(G_1\times G_1)$, which corresponds to Distribution \#9.  This construction does not readily apply to $G_1\times G_2$ (in fact $J(G_1\times G_2)=J(G_1\times G_1)$).  The group $J(G_2\times G_2)$ does give a new distribution (it is the normalizer of $SO(2)\times SO(2)$ in $USp(4)$), but it does not correspond to any we have found.  However, the group $J(G_2\times G_2)$ contains a subgroup $K$ not equal to $G_2\times G_2$ which matches Distribution \#19.  The group $K$ has identity component $SO(2)\times SO(2)$ and a cyclic component group, of order 4 (this determines $K$).

Of all the groups we have constructed, only $K$ does not correspond, in some fashion, to a split Jacobian.  As noted above, Distribution \#19 arises for curves with simple Jacobians and complex multiplication over a quartic $CM$ field, and in fact all nineteen such curves documented by van Wamelen have this distribution \cite{Wamelen:Genus2CM,Wamelen:ProvingCM}.  The only remaining group to consider is $USp(4)$ itself, which of course gives Distribution \#1.

Table \ref{table:USp4groups} gives a complete list of the subgroups of $USp(4)$ we have identified,
one for each distribution in Table \ref{table:g2a1dist}, and also entries for $J(H^6_2)$ and $J(G_2\times G_2)$. These last two distributions appear to
be spurious; see \S~\ref{section:nonexistence}.

\setlength\extrarowheight{3pt}
\begin{table}
\begin{center}
\begin{tabular}{rcccccrrrrrr}
\toprule
\#&$\qquad H\qquad$&$d$&$c(H)$&$\enspace z(H)\enspace$&$M_2$ &$M_4$ &$M_6$ & $M_8$&$M_{10}$\\\midrule
1&$USp(4)$&10&1&0&1&3&14&84&594\\
2&$G_1\times G_1$&6&1&0&2&10&70&588&5544\\
3&$G_1\times G_2$&4&2&0&2&11&90&889&9723\\
4&$H_1^3$&3&3&0&2&12&110&1204&14364\\
5&$H_1$&3&1&0&4&32&320&3584&43008\\
6&$H_1^6$&3&6&1/6&2&12&100&980&10584\\
7&$H_1^4$&3&4&1/4&2&12&100&1008&11424\\
8&$G_2\times G_2$&2&4&1/4&2&12&110&1260&16002\\
9&$J(G_1\times G_1)$&6&2&1/2&1&5&35&294&2772\\
10&$J(H_1^3)$&3&6&1/2&1&6&55&602&7182\\
11&$H_1^-$&3&2&1/2&2&16&160&1792&21504\\
12&$H_2^3$&1&6&1/2&2&18&220&3010&43092\\
13&$H_2$&1&2&1/2&4&48&640&8960&129024\\
14&$J(H_1^6)$&3&12&7/12&1&6&50&490&5292\\
15&$H_2^6$&1&12&7/12&2&18&200&2450&31752\\
16&$J(H_1^4)$&3&8&5/8&1&6&50&504&5712\\
17&$H_2^4$&1&8&5/8&2&18&200&2520&34272\\
18&$J(H_1^-)$&3&4&3/4&1&8&80&896&10752\\
19&$K$&2&4&3/4&1&9&100&1225&15876\\
20&$J(H_2^3)$&1&12&3/4&1&9&110&1505&21546\\
21&$H_2^-$&1&4&3/4&2&24&320&4480&64512\\
22&$J(H_2^4)$&1&16&13/16&1&9&100&1260&17136\\
23&$J(H_2^-)$&1&8&7/8&1&12&160&2240&32256\\\midrule
*&$J(G_2\times G_2)$&2&8&5/8&1&6&55&630&8001\\
*&$J(H_2^6)$&1&24&19/24&1&9&100&1225&15876\\
\bottomrule
\end{tabular}
\vspace{6pt}
\caption{Candidate subgroups of $USp(4)$.}\label{table:USp4groups}
\end{center}
\begin{minipage}{5in}
\small
Row numbers correspond to the distributions in Table \ref{table:g2a1dist}.  Column $d$ lists the real dimension of $H$, and $c(H)$ counts its components.  The column $z(H)$ gives the density of zero traces, and $M_n=\Exp_H[(\tr(B))^n]$ for a random $B\in H$.  The last two rows are not known to match the $L$-polynomial distribution of a genus 2 curve.
\normalsize
\end{minipage}
\end{table}

\subsection{Supporting Evidence}\label{section:evidence}

Aside from closely matching the trace distributions we have found, there is additional data that supports our choice of the subgroups $H$ appearing in Table \ref{table:USp4groups}.  First, we find that these $H$ not only match the distribution of the $a_1$ coefficient in $\barLp$, they also appear to give the correct distribution of $a_2$.  We should note that for the three curves in Table \ref{table:g2a1dist} where $f(x)$ has degree 6, the available methods for computing $L_p(T)$ are much less efficient, so we did not attempt this verification in these three cases.\footnote{The $a_1$ coefficient can be computed reasonably efficiently in the degree 6 case by counting points on $C$ in $\Fp$, but the group law on the Jacobian is much slower.}

More significantly, the disconnected groups in Table \ref{table:USp4groups} also appear to give the correct distribution of $a_1$ (and $a_2$ for $f(x)$ of degree 5) on each of their components.  If we partition the components of $H$ according to their distributions of characteristic polynomials, for a given curve we can typically find a partitioning of primes into sets of corresponding density with matching $\barLp$ distributions.

Taking Distribution \#10 as an example, we have $H=J(H_1^3)$ in Table \ref{table:USp4groups} and the corresponding curve in  Table \ref{table:g2a1dist} is given by
$$y^2=f(x)=x^5-2x^4-3x^3+2x^2+8x=x(x-2)(x^3-3x-4).$$
The cubic $g(x)=x^3-3x-4$ has Galois group $S_3$.  The set of primes $P_3$ for which $g(x)$ splits into three factors in $\Fp[x]$ has density $1/6$, and corresponds to the identity component of $J(H_1^3)$.\footnote{We thank Dan Bump for suggesting this approach.}  The set of primes $P_2$ where $g(x)$ splits into two factors has density 1/2, and corresponds to the nonidentity coset of $H_1^3$ in $J(H_1^3)$, containing three components with identical eigenvalue distributions.  The remaining set of primes $P_1$ for which $g(x)$ is irreducible has density 1/3 and corresponds to the set of elements of $H_1^3$ for which the determinant of the block diagonal matrix $A$ is not 1 (this includes two of the six components of $J(H_1^3)$).  Table \ref{table:components} lists the moment statistics for $a_1$ and $a_2$, restricted to the sets $P_1$, $P_2$, and $P_3$, with values for the corresponding subset of $J(H_1^3)$ beneath.

\setlength\extrarowheight{2pt}
\begin{table}
\begin{center}
\begin{tabular}{lcrrrrrrrr}
\toprule
Set&$\quad X\quad $&$M_1$ & $M_2$ &$M_3$ & $M_4$ &$M_5$ &$M_6$ &$M_7$ &$M_8$\\\midrule
$P_1$&$a_1$&0.000&1.001&-0.002&2.002&-0.007&5.005&-0.027&14.01\\
   &&0&1&0&2&0&5&0&14\\
&$a_2$&0.001&1.001&1.001&3.002&6.005&15.01&36.03&91.10\\
   &  &0&1&1&3&6&15&36&91\\\midrule
$P_2$&$a_1$&0.000&0.000&0.000&0.000&0.000&0.000&0.000&0.000\\
   &&0&0&0&0&0&0&0&0\\
&$a_2$&1.000&2.010&3.001&6.003&10.00&20.01&35.02&70.05\\
   &&1&2&3&6&10&20&35&70\\\midrule
$P_3$&$a_1$&-0.002&3.999&-0.045&31.95&-0.590&319.3&-7.69&3574\\
   &&0&4&0&32&0&320&0&3584\\
   &$a_2$&2.999&9.995&36.97&149.8&652.7&3005&14404&71160\\
   &&3&10&37&150&654&3012&14445&71398\\
\bottomrule
\end{tabular}
\vspace{6pt}
\caption{Component distributions of $a_1$ and $a_2$ for curve $\#10$.}\label{table:components}
\end{center}
\end{table}

A similar analysis can be applied to the other disconnected groups in Table \ref{table:USp4groups}, however the definition of the sets $P_i$ varies.  For Distribution \#11 the group $H_1^{-}$ has two components, and for the curve $y^2=f(x)=x^5+x^3=x$ the correct partitioning of primes simply depends on the value of $p$ modulo 4, not on how $f(x)$ splits in $\Fp[x]$ (in fact, the set of primes where $f(x)$ splits intersects both partitions).  In other cases both a modular constraint and a splitting condition may apply.  In general there is some partitioning of primes which corresponds to a partitioning of the components of $H$, and the corresponding distributions appear to agree.

In addition to verifying the distribution of $\barLp$ over sets of primes, for each group corresponding to a split Jacobian, one can also check whether $L_p(T)$ admits a factorization of the expected type.  For the product groups we expect $L_p(T)$ to factor into two quadratics with coefficients in $\mathbb{Z}$.  For the groups $H_i^k$ we expect a factorization of the form given in (\ref{equation:splitLp}), and a similar form applies to $H_i^-$.\footnote{As previously noted, there are curves with simple Jacobians matching distributions \#2 and \#11 for which this would not apply, but they don't appear in Table \ref{table:g2a1dist}.}  For each curve in Table \ref{table:g2a1dist} with $f(x)$ of degree 5 corresponding to such a group, we verified the existence of the expected factorization for primes $p\le 10^6$.  For groups with multiple components, we partitioned the primes appropriately (note that for groups of the form $J(H)$ we do not expect a factorization for primes corresponding to components with off-diagonal elements).

The final piece of evidence we present is much less precise, and of an entirely different nature.  By computing the group structure of the Jacobian and determining the rank of its $\ell$-Sylow subgroup for many $p$, one can estimate the size of the mod $\ell$ Galois image $G[\ell]$ in $GSp(2g,\mathbb{Z}/\ell\mathbb{Z})$.  The primes for which the $\ell$-Sylow subgroup has rank $2g$ correspond to the identity element in $G[\ell]$.    The group $GSp(2g,\mathbb{Z}/\ell\mathbb{Z})$ has size $O(\ell^{11})$, and this corresponds to the real dimension of $USp(2g)$ which is 10 (the reduction in dimension arises from unitarization).  For all but the first group in Table \ref{table:USp4groups}. the real dimension is at most 6. We should expect correspondingly small $G[\ell]$, at most $O(\ell^7)$ in size.  By computing
\begin{equation}
d = \frac{\log(\#G[\ell_1])- \log(\#G[\ell_2])}{\log\ell_1 - \log\ell_2}
\end{equation}
for various $\ell_1\ne \ell_2$ we obtain a general estimate for $\#G[\ell]=O(\ell^d)$.  The value $d-1$ is then an estimate of the real dimension of $H$.  This is necessarily a rather crude approximation, and one must take care to avoid exceptional values of $\ell$.  We performed this computation for each exceptional curve in Table \ref{table:g2a1dist} where $f(x)$ has degree 5 for $p\le N=2^{24}$ and $\ell$ ranging from 3 to 19.  The results agreed with the corresponding dimensions in Table \ref{table:g2a1dist} to within $\pm 1$.

\subsection{Nonexistence Arguments} \label{section:nonexistence}

We wish to thank 
Jean-Pierre Serre
(private communication) for suggesting the following argument to rule
out the case $H = J(G_2 \times G_2)$. 
Let $E$ be the endomorphism 
algebra of the 
Jacobian of the curve over $\overline{\mathbb{Q}}$.
Then $E \otimes \mathbb{Q}$
must be a commutative $\mathbb{Q}$-algebra of rank 4; more precisely, it
must be either
$K \otimes_{\mathbb{Q}} K'$ for some nonisomorphic imaginary quadratic fields
$K, K'$, or a quartic CM-field. (We cannot have $K \cong K'$
or else $H$ would be forced to have dimension 1 rather than 2.)
In either case, $\Aut(E \otimes \mathbb{Q})$
has at most 4 elements, so the elements of $H$ are defined over
a field of degree at most 4 over $\mathbb{Q}$. However,
$J(G_2 \times G_2)$ has 8 connected
components, so this is impossible.

We do not yet have an analogous argument to rule out $H = J(H_2^6)$.
However, Serre suggests that it should be possible to give such an argument
based on the following data: the algebra $E \otimes \mathbb{Q}$ must be a
(possibly split) quaternion algebra over a quadratic
field, and the image of Galois in $\Aut(E \otimes \mathbb{Q})$ must have order 24.

\section{Conclusion}

Based on the results presented in Section \ref{section:g2Exceptions}, we now state a more explicit form of Conjecture \ref{conjecture:GeneralSatoTate} for genus 2 curves.

\begin{conj}
Let $C$ be a genus $2$ curve.  The distribution of $\barLp$ over $p\le N$ converges (as $N\to\infty)$ to the distribution of $\chi(T)$ in one of the first $23$ subgroups of $USp(4)$ listed in Table \ref{table:USp4groups}.  For almost all curves, this group is $USp(4)$.
\end{conj}


It would be interesting to carry out a similar analysis in genus 3, but it is not immediately clear from the genus 2 results how many exceptional distributions one should expect.  An exhaustive search of the type undertaken in genus 2 may not be computationally feasible in genus 3.

We end by once again thanking Nicholas Katz for his invaluable support throughout this project, and David Vogan for several helpful conversations.  We also thank Zeev Rudnick for his feedback on an early draft of this paper, and 
Jean-Pierre Serre for his remarks in \S~\ref{section:nonexistence}.

\clearpage
\section{Appendix I - Distributions of $s_k$}
This appendix is a gallery of distributions of $s_k$ for hyperelliptic curves of genus 1, 2, and 3 with large Galois image.  The $s_k$ are the $k$th power sums (Newton symmetric function) of the roots of the unitarized $L$-polynomial $\barLp$ defined in Section \ref{section:preliminaries}.

Each figure represents a histogram of approximately $\pi(N)$ values derived from $\barLp$ for $p\le N$ where the curve has good reduction.  The horizontal axis ranges from $-\binom{2g}{k}$ to $\binom{2g}{k}$, divided into approximately $\sqrt{\pi(N)}$ buckets.  The vertical axis is scaled to fit the data, with the height of the uniform distribution indicated by a dotted line.

\subsection{Genus 1}
\vspace{12pt}
$N = 2^{35};\quad y^2 = x^3 + 314159x + 271828$.\vspace{12pt}\\
\vspace{24pt}
\input{g1_distrib}
\clearpage
\subsection{Genus 2}
\vspace{12pt}
$N = 2^{26};$\\
$y^2 = x^5 + 314159x^3 + 271828x^2 + 1644934x + 57721566$.\\
\vspace{24pt}
\input{g2_distrib}
\clearpage
\subsection{Genus 3}
$N = 2^{25};$\\
$y^2 = x^7 + 314159x^5 + 271828x^4 + 1644934x^3 +57721566x^2 +1618034x + 141021.$\\
\vspace{12pt}
\input{g3_distrib}
\clearpage

\section{Appendix II - Distributions of $a_1$ in Genus 2}

This appendix contains $a_1$ histograms for each of the first 12 curves listed in Table \ref{table:g2a1dist}, computed with $N=2^{26}$.  The approximate area of the central spike is given by $z(C)$ in Table \ref{table:g2a1dist}, corresponding to primes for which $a_p = 0$.  The secondary spikes at $a_1=\pm 2$ appearing in distributions \#8 and \#12 have approximately zero area. 
\vspace{24pt}
\small
\input{g2_a1_1}

\hspace{44pt}\#1\hspace{204pt}\#2\\
\vspace{36pt}
\input{g2_a1_2}

\hspace{44pt}\#3\hspace{204pt}\#4\\
\vspace{36pt}
\input{g2_a1_3}

\hspace{44pt}\#5\hspace{204pt}\#6\\
\clearpage
\input{g2_a1_4}

\hspace{44pt}\#7\hspace{204pt}\#8\\
\vspace{36pt}
\input{g2_a1_5}

\hspace{44pt}\#9\hspace{200pt}\#10\\
\vspace{36pt}
\input{g2_a1_6}

\hspace{42pt}\#11\hspace{200pt}\#12\\
\clearpage

\bibliographystyle{amsplain}
\input{agct11.bbl}

\end{document}

%% file: converge.tex
\small
\vspace{30pt}
\input{g1_a1_12_15}

\hspace{36pt}$N=2^{12}$\hspace{183pt}$N=2^{15}$
\vspace{30pt}
\input{g1_a1_18_21}

\hspace{36pt}$N=2^{18}$\hspace{183pt}$N=2^{21}$
\vspace{30pt}
\input{g1_a1_24_27}

\hspace{36pt}$N=2^{24}$\hspace{183pt}$N=2^{27}$
\vspace{30pt}
\input{g1_a1_30_33}

\hspace{36pt}$N=2^{30}$\hspace{183pt}$N=2^{33}$
\vspace{15pt}
\normalsize
\begin{center}
Distribution of $a_1=-a_p/\sqrt{p}$ for $p \le N$ with good reduction.
$$y^2 = x^3+314159x+271828.$$
\end{center}

\pagebreak
\small
\vspace{30pt}
\input{g2_a1_12_14}

\hspace{36pt}$N=2^{12}$\hspace{183pt}$N=2^{14}$
\vspace{30pt}
\input{g2_a1_16_18}

\hspace{36pt}$N=2^{16}$\hspace{183pt}$N=2^{18}$
\vspace{30pt}
\input{g2_a1_20_22}

\hspace{36pt}$N=2^{20}$\hspace{183pt}$N=2^{22}$
\vspace{30pt}
\input{g2_a1_24_26}

\hspace{36pt}$N=2^{24}$\hspace{183pt}$N=2^{26}$
\vspace{15pt}
\normalsize
\begin{center}
Distribution of $a_1=-a_p/\sqrt{p}$ for $p \le N$ with good reduction.
$$y^2 = x^5 + 314159x^3 + 271828x^2 + 1644934x + 57721566.$$
\end{center}
\pagebreak
\small
\vspace{30pt}
\input{g3_a1_11_13}

\hspace{36pt}$N=2^{11}$\hspace{183pt}$N=2^{13}$
\vspace{30pt}
\input{g3_a1_15_17}

\hspace{36pt}$N=2^{15}$\hspace{183pt}$N=2^{17}$
\vspace{30pt}
\input{g3_a1_19_21}

\hspace{36pt}$N=2^{19}$\hspace{183pt}$N=2^{21}$
\vspace{30pt}
\input{g3_a1_23_25}

\hspace{36pt}$N=2^{23}$\hspace{183pt}$N=2^{25}$
\vspace{15pt}
\normalsize
\begin{center}
Distribution of $a_1=-a_p/\sqrt{p}$ for $p \le N$ with good reduction.
$$y^2 = x^7 + 314159x^5 + 271828x^4 + 1644934x^3 +57721566x^2 +1618034x + 141021.$$
\end{center}

%% file: g1_a1_12_15.tex
\hspace{-30pt}
\begin{picture}(2400,640)
\color[rgb]{0.0,0.2,0.8}
\z(24.00,0)(24.00,203.74)(67.43,203.74)(67.43,0)
\z(67.43,0)(67.43,259.30)(110.87,259.30)(110.87,0)
\z(110.87,0)(110.87,463.04)(154.30,463.04)(154.30,0)
\z(154.30,0)(154.30,388.96)(197.74,388.96)(197.74,0)
\z(197.74,0)(197.74,500.09)(241.17,500.09)(241.17,0)
\z(241.17,0)(241.17,518.61)(284.61,518.61)(284.61,0)
\z(284.61,0)(284.61,574.17)(328.04,574.17)(328.04,0)
\z(328.04,0)(328.04,481.57)(371.48,481.57)(371.48,0)
\z(371.48,0)(371.48,426.00)(414.91,426.00)(414.91,0)
\z(414.91,0)(414.91,555.65)(458.35,555.65)(458.35,0)
\z(458.35,0)(458.35,351.91)(501.78,351.91)(501.78,0)
\z(501.78,0)(501.78,640)(545.22,640)(545.22,0)
\z(545.22,0)(545.22,481.57)(588.65,481.57)(588.65,0)
\z(588.65,0)(588.65,555.65)(632.09,555.65)(632.09,0)
\z(632.09,0)(632.09,555.65)(675.52,555.65)(675.52,0)
\z(675.52,0)(675.52,574.17)(718.96,574.17)(718.96,0)
\z(718.96,0)(718.96,463.04)(762.39,463.04)(762.39,0)
\z(762.39,0)(762.39,640)(805.83,640)(805.83,0)
\z(805.83,0)(805.83,444.52)(849.26,444.52)(849.26,0)
\z(849.26,0)(849.26,388.96)(892.70,388.96)(892.70,0)
\z(892.70,0)(892.70,333.39)(936.13,333.39)(936.13,0)
\z(936.13,0)(936.13,222.26)(979.57,222.26)(979.57,0)
\z(979.57,0)(979.57,185.22)(1023.00,185.22)(1023.00,0)
\color{black}
\dottedline{12}(0,426.00)(1000,426.00)
\z(0,0)(0,640)(1000,640)(1000,0)(0,0)
\color[rgb]{0.0,0.2,0.8}
\z(1424.00,0)(1424.00,57.76)(1440.93,57.76)(1440.93,0)
\z(1440.93,0)(1440.93,187.73)(1457.86,187.73)(1457.86,0)
\z(1457.86,0)(1457.86,194.95)(1474.80,194.95)(1474.80,0)
\z(1474.80,0)(1474.80,252.71)(1491.73,252.71)(1491.73,0)
\z(1491.73,0)(1491.73,375.46)(1508.66,375.46)(1508.66,0)
\z(1508.66,0)(1508.66,303.25)(1525.59,303.25)(1525.59,0)
\z(1525.59,0)(1525.59,375.46)(1542.53,375.46)(1542.53,0)
\z(1542.53,0)(1542.53,353.80)(1559.46,353.80)(1559.46,0)
\z(1559.46,0)(1559.46,382.68)(1576.39,382.68)(1576.39,0)
\z(1576.39,0)(1576.39,339.36)(1593.32,339.36)(1593.32,0)
\z(1593.32,0)(1593.32,346.58)(1610.25,346.58)(1610.25,0)
\z(1610.25,0)(1610.25,440.44)(1627.19,440.44)(1627.19,0)
\z(1627.19,0)(1627.19,476.54)(1644.12,476.54)(1644.12,0)
\z(1644.12,0)(1644.12,411.56)(1661.05,411.56)(1661.05,0)
\z(1661.05,0)(1661.05,346.58)(1677.98,346.58)(1677.98,0)
\z(1677.98,0)(1677.98,555.97)(1694.92,555.97)(1694.92,0)
\z(1694.92,0)(1694.92,519.86)(1711.85,519.86)(1711.85,0)
\z(1711.85,0)(1711.85,433.22)(1728.78,433.22)(1728.78,0)
\z(1728.78,0)(1728.78,577.63)(1745.71,577.63)(1745.71,0)
\z(1745.71,0)(1745.71,483.76)(1762.64,483.76)(1762.64,0)
\z(1762.64,0)(1762.64,541.53)(1779.58,541.53)(1779.58,0)
\z(1779.58,0)(1779.58,483.76)(1796.51,483.76)(1796.51,0)
\z(1796.51,0)(1796.51,476.54)(1813.44,476.54)(1813.44,0)
\z(1813.44,0)(1813.44,483.76)(1830.37,483.76)(1830.37,0)
\z(1830.37,0)(1830.37,640)(1847.31,640)(1847.31,0)
\z(1847.31,0)(1847.31,519.86)(1864.24,519.86)(1864.24,0)
\z(1864.24,0)(1864.24,541.53)(1881.17,541.53)(1881.17,0)
\z(1881.17,0)(1881.17,613.73)(1898.10,613.73)(1898.10,0)
\z(1898.10,0)(1898.10,498.20)(1915.03,498.20)(1915.03,0)
\z(1915.03,0)(1915.03,541.53)(1931.97,541.53)(1931.97,0)
\z(1931.97,0)(1931.97,519.86)(1948.90,519.86)(1948.90,0)
\z(1948.90,0)(1948.90,548.75)(1965.83,548.75)(1965.83,0)
\z(1965.83,0)(1965.83,534.31)(1982.76,534.31)(1982.76,0)
\z(1982.76,0)(1982.76,512.64)(1999.69,512.64)(1999.69,0)
\z(1999.69,0)(1999.69,563.19)(2016.63,563.19)(2016.63,0)
\z(2016.63,0)(2016.63,635.39)(2033.56,635.39)(2033.56,0)
\z(2033.56,0)(2033.56,426.00)(2050.49,426.00)(2050.49,0)
\z(2050.49,0)(2050.49,620.95)(2067.42,620.95)(2067.42,0)
\z(2067.42,0)(2067.42,527.08)(2084.36,527.08)(2084.36,0)
\z(2084.36,0)(2084.36,411.56)(2101.29,411.56)(2101.29,0)
\z(2101.29,0)(2101.29,635.39)(2118.22,635.39)(2118.22,0)
\z(2118.22,0)(2118.22,476.54)(2135.15,476.54)(2135.15,0)
\z(2135.15,0)(2135.15,397.12)(2152.08,397.12)(2152.08,0)
\z(2152.08,0)(2152.08,548.75)(2169.02,548.75)(2169.02,0)
\z(2169.02,0)(2169.02,512.64)(2185.95,512.64)(2185.95,0)
\z(2185.95,0)(2185.95,563.19)(2202.88,563.19)(2202.88,0)
\z(2202.88,0)(2202.88,469.32)(2219.81,469.32)(2219.81,0)
\z(2219.81,0)(2219.81,440.44)(2236.75,440.44)(2236.75,0)
\z(2236.75,0)(2236.75,426.00)(2253.68,426.00)(2253.68,0)
\z(2253.68,0)(2253.68,368.24)(2270.61,368.24)(2270.61,0)
\z(2270.61,0)(2270.61,368.24)(2287.54,368.24)(2287.54,0)
\z(2287.54,0)(2287.54,418.78)(2304.47,418.78)(2304.47,0)
\z(2304.47,0)(2304.47,303.25)(2321.41,303.25)(2321.41,0)
\z(2321.41,0)(2321.41,296.03)(2338.34,296.03)(2338.34,0)
\z(2338.34,0)(2338.34,332.14)(2355.27,332.14)(2355.27,0)
\z(2355.27,0)(2355.27,223.83)(2372.20,223.83)(2372.20,0)
\z(2372.20,0)(2372.20,173.29)(2389.14,173.29)(2389.14,0)
\z(2389.14,0)(2389.14,194.95)(2406.07,194.95)(2406.07,0)
\z(2406.07,0)(2406.07,101.08)(2423.00,101.08)(2423.00,0)
\color{black}
\dottedline{12}(1400,426.00)(2400,426.00)
\z(1400,0)(1400,640)(2400,640)(2400,0)(1400,0)
\end{picture}

%% file: g1_a1_18_21.tex
\hspace{-30pt}
\begin{picture}(2400,640)
\color[rgb]{0.0,0.2,0.8}
\z(24.00,0)(24.00,45.14)(30.62,45.14)(30.62,0)
\z(30.62,0)(30.62,78.99)(37.23,78.99)(37.23,0)
\z(37.23,0)(37.23,132.60)(43.85,132.60)(43.85,0)
\z(43.85,0)(43.85,177.74)(50.46,177.74)(50.46,0)
\z(50.46,0)(50.46,197.48)(57.08,197.48)(57.08,0)
\z(57.08,0)(57.08,186.20)(63.70,186.20)(63.70,0)
\z(63.70,0)(63.70,253.91)(70.31,253.91)(70.31,0)
\z(70.31,0)(70.31,211.59)(76.93,211.59)(76.93,0)
\z(76.93,0)(76.93,228.52)(83.54,228.52)(83.54,0)
\z(83.54,0)(83.54,239.80)(90.16,239.80)(90.16,0)
\z(90.16,0)(90.16,287.76)(96.77,287.76)(96.77,0)
\z(96.77,0)(96.77,287.76)(103.39,287.76)(103.39,0)
\z(103.39,0)(103.39,296.23)(110.01,296.23)(110.01,0)
\z(110.01,0)(110.01,321.62)(116.62,321.62)(116.62,0)
\z(116.62,0)(116.62,358.29)(123.24,358.29)(123.24,0)
\z(123.24,0)(123.24,287.76)(129.85,287.76)(129.85,0)
\z(129.85,0)(129.85,378.04)(136.47,378.04)(136.47,0)
\z(136.47,0)(136.47,363.93)(143.09,363.93)(143.09,0)
\z(143.09,0)(143.09,335.72)(149.70,335.72)(149.70,0)
\z(149.70,0)(149.70,344.19)(156.32,344.19)(156.32,0)
\z(156.32,0)(156.32,318.79)(162.93,318.79)(162.93,0)
\z(162.93,0)(162.93,378.04)(169.55,378.04)(169.55,0)
\z(169.55,0)(169.55,335.72)(176.17,335.72)(176.17,0)
\z(176.17,0)(176.17,403.43)(182.78,403.43)(182.78,0)
\z(182.78,0)(182.78,420.36)(189.40,420.36)(189.40,0)
\z(189.40,0)(189.40,417.54)(196.01,417.54)(196.01,0)
\z(196.01,0)(196.01,375.22)(202.63,375.22)(202.63,0)
\z(202.63,0)(202.63,451.39)(209.25,451.39)(209.25,0)
\z(209.25,0)(209.25,406.25)(215.86,406.25)(215.86,0)
\z(215.86,0)(215.86,454.21)(222.48,454.21)(222.48,0)
\z(222.48,0)(222.48,420.36)(229.09,420.36)(229.09,0)
\z(229.09,0)(229.09,389.32)(235.71,389.32)(235.71,0)
\z(235.71,0)(235.71,527.56)(242.32,527.56)(242.32,0)
\z(242.32,0)(242.32,440.11)(248.94,440.11)(248.94,0)
\z(248.94,0)(248.94,431.64)(255.56,431.64)(255.56,0)
\z(255.56,0)(255.56,479.60)(262.17,479.60)(262.17,0)
\z(262.17,0)(262.17,409.07)(268.79,409.07)(268.79,0)
\z(268.79,0)(268.79,465.50)(275.40,465.50)(275.40,0)
\z(275.40,0)(275.40,442.93)(282.02,442.93)(282.02,0)
\z(282.02,0)(282.02,502.17)(288.64,502.17)(288.64,0)
\z(288.64,0)(288.64,541.67)(295.25,541.67)(295.25,0)
\z(295.25,0)(295.25,490.89)(301.87,490.89)(301.87,0)
\z(301.87,0)(301.87,524.74)(308.48,524.74)(308.48,0)
\z(308.48,0)(308.48,457.03)(315.10,457.03)(315.10,0)
\z(315.10,0)(315.10,502.17)(321.72,502.17)(321.72,0)
\z(321.72,0)(321.72,451.39)(328.33,451.39)(328.33,0)
\z(328.33,0)(328.33,567.06)(334.95,567.06)(334.95,0)
\z(334.95,0)(334.95,519.10)(341.56,519.10)(341.56,0)
\z(341.56,0)(341.56,510.64)(348.18,510.64)(348.18,0)
\z(348.18,0)(348.18,493.71)(354.79,493.71)(354.79,0)
\z(354.79,0)(354.79,448.57)(361.41,448.57)(361.41,0)
\z(361.41,0)(361.41,493.71)(368.03,493.71)(368.03,0)
\z(368.03,0)(368.03,485.25)(374.64,485.25)(374.64,0)
\z(374.64,0)(374.64,507.81)(381.26,507.81)(381.26,0)
\z(381.26,0)(381.26,504.99)(387.87,504.99)(387.87,0)
\z(387.87,0)(387.87,552.95)(394.49,552.95)(394.49,0)
\z(394.49,0)(394.49,598.09)(401.11,598.09)(401.11,0)
\z(401.11,0)(401.11,567.06)(407.72,567.06)(407.72,0)
\z(407.72,0)(407.72,595.27)(414.34,595.27)(414.34,0)
\z(414.34,0)(414.34,561.42)(420.95,561.42)(420.95,0)
\z(420.95,0)(420.95,524.74)(427.57,524.74)(427.57,0)
\z(427.57,0)(427.57,544.49)(434.19,544.49)(434.19,0)
\z(434.19,0)(434.19,589.63)(440.80,589.63)(440.80,0)
\z(440.80,0)(440.80,600.91)(447.42,600.91)(447.42,0)
\z(447.42,0)(447.42,468.32)(454.03,468.32)(454.03,0)
\z(454.03,0)(454.03,524.74)(460.65,524.74)(460.65,0)
\z(460.65,0)(460.65,536.03)(467.26,536.03)(467.26,0)
\z(467.26,0)(467.26,544.49)(473.88,544.49)(473.88,0)
\z(473.88,0)(473.88,536.03)(480.50,536.03)(480.50,0)
\z(480.50,0)(480.50,516.28)(487.11,516.28)(487.11,0)
\z(487.11,0)(487.11,572.70)(493.73,572.70)(493.73,0)
\z(493.73,0)(493.73,499.35)(500.34,499.35)(500.34,0)
\z(500.34,0)(500.34,442.93)(506.96,442.93)(506.96,0)
\z(506.96,0)(506.96,626.30)(513.58,626.30)(513.58,0)
\z(513.58,0)(513.58,527.56)(520.19,527.56)(520.19,0)
\z(520.19,0)(520.19,626.30)(526.81,626.30)(526.81,0)
\z(526.81,0)(526.81,541.67)(533.42,541.67)(533.42,0)
\z(533.42,0)(533.42,575.52)(540.04,575.52)(540.04,0)
\z(540.04,0)(540.04,524.74)(546.66,524.74)(546.66,0)
\z(546.66,0)(546.66,550.13)(553.27,550.13)(553.27,0)
\z(553.27,0)(553.27,564.24)(559.89,564.24)(559.89,0)
\z(559.89,0)(559.89,538.85)(566.50,538.85)(566.50,0)
\z(566.50,0)(566.50,524.74)(573.12,524.74)(573.12,0)
\z(573.12,0)(573.12,589.63)(579.74,589.63)(579.74,0)
\z(579.74,0)(579.74,569.88)(586.35,569.88)(586.35,0)
\z(586.35,0)(586.35,623.48)(592.97,623.48)(592.97,0)
\z(592.97,0)(592.97,502.17)(599.58,502.17)(599.58,0)
\z(599.58,0)(599.58,527.56)(606.20,527.56)(606.20,0)
\z(606.20,0)(606.20,516.28)(612.81,516.28)(612.81,0)
\z(612.81,0)(612.81,595.27)(619.43,595.27)(619.43,0)
\z(619.43,0)(619.43,561.42)(626.05,561.42)(626.05,0)
\z(626.05,0)(626.05,510.64)(632.66,510.64)(632.66,0)
\z(632.66,0)(632.66,454.21)(639.28,454.21)(639.28,0)
\z(639.28,0)(639.28,547.31)(645.89,547.31)(645.89,0)
\z(645.89,0)(645.89,536.03)(652.51,536.03)(652.51,0)
\z(652.51,0)(652.51,536.03)(659.13,536.03)(659.13,0)
\z(659.13,0)(659.13,541.67)(665.74,541.67)(665.74,0)
\z(665.74,0)(665.74,555.77)(672.36,555.77)(672.36,0)
\z(672.36,0)(672.36,499.35)(678.97,499.35)(678.97,0)
\z(678.97,0)(678.97,558.60)(685.59,558.60)(685.59,0)
\z(685.59,0)(685.59,547.31)(692.21,547.31)(692.21,0)
\z(692.21,0)(692.21,440.11)(698.82,440.11)(698.82,0)
\z(698.82,0)(698.82,516.28)(705.44,516.28)(705.44,0)
\z(705.44,0)(705.44,567.06)(712.05,567.06)(712.05,0)
\z(712.05,0)(712.05,504.99)(718.67,504.99)(718.67,0)
\z(718.67,0)(718.67,479.60)(725.28,479.60)(725.28,0)
\z(725.28,0)(725.28,476.78)(731.90,476.78)(731.90,0)
\z(731.90,0)(731.90,507.81)(738.52,507.81)(738.52,0)
\z(738.52,0)(738.52,454.21)(745.13,454.21)(745.13,0)
\z(745.13,0)(745.13,437.28)(751.75,437.28)(751.75,0)
\z(751.75,0)(751.75,485.25)(758.36,485.25)(758.36,0)
\z(758.36,0)(758.36,527.56)(764.98,527.56)(764.98,0)
\z(764.98,0)(764.98,555.77)(771.60,555.77)(771.60,0)
\z(771.60,0)(771.60,507.81)(778.21,507.81)(778.21,0)
\z(778.21,0)(778.21,468.32)(784.83,468.32)(784.83,0)
\z(784.83,0)(784.83,423.18)(791.44,423.18)(791.44,0)
\z(791.44,0)(791.44,555.77)(798.06,555.77)(798.06,0)
\z(798.06,0)(798.06,457.03)(804.68,457.03)(804.68,0)
\z(804.68,0)(804.68,383.68)(811.29,383.68)(811.29,0)
\z(811.29,0)(811.29,507.81)(817.91,507.81)(817.91,0)
\z(817.91,0)(817.91,476.78)(824.52,476.78)(824.52,0)
\z(824.52,0)(824.52,431.64)(831.14,431.64)(831.14,0)
\z(831.14,0)(831.14,431.64)(837.75,431.64)(837.75,0)
\z(837.75,0)(837.75,392.15)(844.37,392.15)(844.37,0)
\z(844.37,0)(844.37,442.93)(850.99,442.93)(850.99,0)
\z(850.99,0)(850.99,397.79)(857.60,397.79)(857.60,0)
\z(857.60,0)(857.60,397.79)(864.22,397.79)(864.22,0)
\z(864.22,0)(864.22,372.40)(870.83,372.40)(870.83,0)
\z(870.83,0)(870.83,400.61)(877.45,400.61)(877.45,0)
\z(877.45,0)(877.45,411.89)(884.07,411.89)(884.07,0)
\z(884.07,0)(884.07,327.26)(890.68,327.26)(890.68,0)
\z(890.68,0)(890.68,301.87)(897.30,301.87)(897.30,0)
\z(897.30,0)(897.30,411.89)(903.91,411.89)(903.91,0)
\z(903.91,0)(903.91,304.69)(910.53,304.69)(910.53,0)
\z(910.53,0)(910.53,299.05)(917.15,299.05)(917.15,0)
\z(917.15,0)(917.15,358.29)(923.76,358.29)(923.76,0)
\z(923.76,0)(923.76,276.48)(930.38,276.48)(930.38,0)
\z(930.38,0)(930.38,307.51)(936.99,307.51)(936.99,0)
\z(936.99,0)(936.99,335.72)(943.61,335.72)(943.61,0)
\z(943.61,0)(943.61,304.69)(950.23,304.69)(950.23,0)
\z(950.23,0)(950.23,290.58)(956.84,290.58)(956.84,0)
\z(956.84,0)(956.84,234.16)(963.46,234.16)(963.46,0)
\z(963.46,0)(963.46,270.83)(970.07,270.83)(970.07,0)
\z(970.07,0)(970.07,217.23)(976.69,217.23)(976.69,0)
\z(976.69,0)(976.69,231.34)(983.30,231.34)(983.30,0)
\z(983.30,0)(983.30,200.30)(989.92,200.30)(989.92,0)
\z(989.92,0)(989.92,239.80)(996.54,239.80)(996.54,0)
\z(996.54,0)(996.54,141.06)(1003.15,141.06)(1003.15,0)
\z(1003.15,0)(1003.15,118.49)(1009.77,118.49)(1009.77,0)
\z(1009.77,0)(1009.77,118.49)(1016.38,118.49)(1016.38,0)
\z(1016.38,0)(1016.38,47.96)(1023.00,47.96)(1023.00,0)
\color{black}
\dottedline{12}(0,426.00)(1000,426.00)
\z(0,0)(0,640)(1000,640)(1000,0)(0,0)
\color[rgb]{0.0,0.2,0.8}
\z(1424.00,0)(1424.00,28.11)(1426.54,28.11)(1426.54,0)
\z(1426.54,0)(1426.54,55.14)(1429.07,55.14)(1429.07,0)
\z(1429.07,0)(1429.07,91.90)(1431.61,91.90)(1431.61,0)
\z(1431.61,0)(1431.61,101.63)(1434.14,101.63)(1434.14,0)
\z(1434.14,0)(1434.14,108.12)(1436.68,108.12)(1436.68,0)
\z(1436.68,0)(1436.68,110.28)(1439.21,110.28)(1439.21,0)
\z(1439.21,0)(1439.21,134.07)(1441.75,134.07)(1441.75,0)
\z(1441.75,0)(1441.75,158.94)(1444.28,158.94)(1444.28,0)
\z(1444.28,0)(1444.28,154.61)(1446.82,154.61)(1446.82,0)
\z(1446.82,0)(1446.82,177.32)(1449.36,177.32)(1449.36,0)
\z(1449.36,0)(1449.36,184.89)(1451.89,184.89)(1451.89,0)
\z(1451.89,0)(1451.89,188.13)(1454.43,188.13)(1454.43,0)
\z(1454.43,0)(1454.43,217.32)(1456.96,217.32)(1456.96,0)
\z(1456.96,0)(1456.96,168.67)(1459.50,168.67)(1459.50,0)
\z(1459.50,0)(1459.50,197.86)(1462.03,197.86)(1462.03,0)
\z(1462.03,0)(1462.03,233.54)(1464.57,233.54)(1464.57,0)
\z(1464.57,0)(1464.57,213.00)(1467.10,213.00)(1467.10,0)
\z(1467.10,0)(1467.10,242.19)(1469.64,242.19)(1469.64,0)
\z(1469.64,0)(1469.64,214.08)(1472.18,214.08)(1472.18,0)
\z(1472.18,0)(1472.18,261.65)(1474.71,261.65)(1474.71,0)
\z(1474.71,0)(1474.71,225.97)(1477.25,225.97)(1477.25,0)
\z(1477.25,0)(1477.25,249.76)(1479.78,249.76)(1479.78,0)
\z(1479.78,0)(1479.78,248.68)(1482.32,248.68)(1482.32,0)
\z(1482.32,0)(1482.32,275.71)(1484.85,275.71)(1484.85,0)
\z(1484.85,0)(1484.85,269.22)(1487.39,269.22)(1487.39,0)
\z(1487.39,0)(1487.39,291.93)(1489.92,291.93)(1489.92,0)
\z(1489.92,0)(1489.92,277.87)(1492.46,277.87)(1492.46,0)
\z(1492.46,0)(1492.46,281.12)(1494.99,281.12)(1494.99,0)
\z(1494.99,0)(1494.99,295.17)(1497.53,295.17)(1497.53,0)
\z(1497.53,0)(1497.53,308.15)(1500.07,308.15)(1500.07,0)
\z(1500.07,0)(1500.07,272.47)(1502.60,272.47)(1502.60,0)
\z(1502.60,0)(1502.60,329.77)(1505.14,329.77)(1505.14,0)
\z(1505.14,0)(1505.14,274.63)(1507.67,274.63)(1507.67,0)
\z(1507.67,0)(1507.67,317.88)(1510.21,317.88)(1510.21,0)
\z(1510.21,0)(1510.21,331.93)(1512.74,331.93)(1512.74,0)
\z(1512.74,0)(1512.74,320.04)(1515.28,320.04)(1515.28,0)
\z(1515.28,0)(1515.28,286.52)(1517.81,286.52)(1517.81,0)
\z(1517.81,0)(1517.81,315.72)(1520.35,315.72)(1520.35,0)
\z(1520.35,0)(1520.35,328.69)(1522.89,328.69)(1522.89,0)
\z(1522.89,0)(1522.89,323.28)(1525.42,323.28)(1525.42,0)
\z(1525.42,0)(1525.42,320.04)(1527.96,320.04)(1527.96,0)
\z(1527.96,0)(1527.96,345.99)(1530.49,345.99)(1530.49,0)
\z(1530.49,0)(1530.49,336.26)(1533.03,336.26)(1533.03,0)
\z(1533.03,0)(1533.03,344.91)(1535.56,344.91)(1535.56,0)
\z(1535.56,0)(1535.56,363.29)(1538.10,363.29)(1538.10,0)
\z(1538.10,0)(1538.10,342.75)(1540.63,342.75)(1540.63,0)
\z(1540.63,0)(1540.63,363.29)(1543.17,363.29)(1543.17,0)
\z(1543.17,0)(1543.17,335.18)(1545.71,335.18)(1545.71,0)
\z(1545.71,0)(1545.71,334.10)(1548.24,334.10)(1548.24,0)
\z(1548.24,0)(1548.24,390.32)(1550.78,390.32)(1550.78,0)
\z(1550.78,0)(1550.78,348.15)(1553.31,348.15)(1553.31,0)
\z(1553.31,0)(1553.31,347.07)(1555.85,347.07)(1555.85,0)
\z(1555.85,0)(1555.85,408.70)(1558.38,408.70)(1558.38,0)
\z(1558.38,0)(1558.38,362.21)(1560.92,362.21)(1560.92,0)
\z(1560.92,0)(1560.92,349.23)(1563.45,349.23)(1563.45,0)
\z(1563.45,0)(1563.45,390.32)(1565.99,390.32)(1565.99,0)
\z(1565.99,0)(1565.99,374.10)(1568.53,374.10)(1568.53,0)
\z(1568.53,0)(1568.53,363.29)(1571.06,363.29)(1571.06,0)
\z(1571.06,0)(1571.06,367.61)(1573.60,367.61)(1573.60,0)
\z(1573.60,0)(1573.60,358.96)(1576.13,358.96)(1576.13,0)
\z(1576.13,0)(1576.13,448.71)(1578.67,448.71)(1578.67,0)
\z(1578.67,0)(1578.67,382.75)(1581.20,382.75)(1581.20,0)
\z(1581.20,0)(1581.20,407.62)(1583.74,407.62)(1583.74,0)
\z(1583.74,0)(1583.74,392.48)(1586.27,392.48)(1586.27,0)
\z(1586.27,0)(1586.27,427.08)(1588.81,427.08)(1588.81,0)
\z(1588.81,0)(1588.81,433.57)(1591.35,433.57)(1591.35,0)
\z(1591.35,0)(1591.35,415.19)(1593.88,415.19)(1593.88,0)
\z(1593.88,0)(1593.88,404.38)(1596.42,404.38)(1596.42,0)
\z(1596.42,0)(1596.42,375.18)(1598.95,375.18)(1598.95,0)
\z(1598.95,0)(1598.95,383.83)(1601.49,383.83)(1601.49,0)
\z(1601.49,0)(1601.49,423.84)(1604.02,423.84)(1604.02,0)
\z(1604.02,0)(1604.02,408.70)(1606.56,408.70)(1606.56,0)
\z(1606.56,0)(1606.56,388.16)(1609.09,388.16)(1609.09,0)
\z(1609.09,0)(1609.09,420.59)(1611.63,420.59)(1611.63,0)
\z(1611.63,0)(1611.63,413.03)(1614.16,413.03)(1614.16,0)
\z(1614.16,0)(1614.16,404.38)(1616.70,404.38)(1616.70,0)
\z(1616.70,0)(1616.70,414.11)(1619.24,414.11)(1619.24,0)
\z(1619.24,0)(1619.24,461.68)(1621.77,461.68)(1621.77,0)
\z(1621.77,0)(1621.77,438.97)(1624.31,438.97)(1624.31,0)
\z(1624.31,0)(1624.31,445.46)(1626.84,445.46)(1626.84,0)
\z(1626.84,0)(1626.84,435.73)(1629.38,435.73)(1629.38,0)
\z(1629.38,0)(1629.38,494.12)(1631.91,494.12)(1631.91,0)
\z(1631.91,0)(1631.91,411.94)(1634.45,411.94)(1634.45,0)
\z(1634.45,0)(1634.45,476.82)(1636.98,476.82)(1636.98,0)
\z(1636.98,0)(1636.98,454.11)(1639.52,454.11)(1639.52,0)
\z(1639.52,0)(1639.52,450.87)(1642.06,450.87)(1642.06,0)
\z(1642.06,0)(1642.06,495.20)(1644.59,495.20)(1644.59,0)
\z(1644.59,0)(1644.59,469.25)(1647.13,469.25)(1647.13,0)
\z(1647.13,0)(1647.13,469.25)(1649.66,469.25)(1649.66,0)
\z(1649.66,0)(1649.66,493.04)(1652.20,493.04)(1652.20,0)
\z(1652.20,0)(1652.20,447.62)(1654.73,447.62)(1654.73,0)
\z(1654.73,0)(1654.73,441.14)(1657.27,441.14)(1657.27,0)
\z(1657.27,0)(1657.27,464.92)(1659.80,464.92)(1659.80,0)
\z(1659.80,0)(1659.80,493.04)(1662.34,493.04)(1662.34,0)
\z(1662.34,0)(1662.34,482.22)(1664.88,482.22)(1664.88,0)
\z(1664.88,0)(1664.88,480.06)(1667.41,480.06)(1667.41,0)
\z(1667.41,0)(1667.41,434.65)(1669.95,434.65)(1669.95,0)
\z(1669.95,0)(1669.95,474.65)(1672.48,474.65)(1672.48,0)
\z(1672.48,0)(1672.48,426.00)(1675.02,426.00)(1675.02,0)
\z(1675.02,0)(1675.02,470.33)(1677.55,470.33)(1677.55,0)
\z(1677.55,0)(1677.55,470.33)(1680.09,470.33)(1680.09,0)
\z(1680.09,0)(1680.09,475.74)(1682.62,475.74)(1682.62,0)
\z(1682.62,0)(1682.62,517.90)(1685.16,517.90)(1685.16,0)
\z(1685.16,0)(1685.16,486.55)(1687.70,486.55)(1687.70,0)
\z(1687.70,0)(1687.70,489.79)(1690.23,489.79)(1690.23,0)
\z(1690.23,0)(1690.23,464.92)(1692.77,464.92)(1692.77,0)
\z(1692.77,0)(1692.77,487.63)(1695.30,487.63)(1695.30,0)
\z(1695.30,0)(1695.30,451.95)(1697.84,451.95)(1697.84,0)
\z(1697.84,0)(1697.84,462.76)(1700.37,462.76)(1700.37,0)
\z(1700.37,0)(1700.37,475.74)(1702.91,475.74)(1702.91,0)
\z(1702.91,0)(1702.91,520.07)(1705.44,520.07)(1705.44,0)
\z(1705.44,0)(1705.44,473.57)(1707.98,473.57)(1707.98,0)
\z(1707.98,0)(1707.98,522.23)(1710.52,522.23)(1710.52,0)
\z(1710.52,0)(1710.52,498.44)(1713.05,498.44)(1713.05,0)
\z(1713.05,0)(1713.05,483.30)(1715.59,483.30)(1715.59,0)
\z(1715.59,0)(1715.59,518.98)(1718.12,518.98)(1718.12,0)
\z(1718.12,0)(1718.12,485.47)(1720.66,485.47)(1720.66,0)
\z(1720.66,0)(1720.66,489.79)(1723.19,489.79)(1723.19,0)
\z(1723.19,0)(1723.19,478.98)(1725.73,478.98)(1725.73,0)
\z(1725.73,0)(1725.73,473.57)(1728.26,473.57)(1728.26,0)
\z(1728.26,0)(1728.26,470.33)(1730.80,470.33)(1730.80,0)
\z(1730.80,0)(1730.80,506.01)(1733.34,506.01)(1733.34,0)
\z(1733.34,0)(1733.34,481.14)(1735.87,481.14)(1735.87,0)
\z(1735.87,0)(1735.87,491.95)(1738.41,491.95)(1738.41,0)
\z(1738.41,0)(1738.41,510.34)(1740.94,510.34)(1740.94,0)
\z(1740.94,0)(1740.94,468.17)(1743.48,468.17)(1743.48,0)
\z(1743.48,0)(1743.48,496.28)(1746.01,496.28)(1746.01,0)
\z(1746.01,0)(1746.01,457.36)(1748.55,457.36)(1748.55,0)
\z(1748.55,0)(1748.55,509.25)(1751.08,509.25)(1751.08,0)
\z(1751.08,0)(1751.08,538.45)(1753.62,538.45)(1753.62,0)
\z(1753.62,0)(1753.62,521.15)(1756.15,521.15)(1756.15,0)
\z(1756.15,0)(1756.15,512.50)(1758.69,512.50)(1758.69,0)
\z(1758.69,0)(1758.69,500.60)(1761.23,500.60)(1761.23,0)
\z(1761.23,0)(1761.23,501.69)(1763.76,501.69)(1763.76,0)
\z(1763.76,0)(1763.76,544.93)(1766.30,544.93)(1766.30,0)
\z(1766.30,0)(1766.30,496.28)(1768.83,496.28)(1768.83,0)
\z(1768.83,0)(1768.83,489.79)(1771.37,489.79)(1771.37,0)
\z(1771.37,0)(1771.37,504.93)(1773.90,504.93)(1773.90,0)
\z(1773.90,0)(1773.90,541.69)(1776.44,541.69)(1776.44,0)
\z(1776.44,0)(1776.44,528.72)(1778.97,528.72)(1778.97,0)
\z(1778.97,0)(1778.97,500.60)(1781.51,500.60)(1781.51,0)
\z(1781.51,0)(1781.51,485.47)(1784.05,485.47)(1784.05,0)
\z(1784.05,0)(1784.05,516.82)(1786.58,516.82)(1786.58,0)
\z(1786.58,0)(1786.58,491.95)(1789.12,491.95)(1789.12,0)
\z(1789.12,0)(1789.12,530.88)(1791.65,530.88)(1791.65,0)
\z(1791.65,0)(1791.65,533.04)(1794.19,533.04)(1794.19,0)
\z(1794.19,0)(1794.19,541.69)(1796.72,541.69)(1796.72,0)
\z(1796.72,0)(1796.72,548.18)(1799.26,548.18)(1799.26,0)
\z(1799.26,0)(1799.26,525.47)(1801.79,525.47)(1801.79,0)
\z(1801.79,0)(1801.79,496.28)(1804.33,496.28)(1804.33,0)
\z(1804.33,0)(1804.33,567.64)(1806.87,567.64)(1806.87,0)
\z(1806.87,0)(1806.87,571.96)(1809.40,571.96)(1809.40,0)
\z(1809.40,0)(1809.40,523.31)(1811.94,523.31)(1811.94,0)
\z(1811.94,0)(1811.94,544.93)(1814.47,544.93)(1814.47,0)
\z(1814.47,0)(1814.47,550.34)(1817.01,550.34)(1817.01,0)
\z(1817.01,0)(1817.01,552.50)(1819.54,552.50)(1819.54,0)
\z(1819.54,0)(1819.54,526.55)(1822.08,526.55)(1822.08,0)
\z(1822.08,0)(1822.08,504.93)(1824.61,504.93)(1824.61,0)
\z(1824.61,0)(1824.61,544.93)(1827.15,544.93)(1827.15,0)
\z(1827.15,0)(1827.15,507.09)(1829.69,507.09)(1829.69,0)
\z(1829.69,0)(1829.69,544.93)(1832.22,544.93)(1832.22,0)
\z(1832.22,0)(1832.22,534.12)(1834.76,534.12)(1834.76,0)
\z(1834.76,0)(1834.76,569.80)(1837.29,569.80)(1837.29,0)
\z(1837.29,0)(1837.29,531.96)(1839.83,531.96)(1839.83,0)
\z(1839.83,0)(1839.83,522.23)(1842.36,522.23)(1842.36,0)
\z(1842.36,0)(1842.36,554.66)(1844.90,554.66)(1844.90,0)
\z(1844.90,0)(1844.90,537.37)(1847.43,537.37)(1847.43,0)
\z(1847.43,0)(1847.43,568.72)(1849.97,568.72)(1849.97,0)
\z(1849.97,0)(1849.97,504.93)(1852.51,504.93)(1852.51,0)
\z(1852.51,0)(1852.51,557.91)(1855.04,557.91)(1855.04,0)
\z(1855.04,0)(1855.04,520.07)(1857.58,520.07)(1857.58,0)
\z(1857.58,0)(1857.58,549.26)(1860.11,549.26)(1860.11,0)
\z(1860.11,0)(1860.11,563.31)(1862.65,563.31)(1862.65,0)
\z(1862.65,0)(1862.65,538.45)(1865.18,538.45)(1865.18,0)
\z(1865.18,0)(1865.18,555.75)(1867.72,555.75)(1867.72,0)
\z(1867.72,0)(1867.72,530.88)(1870.25,530.88)(1870.25,0)
\z(1870.25,0)(1870.25,502.77)(1872.79,502.77)(1872.79,0)
\z(1872.79,0)(1872.79,536.28)(1875.32,536.28)(1875.32,0)
\z(1875.32,0)(1875.32,541.69)(1877.86,541.69)(1877.86,0)
\z(1877.86,0)(1877.86,525.47)(1880.40,525.47)(1880.40,0)
\z(1880.40,0)(1880.40,509.25)(1882.93,509.25)(1882.93,0)
\z(1882.93,0)(1882.93,570.88)(1885.47,570.88)(1885.47,0)
\z(1885.47,0)(1885.47,563.31)(1888.00,563.31)(1888.00,0)
\z(1888.00,0)(1888.00,557.91)(1890.54,557.91)(1890.54,0)
\z(1890.54,0)(1890.54,525.47)(1893.07,525.47)(1893.07,0)
\z(1893.07,0)(1893.07,541.69)(1895.61,541.69)(1895.61,0)
\z(1895.61,0)(1895.61,511.42)(1898.14,511.42)(1898.14,0)
\z(1898.14,0)(1898.14,528.72)(1900.68,528.72)(1900.68,0)
\z(1900.68,0)(1900.68,521.15)(1903.22,521.15)(1903.22,0)
\z(1903.22,0)(1903.22,523.31)(1905.75,523.31)(1905.75,0)
\z(1905.75,0)(1905.75,511.42)(1908.29,511.42)(1908.29,0)
\z(1908.29,0)(1908.29,547.10)(1910.82,547.10)(1910.82,0)
\z(1910.82,0)(1910.82,555.75)(1913.36,555.75)(1913.36,0)
\z(1913.36,0)(1913.36,560.07)(1915.89,560.07)(1915.89,0)
\z(1915.89,0)(1915.89,536.28)(1918.43,536.28)(1918.43,0)
\z(1918.43,0)(1918.43,539.53)(1920.96,539.53)(1920.96,0)
\z(1920.96,0)(1920.96,490.87)(1923.50,490.87)(1923.50,0)
\z(1923.50,0)(1923.50,640)(1926.04,640)(1926.04,0)
\z(1926.04,0)(1926.04,536.28)(1928.57,536.28)(1928.57,0)
\z(1928.57,0)(1928.57,502.77)(1931.11,502.77)(1931.11,0)
\z(1931.11,0)(1931.11,576.29)(1933.64,576.29)(1933.64,0)
\z(1933.64,0)(1933.64,558.99)(1936.18,558.99)(1936.18,0)
\z(1936.18,0)(1936.18,540.61)(1938.71,540.61)(1938.71,0)
\z(1938.71,0)(1938.71,520.07)(1941.25,520.07)(1941.25,0)
\z(1941.25,0)(1941.25,543.85)(1943.78,543.85)(1943.78,0)
\z(1943.78,0)(1943.78,537.37)(1946.32,537.37)(1946.32,0)
\z(1946.32,0)(1946.32,512.50)(1948.86,512.50)(1948.86,0)
\z(1948.86,0)(1948.86,531.96)(1951.39,531.96)(1951.39,0)
\z(1951.39,0)(1951.39,535.20)(1953.93,535.20)(1953.93,0)
\z(1953.93,0)(1953.93,556.83)(1956.46,556.83)(1956.46,0)
\z(1956.46,0)(1956.46,509.25)(1959.00,509.25)(1959.00,0)
\z(1959.00,0)(1959.00,548.18)(1961.53,548.18)(1961.53,0)
\z(1961.53,0)(1961.53,570.88)(1964.07,570.88)(1964.07,0)
\z(1964.07,0)(1964.07,529.80)(1966.60,529.80)(1966.60,0)
\z(1966.60,0)(1966.60,525.47)(1969.14,525.47)(1969.14,0)
\z(1969.14,0)(1969.14,533.04)(1971.68,533.04)(1971.68,0)
\z(1971.68,0)(1971.68,543.85)(1974.21,543.85)(1974.21,0)
\z(1974.21,0)(1974.21,562.23)(1976.75,562.23)(1976.75,0)
\z(1976.75,0)(1976.75,567.64)(1979.28,567.64)(1979.28,0)
\z(1979.28,0)(1979.28,587.10)(1981.82,587.10)(1981.82,0)
\z(1981.82,0)(1981.82,528.72)(1984.35,528.72)(1984.35,0)
\z(1984.35,0)(1984.35,560.07)(1986.89,560.07)(1986.89,0)
\z(1986.89,0)(1986.89,550.34)(1989.42,550.34)(1989.42,0)
\z(1989.42,0)(1989.42,488.71)(1991.96,488.71)(1991.96,0)
\z(1991.96,0)(1991.96,504.93)(1994.49,504.93)(1994.49,0)
\z(1994.49,0)(1994.49,527.63)(1997.03,527.63)(1997.03,0)
\z(1997.03,0)(1997.03,496.28)(1999.57,496.28)(1999.57,0)
\z(1999.57,0)(1999.57,541.69)(2002.10,541.69)(2002.10,0)
\z(2002.10,0)(2002.10,521.15)(2004.64,521.15)(2004.64,0)
\z(2004.64,0)(2004.64,533.04)(2007.17,533.04)(2007.17,0)
\z(2007.17,0)(2007.17,542.77)(2009.71,542.77)(2009.71,0)
\z(2009.71,0)(2009.71,503.85)(2012.24,503.85)(2012.24,0)
\z(2012.24,0)(2012.24,539.53)(2014.78,539.53)(2014.78,0)
\z(2014.78,0)(2014.78,530.88)(2017.31,530.88)(2017.31,0)
\z(2017.31,0)(2017.31,528.72)(2019.85,528.72)(2019.85,0)
\z(2019.85,0)(2019.85,557.91)(2022.39,557.91)(2022.39,0)
\z(2022.39,0)(2022.39,563.31)(2024.92,563.31)(2024.92,0)
\z(2024.92,0)(2024.92,526.55)(2027.46,526.55)(2027.46,0)
\z(2027.46,0)(2027.46,549.26)(2029.99,549.26)(2029.99,0)
\z(2029.99,0)(2029.99,544.93)(2032.53,544.93)(2032.53,0)
\z(2032.53,0)(2032.53,511.42)(2035.06,511.42)(2035.06,0)
\z(2035.06,0)(2035.06,506.01)(2037.60,506.01)(2037.60,0)
\z(2037.60,0)(2037.60,553.58)(2040.13,553.58)(2040.13,0)
\z(2040.13,0)(2040.13,553.58)(2042.67,553.58)(2042.67,0)
\z(2042.67,0)(2042.67,539.53)(2045.21,539.53)(2045.21,0)
\z(2045.21,0)(2045.21,531.96)(2047.74,531.96)(2047.74,0)
\z(2047.74,0)(2047.74,456.27)(2050.28,456.27)(2050.28,0)
\z(2050.28,0)(2050.28,560.07)(2052.81,560.07)(2052.81,0)
\z(2052.81,0)(2052.81,542.77)(2055.35,542.77)(2055.35,0)
\z(2055.35,0)(2055.35,570.88)(2057.88,570.88)(2057.88,0)
\z(2057.88,0)(2057.88,531.96)(2060.42,531.96)(2060.42,0)
\z(2060.42,0)(2060.42,494.12)(2062.95,494.12)(2062.95,0)
\z(2062.95,0)(2062.95,567.64)(2065.49,567.64)(2065.49,0)
\z(2065.49,0)(2065.49,509.25)(2068.03,509.25)(2068.03,0)
\z(2068.03,0)(2068.03,526.55)(2070.56,526.55)(2070.56,0)
\z(2070.56,0)(2070.56,535.20)(2073.10,535.20)(2073.10,0)
\z(2073.10,0)(2073.10,543.85)(2075.63,543.85)(2075.63,0)
\z(2075.63,0)(2075.63,531.96)(2078.17,531.96)(2078.17,0)
\z(2078.17,0)(2078.17,487.63)(2080.70,487.63)(2080.70,0)
\z(2080.70,0)(2080.70,506.01)(2083.24,506.01)(2083.24,0)
\z(2083.24,0)(2083.24,565.48)(2085.77,565.48)(2085.77,0)
\z(2085.77,0)(2085.77,533.04)(2088.31,533.04)(2088.31,0)
\z(2088.31,0)(2088.31,509.25)(2090.85,509.25)(2090.85,0)
\z(2090.85,0)(2090.85,517.90)(2093.38,517.90)(2093.38,0)
\z(2093.38,0)(2093.38,498.44)(2095.92,498.44)(2095.92,0)
\z(2095.92,0)(2095.92,528.72)(2098.45,528.72)(2098.45,0)
\z(2098.45,0)(2098.45,462.76)(2100.99,462.76)(2100.99,0)
\z(2100.99,0)(2100.99,494.12)(2103.52,494.12)(2103.52,0)
\z(2103.52,0)(2103.52,499.52)(2106.06,499.52)(2106.06,0)
\z(2106.06,0)(2106.06,508.17)(2108.59,508.17)(2108.59,0)
\z(2108.59,0)(2108.59,489.79)(2111.13,489.79)(2111.13,0)
\z(2111.13,0)(2111.13,502.77)(2113.66,502.77)(2113.66,0)
\z(2113.66,0)(2113.66,533.04)(2116.20,533.04)(2116.20,0)
\z(2116.20,0)(2116.20,518.98)(2118.74,518.98)(2118.74,0)
\z(2118.74,0)(2118.74,506.01)(2121.27,506.01)(2121.27,0)
\z(2121.27,0)(2121.27,478.98)(2123.81,478.98)(2123.81,0)
\z(2123.81,0)(2123.81,506.01)(2126.34,506.01)(2126.34,0)
\z(2126.34,0)(2126.34,470.33)(2128.88,470.33)(2128.88,0)
\z(2128.88,0)(2128.88,499.52)(2131.41,499.52)(2131.41,0)
\z(2131.41,0)(2131.41,489.79)(2133.95,489.79)(2133.95,0)
\z(2133.95,0)(2133.95,459.52)(2136.48,459.52)(2136.48,0)
\z(2136.48,0)(2136.48,544.93)(2139.02,544.93)(2139.02,0)
\z(2139.02,0)(2139.02,528.72)(2141.56,528.72)(2141.56,0)
\z(2141.56,0)(2141.56,508.17)(2144.09,508.17)(2144.09,0)
\z(2144.09,0)(2144.09,478.98)(2146.63,478.98)(2146.63,0)
\z(2146.63,0)(2146.63,458.44)(2149.16,458.44)(2149.16,0)
\z(2149.16,0)(2149.16,534.12)(2151.70,534.12)(2151.70,0)
\z(2151.70,0)(2151.70,474.65)(2154.23,474.65)(2154.23,0)
\z(2154.23,0)(2154.23,453.03)(2156.77,453.03)(2156.77,0)
\z(2156.77,0)(2156.77,496.28)(2159.30,496.28)(2159.30,0)
\z(2159.30,0)(2159.30,506.01)(2161.84,506.01)(2161.84,0)
\z(2161.84,0)(2161.84,513.58)(2164.38,513.58)(2164.38,0)
\z(2164.38,0)(2164.38,485.47)(2166.91,485.47)(2166.91,0)
\z(2166.91,0)(2166.91,477.90)(2169.45,477.90)(2169.45,0)
\z(2169.45,0)(2169.45,433.57)(2171.98,433.57)(2171.98,0)
\z(2171.98,0)(2171.98,486.55)(2174.52,486.55)(2174.52,0)
\z(2174.52,0)(2174.52,483.30)(2177.05,483.30)(2177.05,0)
\z(2177.05,0)(2177.05,433.57)(2179.59,433.57)(2179.59,0)
\z(2179.59,0)(2179.59,416.27)(2182.12,416.27)(2182.12,0)
\z(2182.12,0)(2182.12,417.35)(2184.66,417.35)(2184.66,0)
\z(2184.66,0)(2184.66,434.65)(2187.20,434.65)(2187.20,0)
\z(2187.20,0)(2187.20,427.08)(2189.73,427.08)(2189.73,0)
\z(2189.73,0)(2189.73,518.98)(2192.27,518.98)(2192.27,0)
\z(2192.27,0)(2192.27,480.06)(2194.80,480.06)(2194.80,0)
\z(2194.80,0)(2194.80,497.36)(2197.34,497.36)(2197.34,0)
\z(2197.34,0)(2197.34,429.24)(2199.87,429.24)(2199.87,0)
\z(2199.87,0)(2199.87,437.89)(2202.41,437.89)(2202.41,0)
\z(2202.41,0)(2202.41,448.71)(2204.94,448.71)(2204.94,0)
\z(2204.94,0)(2204.94,448.71)(2207.48,448.71)(2207.48,0)
\z(2207.48,0)(2207.48,459.52)(2210.02,459.52)(2210.02,0)
\z(2210.02,0)(2210.02,431.41)(2212.55,431.41)(2212.55,0)
\z(2212.55,0)(2212.55,449.79)(2215.09,449.79)(2215.09,0)
\z(2215.09,0)(2215.09,411.94)(2217.62,411.94)(2217.62,0)
\z(2217.62,0)(2217.62,417.35)(2220.16,417.35)(2220.16,0)
\z(2220.16,0)(2220.16,453.03)(2222.69,453.03)(2222.69,0)
\z(2222.69,0)(2222.69,398.97)(2225.23,398.97)(2225.23,0)
\z(2225.23,0)(2225.23,436.81)(2227.76,436.81)(2227.76,0)
\z(2227.76,0)(2227.76,405.46)(2230.30,405.46)(2230.30,0)
\z(2230.30,0)(2230.30,443.30)(2232.84,443.30)(2232.84,0)
\z(2232.84,0)(2232.84,434.65)(2235.37,434.65)(2235.37,0)
\z(2235.37,0)(2235.37,393.56)(2237.91,393.56)(2237.91,0)
\z(2237.91,0)(2237.91,384.91)(2240.44,384.91)(2240.44,0)
\z(2240.44,0)(2240.44,401.13)(2242.98,401.13)(2242.98,0)
\z(2242.98,0)(2242.98,383.83)(2245.51,383.83)(2245.51,0)
\z(2245.51,0)(2245.51,442.22)(2248.05,442.22)(2248.05,0)
\z(2248.05,0)(2248.05,419.51)(2250.58,419.51)(2250.58,0)
\z(2250.58,0)(2250.58,413.03)(2253.12,413.03)(2253.12,0)
\z(2253.12,0)(2253.12,413.03)(2255.65,413.03)(2255.65,0)
\z(2255.65,0)(2255.65,364.37)(2258.19,364.37)(2258.19,0)
\z(2258.19,0)(2258.19,417.35)(2260.73,417.35)(2260.73,0)
\z(2260.73,0)(2260.73,418.43)(2263.26,418.43)(2263.26,0)
\z(2263.26,0)(2263.26,385.99)(2265.80,385.99)(2265.80,0)
\z(2265.80,0)(2265.80,406.54)(2268.33,406.54)(2268.33,0)
\z(2268.33,0)(2268.33,400.05)(2270.87,400.05)(2270.87,0)
\z(2270.87,0)(2270.87,367.61)(2273.40,367.61)(2273.40,0)
\z(2273.40,0)(2273.40,417.35)(2275.94,417.35)(2275.94,0)
\z(2275.94,0)(2275.94,397.89)(2278.47,397.89)(2278.47,0)
\z(2278.47,0)(2278.47,406.54)(2281.01,406.54)(2281.01,0)
\z(2281.01,0)(2281.01,368.70)(2283.55,368.70)(2283.55,0)
\z(2283.55,0)(2283.55,422.76)(2286.08,422.76)(2286.08,0)
\z(2286.08,0)(2286.08,347.07)(2288.62,347.07)(2288.62,0)
\z(2288.62,0)(2288.62,361.13)(2291.15,361.13)(2291.15,0)
\z(2291.15,0)(2291.15,347.07)(2293.69,347.07)(2293.69,0)
\z(2293.69,0)(2293.69,342.75)(2296.22,342.75)(2296.22,0)
\z(2296.22,0)(2296.22,371.94)(2298.76,371.94)(2298.76,0)
\z(2298.76,0)(2298.76,337.34)(2301.29,337.34)(2301.29,0)
\z(2301.29,0)(2301.29,358.96)(2303.83,358.96)(2303.83,0)
\z(2303.83,0)(2303.83,347.07)(2306.37,347.07)(2306.37,0)
\z(2306.37,0)(2306.37,324.37)(2308.90,324.37)(2308.90,0)
\z(2308.90,0)(2308.90,314.63)(2311.44,314.63)(2311.44,0)
\z(2311.44,0)(2311.44,299.50)(2313.97,299.50)(2313.97,0)
\z(2313.97,0)(2313.97,336.26)(2316.51,336.26)(2316.51,0)
\z(2316.51,0)(2316.51,365.45)(2319.04,365.45)(2319.04,0)
\z(2319.04,0)(2319.04,358.96)(2321.58,358.96)(2321.58,0)
\z(2321.58,0)(2321.58,311.39)(2324.11,311.39)(2324.11,0)
\z(2324.11,0)(2324.11,328.69)(2326.65,328.69)(2326.65,0)
\z(2326.65,0)(2326.65,303.82)(2329.19,303.82)(2329.19,0)
\z(2329.19,0)(2329.19,294.09)(2331.72,294.09)(2331.72,0)
\z(2331.72,0)(2331.72,313.55)(2334.26,313.55)(2334.26,0)
\z(2334.26,0)(2334.26,327.61)(2336.79,327.61)(2336.79,0)
\z(2336.79,0)(2336.79,286.52)(2339.33,286.52)(2339.33,0)
\z(2339.33,0)(2339.33,295.17)(2341.86,295.17)(2341.86,0)
\z(2341.86,0)(2341.86,313.55)(2344.40,313.55)(2344.40,0)
\z(2344.40,0)(2344.40,334.10)(2346.93,334.10)(2346.93,0)
\z(2346.93,0)(2346.93,262.74)(2349.47,262.74)(2349.47,0)
\z(2349.47,0)(2349.47,255.17)(2352.01,255.17)(2352.01,0)
\z(2352.01,0)(2352.01,287.60)(2354.54,287.60)(2354.54,0)
\z(2354.54,0)(2354.54,274.63)(2357.08,274.63)(2357.08,0)
\z(2357.08,0)(2357.08,256.25)(2359.61,256.25)(2359.61,0)
\z(2359.61,0)(2359.61,267.06)(2362.15,267.06)(2362.15,0)
\z(2362.15,0)(2362.15,270.30)(2364.68,270.30)(2364.68,0)
\z(2364.68,0)(2364.68,240.03)(2367.22,240.03)(2367.22,0)
\z(2367.22,0)(2367.22,245.44)(2369.75,245.44)(2369.75,0)
\z(2369.75,0)(2369.75,241.11)(2372.29,241.11)(2372.29,0)
\z(2372.29,0)(2372.29,254.09)(2374.82,254.09)(2374.82,0)
\z(2374.82,0)(2374.82,209.76)(2377.36,209.76)(2377.36,0)
\z(2377.36,0)(2377.36,228.14)(2379.90,228.14)(2379.90,0)
\z(2379.90,0)(2379.90,218.41)(2382.43,218.41)(2382.43,0)
\z(2382.43,0)(2382.43,229.22)(2384.97,229.22)(2384.97,0)
\z(2384.97,0)(2384.97,217.32)(2387.50,217.32)(2387.50,0)
\z(2387.50,0)(2387.50,223.81)(2390.04,223.81)(2390.04,0)
\z(2390.04,0)(2390.04,222.73)(2392.57,222.73)(2392.57,0)
\z(2392.57,0)(2392.57,174.08)(2395.11,174.08)(2395.11,0)
\z(2395.11,0)(2395.11,174.08)(2397.64,174.08)(2397.64,0)
\z(2397.64,0)(2397.64,187.05)(2400.18,187.05)(2400.18,0)
\z(2400.18,0)(2400.18,175.16)(2402.72,175.16)(2402.72,0)
\z(2402.72,0)(2402.72,167.59)(2405.25,167.59)(2405.25,0)
\z(2405.25,0)(2405.25,139.48)(2407.79,139.48)(2407.79,0)
\z(2407.79,0)(2407.79,117.85)(2410.32,117.85)(2410.32,0)
\z(2410.32,0)(2410.32,99.47)(2412.86,99.47)(2412.86,0)
\z(2412.86,0)(2412.86,99.47)(2415.39,99.47)(2415.39,0)
\z(2415.39,0)(2415.39,108.12)(2417.93,108.12)(2417.93,0)
\z(2417.93,0)(2417.93,64.87)(2420.46,64.87)(2420.46,0)
\z(2420.46,0)(2420.46,28.11)(2423.00,28.11)(2423.00,0)
\color{black}
\dottedline{12}(1400,426.00)(2400,426.00)
\z(1400,0)(1400,640)(2400,640)(2400,0)(1400,0)
\end{picture}

%% file: g2_a1_12_14.tex
\hspace{-30pt}
\begin{picture}(2400,640)
\color[rgb]{0.0,0.2,0.8}
\z(110.87,0)(110.87,6.17)(154.30,6.17)(154.30,0)
\z(154.30,0)(154.30,12.35)(197.74,12.35)(197.74,0)
\z(197.74,0)(197.74,30.87)(241.17,30.87)(241.17,0)
\z(241.17,0)(241.17,104.96)(284.61,104.96)(284.61,0)
\z(284.61,0)(284.61,98.78)(328.04,98.78)(328.04,0)
\z(328.04,0)(328.04,203.74)(371.48,203.74)(371.48,0)
\z(371.48,0)(371.48,271.65)(414.91,271.65)(414.91,0)
\z(414.91,0)(414.91,388.96)(458.35,388.96)(458.35,0)
\z(458.35,0)(458.35,407.48)(501.78,407.48)(501.78,0)
\z(501.78,0)(501.78,524.78)(545.22,524.78)(545.22,0)
\z(545.22,0)(545.22,395.13)(588.65,395.13)(588.65,0)
\z(588.65,0)(588.65,333.39)(632.09,333.39)(632.09,0)
\z(632.09,0)(632.09,284.00)(675.52,284.00)(675.52,0)
\z(675.52,0)(675.52,172.87)(718.96,172.87)(718.96,0)
\z(718.96,0)(718.96,111.13)(762.39,111.13)(762.39,0)
\z(762.39,0)(762.39,67.91)(805.83,67.91)(805.83,0)
\z(805.83,0)(805.83,37.04)(849.26,37.04)(849.26,0)
\z(849.26,0)(849.26,18.52)(892.70,18.52)(892.70,0)
\color{black}
\dottedline{12}(0,142.00)(1000,142.00)
\z(0,0)(0,640)(1000,640)(1000,0)(0,0)
\color[rgb]{0.0,0.2,0.8}
\z(1470.47,0)(1470.47,3.30)(1493.70,3.30)(1493.70,0)
\z(1493.70,0)(1493.70,9.91)(1516.93,9.91)(1516.93,0)
\z(1516.93,0)(1516.93,3.30)(1540.16,3.30)(1540.16,0)
\z(1540.16,0)(1540.16,3.30)(1563.40,3.30)(1563.40,0)
\z(1563.40,0)(1563.40,6.60)(1586.63,6.60)(1586.63,0)
\z(1586.63,0)(1586.63,36.33)(1609.86,36.33)(1609.86,0)
\z(1609.86,0)(1609.86,13.21)(1633.09,13.21)(1633.09,0)
\z(1633.09,0)(1633.09,46.23)(1656.33,46.23)(1656.33,0)
\z(1656.33,0)(1656.33,79.26)(1679.56,79.26)(1679.56,0)
\z(1679.56,0)(1679.56,108.98)(1702.79,108.98)(1702.79,0)
\z(1702.79,0)(1702.79,122.19)(1726.02,122.19)(1726.02,0)
\z(1726.02,0)(1726.02,194.84)(1749.26,194.84)(1749.26,0)
\z(1749.26,0)(1749.26,211.35)(1772.49,211.35)(1772.49,0)
\z(1772.49,0)(1772.49,260.88)(1795.72,260.88)(1795.72,0)
\z(1795.72,0)(1795.72,257.58)(1818.95,257.58)(1818.95,0)
\z(1818.95,0)(1818.95,300.51)(1842.19,300.51)(1842.19,0)
\z(1842.19,0)(1842.19,363.26)(1865.42,363.26)(1865.42,0)
\z(1865.42,0)(1865.42,485.44)(1888.65,485.44)(1888.65,0)
\z(1888.65,0)(1888.65,402.88)(1911.88,402.88)(1911.88,0)
\z(1911.88,0)(1911.88,544.88)(1935.12,544.88)(1935.12,0)
\z(1935.12,0)(1935.12,511.86)(1958.35,511.86)(1958.35,0)
\z(1958.35,0)(1958.35,396.28)(1981.58,396.28)(1981.58,0)
\z(1981.58,0)(1981.58,369.86)(2004.81,369.86)(2004.81,0)
\z(2004.81,0)(2004.81,336.84)(2028.05,336.84)(2028.05,0)
\z(2028.05,0)(2028.05,284.00)(2051.28,284.00)(2051.28,0)
\z(2051.28,0)(2051.28,250.98)(2074.51,250.98)(2074.51,0)
\z(2074.51,0)(2074.51,224.56)(2097.74,224.56)(2097.74,0)
\z(2097.74,0)(2097.74,122.19)(2120.98,122.19)(2120.98,0)
\z(2120.98,0)(2120.98,99.07)(2144.21,99.07)(2144.21,0)
\z(2144.21,0)(2144.21,62.74)(2167.44,62.74)(2167.44,0)
\z(2167.44,0)(2167.44,39.63)(2190.67,39.63)(2190.67,0)
\z(2190.67,0)(2190.67,52.84)(2213.91,52.84)(2213.91,0)
\z(2213.91,0)(2213.91,29.72)(2237.14,29.72)(2237.14,0)
\z(2237.14,0)(2237.14,13.21)(2260.37,13.21)(2260.37,0)
\z(2260.37,0)(2260.37,9.91)(2283.60,9.91)(2283.60,0)
\z(2283.60,0)(2283.60,6.60)(2306.84,6.60)(2306.84,0)
\color{black}
\dottedline{12}(1400,142.00)(2400,142.00)
\z(1400,0)(1400,640)(2400,640)(2400,0)(1400,0)
\end{picture}

%% file: g2_a1_16_18.tex
\hspace{-30pt}
\begin{picture}(2400,640)
\color[rgb]{0.0,0.2,0.8}
\z(86.44,0)(86.44,1.77)(98.92,1.77)(98.92,0)
\z(98.92,0)(98.92,3.55)(111.41,3.55)(111.41,0)
\z(111.41,0)(111.41,7.10)(123.90,7.10)(123.90,0)
\z(123.90,0)(123.90,1.77)(136.39,1.77)(136.39,0)
\z(136.39,0)(136.39,7.10)(148.88,7.10)(148.88,0)
\z(161.36,0)(161.36,3.55)(173.85,3.55)(173.85,0)
\z(173.85,0)(173.85,14.20)(186.34,14.20)(186.34,0)
\z(186.34,0)(186.34,21.30)(198.83,21.30)(198.83,0)
\z(198.83,0)(198.83,21.30)(211.31,21.30)(211.31,0)
\z(211.31,0)(211.31,21.30)(223.80,21.30)(223.80,0)
\z(223.80,0)(223.80,24.85)(236.29,24.85)(236.29,0)
\z(236.29,0)(236.29,46.15)(248.78,46.15)(248.78,0)
\z(248.78,0)(248.78,49.70)(261.26,49.70)(261.26,0)
\z(261.26,0)(261.26,72.77)(273.75,72.77)(273.75,0)
\z(273.75,0)(273.75,76.33)(286.24,76.33)(286.24,0)
\z(286.24,0)(286.24,85.20)(298.73,85.20)(298.73,0)
\z(298.73,0)(298.73,143.78)(311.21,143.78)(311.21,0)
\z(311.21,0)(311.21,140.22)(323.70,140.22)(323.70,0)
\z(323.70,0)(323.70,161.53)(336.19,161.53)(336.19,0)
\z(336.19,0)(336.19,173.95)(348.68,173.95)(348.68,0)
\z(348.68,0)(348.68,207.67)(361.16,207.67)(361.16,0)
\z(361.16,0)(361.16,198.80)(373.65,198.80)(373.65,0)
\z(373.65,0)(373.65,239.62)(386.14,239.62)(386.14,0)
\z(386.14,0)(386.14,228.97)(398.63,228.97)(398.63,0)
\z(398.63,0)(398.63,241.40)(411.11,241.40)(411.11,0)
\z(411.11,0)(411.11,321.27)(423.60,321.27)(423.60,0)
\z(423.60,0)(423.60,315.95)(436.09,315.95)(436.09,0)
\z(436.09,0)(436.09,351.45)(448.58,351.45)(448.58,0)
\z(448.58,0)(448.58,390.50)(461.06,390.50)(461.06,0)
\z(461.06,0)(461.06,443.75)(473.55,443.75)(473.55,0)
\z(473.55,0)(473.55,431.32)(486.04,431.32)(486.04,0)
\z(486.04,0)(486.04,454.40)(498.53,454.40)(498.53,0)
\z(498.53,0)(498.53,472.15)(511.01,472.15)(511.01,0)
\z(511.01,0)(511.01,484.57)(523.50,484.57)(523.50,0)
\z(523.50,0)(523.50,539.60)(535.99,539.60)(535.99,0)
\z(535.99,0)(535.99,498.77)(548.48,498.77)(548.48,0)
\z(548.48,0)(548.48,488.12)(560.96,488.12)(560.96,0)
\z(560.96,0)(560.96,404.70)(573.45,404.70)(573.45,0)
\z(573.45,0)(573.45,383.40)(585.94,383.40)(585.94,0)
\z(585.94,0)(585.94,381.62)(598.43,381.62)(598.43,0)
\z(598.43,0)(598.43,356.77)(610.91,356.77)(610.91,0)
\z(610.91,0)(610.91,305.30)(623.40,305.30)(623.40,0)
\z(623.40,0)(623.40,319.50)(635.89,319.50)(635.89,0)
\z(635.89,0)(635.89,330.15)(648.37,330.15)(648.37,0)
\z(648.37,0)(648.37,250.27)(660.86,250.27)(660.86,0)
\z(660.86,0)(660.86,202.35)(673.35,202.35)(673.35,0)
\z(673.35,0)(673.35,211.22)(685.84,211.22)(685.84,0)
\z(685.84,0)(685.84,143.78)(698.32,143.78)(698.32,0)
\z(698.32,0)(698.32,168.62)(710.81,168.62)(710.81,0)
\z(710.81,0)(710.81,124.25)(723.30,124.25)(723.30,0)
\z(723.30,0)(723.30,94.07)(735.79,94.07)(735.79,0)
\z(735.79,0)(735.79,94.07)(748.27,94.07)(748.27,0)
\z(748.27,0)(748.27,94.07)(760.76,94.07)(760.76,0)
\z(760.76,0)(760.76,71.00)(773.25,71.00)(773.25,0)
\z(773.25,0)(773.25,51.47)(785.74,51.47)(785.74,0)
\z(785.74,0)(785.74,40.82)(798.22,40.82)(798.22,0)
\z(798.22,0)(798.22,51.47)(810.71,51.47)(810.71,0)
\z(810.71,0)(810.71,40.82)(823.20,40.82)(823.20,0)
\z(823.20,0)(823.20,26.62)(835.69,26.62)(835.69,0)
\z(835.69,0)(835.69,17.75)(848.17,17.75)(848.17,0)
\z(848.17,0)(848.17,10.65)(860.66,10.65)(860.66,0)
\z(860.66,0)(860.66,14.20)(873.15,14.20)(873.15,0)
\z(873.15,0)(873.15,8.88)(885.64,8.88)(885.64,0)
\z(885.64,0)(885.64,10.65)(898.12,10.65)(898.12,0)
\z(898.12,0)(898.12,3.55)(910.61,3.55)(910.61,0)
\z(910.61,0)(910.61,5.32)(923.10,5.32)(923.10,0)
\z(923.10,0)(923.10,1.77)(935.59,1.77)(935.59,0)
\z(948.07,0)(948.07,1.77)(960.56,1.77)(960.56,0)
\color{black}
\dottedline{12}(0,142.00)(1000,142.00)
\z(0,0)(0,640)(1000,640)(1000,0)(0,0)
\color[rgb]{0.0,0.2,0.8}
\z(1457.08,0)(1457.08,0.94)(1463.70,0.94)(1463.70,0)
\z(1483.54,0)(1483.54,0.94)(1490.16,0.94)(1490.16,0)
\z(1496.77,0)(1496.77,2.82)(1503.39,2.82)(1503.39,0)
\z(1503.39,0)(1503.39,1.88)(1510.01,1.88)(1510.01,0)
\z(1510.01,0)(1510.01,3.76)(1516.62,3.76)(1516.62,0)
\z(1516.62,0)(1516.62,2.82)(1523.24,2.82)(1523.24,0)
\z(1523.24,0)(1523.24,1.88)(1529.85,1.88)(1529.85,0)
\z(1529.85,0)(1529.85,2.82)(1536.47,2.82)(1536.47,0)
\z(1536.47,0)(1536.47,2.82)(1543.09,2.82)(1543.09,0)
\z(1543.09,0)(1543.09,5.64)(1549.70,5.64)(1549.70,0)
\z(1549.70,0)(1549.70,1.88)(1556.32,1.88)(1556.32,0)
\z(1556.32,0)(1556.32,7.52)(1562.93,7.52)(1562.93,0)
\z(1562.93,0)(1562.93,6.58)(1569.55,6.58)(1569.55,0)
\z(1569.55,0)(1569.55,10.34)(1576.17,10.34)(1576.17,0)
\z(1576.17,0)(1576.17,13.17)(1582.78,13.17)(1582.78,0)
\z(1582.78,0)(1582.78,16.93)(1589.40,16.93)(1589.40,0)
\z(1589.40,0)(1589.40,19.75)(1596.01,19.75)(1596.01,0)
\z(1596.01,0)(1596.01,20.69)(1602.63,20.69)(1602.63,0)
\z(1602.63,0)(1602.63,21.63)(1609.25,21.63)(1609.25,0)
\z(1609.25,0)(1609.25,15.05)(1615.86,15.05)(1615.86,0)
\z(1615.86,0)(1615.86,25.39)(1622.48,25.39)(1622.48,0)
\z(1622.48,0)(1622.48,38.56)(1629.09,38.56)(1629.09,0)
\z(1629.09,0)(1629.09,25.39)(1635.71,25.39)(1635.71,0)
\z(1635.71,0)(1635.71,41.38)(1642.32,41.38)(1642.32,0)
\z(1642.32,0)(1642.32,34.79)(1648.94,34.79)(1648.94,0)
\z(1648.94,0)(1648.94,47.96)(1655.56,47.96)(1655.56,0)
\z(1655.56,0)(1655.56,50.78)(1662.17,50.78)(1662.17,0)
\z(1662.17,0)(1662.17,51.72)(1668.79,51.72)(1668.79,0)
\z(1668.79,0)(1668.79,72.41)(1675.40,72.41)(1675.40,0)
\z(1675.40,0)(1675.40,70.53)(1682.02,70.53)(1682.02,0)
\z(1682.02,0)(1682.02,94.04)(1688.64,94.04)(1688.64,0)
\z(1688.64,0)(1688.64,76.17)(1695.25,76.17)(1695.25,0)
\z(1695.25,0)(1695.25,101.56)(1701.87,101.56)(1701.87,0)
\z(1701.87,0)(1701.87,111.91)(1708.48,111.91)(1708.48,0)
\z(1708.48,0)(1708.48,104.38)(1715.10,104.38)(1715.10,0)
\z(1715.10,0)(1715.10,126.95)(1721.72,126.95)(1721.72,0)
\z(1721.72,0)(1721.72,133.54)(1728.33,133.54)(1728.33,0)
\z(1728.33,0)(1728.33,154.23)(1734.95,154.23)(1734.95,0)
\z(1734.95,0)(1734.95,133.54)(1741.56,133.54)(1741.56,0)
\z(1741.56,0)(1741.56,174.91)(1748.18,174.91)(1748.18,0)
\z(1748.18,0)(1748.18,155.17)(1754.79,155.17)(1754.79,0)
\z(1754.79,0)(1754.79,213.47)(1761.41,213.47)(1761.41,0)
\z(1761.41,0)(1761.41,185.26)(1768.03,185.26)(1768.03,0)
\z(1768.03,0)(1768.03,207.83)(1774.64,207.83)(1774.64,0)
\z(1774.64,0)(1774.64,211.59)(1781.26,211.59)(1781.26,0)
\z(1781.26,0)(1781.26,224.75)(1787.87,224.75)(1787.87,0)
\z(1787.87,0)(1787.87,241.68)(1794.49,241.68)(1794.49,0)
\z(1794.49,0)(1794.49,254.85)(1801.11,254.85)(1801.11,0)
\z(1801.11,0)(1801.11,269.89)(1807.72,269.89)(1807.72,0)
\z(1807.72,0)(1807.72,267.07)(1814.34,267.07)(1814.34,0)
\z(1814.34,0)(1814.34,319.74)(1820.95,319.74)(1820.95,0)
\z(1820.95,0)(1820.95,315.97)(1827.57,315.97)(1827.57,0)
\z(1827.57,0)(1827.57,333.84)(1834.19,333.84)(1834.19,0)
\z(1834.19,0)(1834.19,348.89)(1840.80,348.89)(1840.80,0)
\z(1840.80,0)(1840.80,339.48)(1847.42,339.48)(1847.42,0)
\z(1847.42,0)(1847.42,394.03)(1854.03,394.03)(1854.03,0)
\z(1854.03,0)(1854.03,382.74)(1860.65,382.74)(1860.65,0)
\z(1860.65,0)(1860.65,400.61)(1867.26,400.61)(1867.26,0)
\z(1867.26,0)(1867.26,406.25)(1873.88,406.25)(1873.88,0)
\z(1873.88,0)(1873.88,426.00)(1880.50,426.00)(1880.50,0)
\z(1880.50,0)(1880.50,491.83)(1887.11,491.83)(1887.11,0)
\z(1887.11,0)(1887.11,486.19)(1893.73,486.19)(1893.73,0)
\z(1893.73,0)(1893.73,464.56)(1900.34,464.56)(1900.34,0)
\z(1900.34,0)(1900.34,507.81)(1906.96,507.81)(1906.96,0)
\z(1906.96,0)(1906.96,468.32)(1913.58,468.32)(1913.58,0)
\z(1913.58,0)(1913.58,506.87)(1920.19,506.87)(1920.19,0)
\z(1920.19,0)(1920.19,508.75)(1926.81,508.75)(1926.81,0)
\z(1926.81,0)(1926.81,545.43)(1933.42,545.43)(1933.42,0)
\z(1933.42,0)(1933.42,450.45)(1940.04,450.45)(1940.04,0)
\z(1940.04,0)(1940.04,479.60)(1946.66,479.60)(1946.66,0)
\z(1946.66,0)(1946.66,485.25)(1953.27,485.25)(1953.27,0)
\z(1953.27,0)(1953.27,484.30)(1959.89,484.30)(1959.89,0)
\z(1959.89,0)(1959.89,418.48)(1966.50,418.48)(1966.50,0)
\z(1966.50,0)(1966.50,436.34)(1973.12,436.34)(1973.12,0)
\z(1973.12,0)(1973.12,376.16)(1979.74,376.16)(1979.74,0)
\z(1979.74,0)(1979.74,421.30)(1986.35,421.30)(1986.35,0)
\z(1986.35,0)(1986.35,386.50)(1992.97,386.50)(1992.97,0)
\z(1992.97,0)(1992.97,379.92)(1999.58,379.92)(1999.58,0)
\z(1999.58,0)(1999.58,363.93)(2006.20,363.93)(2006.20,0)
\z(2006.20,0)(2006.20,353.59)(2012.81,353.59)(2012.81,0)
\z(2012.81,0)(2012.81,303.75)(2019.43,303.75)(2019.43,0)
\z(2019.43,0)(2019.43,268.01)(2026.05,268.01)(2026.05,0)
\z(2026.05,0)(2026.05,268.95)(2032.66,268.95)(2032.66,0)
\z(2032.66,0)(2032.66,300.93)(2039.28,300.93)(2039.28,0)
\z(2039.28,0)(2039.28,276.48)(2045.89,276.48)(2045.89,0)
\z(2045.89,0)(2045.89,262.37)(2052.51,262.37)(2052.51,0)
\z(2052.51,0)(2052.51,236.98)(2059.13,236.98)(2059.13,0)
\z(2059.13,0)(2059.13,214.41)(2065.74,214.41)(2065.74,0)
\z(2065.74,0)(2065.74,242.62)(2072.36,242.62)(2072.36,0)
\z(2072.36,0)(2072.36,203.13)(2078.97,203.13)(2078.97,0)
\z(2078.97,0)(2078.97,191.84)(2085.59,191.84)(2085.59,0)
\z(2085.59,0)(2085.59,199.36)(2092.21,199.36)(2092.21,0)
\z(2092.21,0)(2092.21,169.27)(2098.82,169.27)(2098.82,0)
\z(2098.82,0)(2098.82,142.00)(2105.44,142.00)(2105.44,0)
\z(2105.44,0)(2105.44,151.40)(2112.05,151.40)(2112.05,0)
\z(2112.05,0)(2112.05,128.83)(2118.67,128.83)(2118.67,0)
\z(2118.67,0)(2118.67,147.64)(2125.28,147.64)(2125.28,0)
\z(2125.28,0)(2125.28,113.79)(2131.90,113.79)(2131.90,0)
\z(2131.90,0)(2131.90,100.62)(2138.52,100.62)(2138.52,0)
\z(2138.52,0)(2138.52,79.93)(2145.13,79.93)(2145.13,0)
\z(2145.13,0)(2145.13,94.98)(2151.75,94.98)(2151.75,0)
\z(2151.75,0)(2151.75,77.11)(2158.36,77.11)(2158.36,0)
\z(2158.36,0)(2158.36,86.52)(2164.98,86.52)(2164.98,0)
\z(2164.98,0)(2164.98,81.81)(2171.60,81.81)(2171.60,0)
\z(2171.60,0)(2171.60,62.07)(2178.21,62.07)(2178.21,0)
\z(2178.21,0)(2178.21,51.72)(2184.83,51.72)(2184.83,0)
\z(2184.83,0)(2184.83,47.96)(2191.44,47.96)(2191.44,0)
\z(2191.44,0)(2191.44,40.44)(2198.06,40.44)(2198.06,0)
\z(2198.06,0)(2198.06,30.09)(2204.68,30.09)(2204.68,0)
\z(2204.68,0)(2204.68,45.14)(2211.29,45.14)(2211.29,0)
\z(2211.29,0)(2211.29,33.85)(2217.91,33.85)(2217.91,0)
\z(2217.91,0)(2217.91,38.56)(2224.52,38.56)(2224.52,0)
\z(2224.52,0)(2224.52,31.97)(2231.14,31.97)(2231.14,0)
\z(2231.14,0)(2231.14,22.57)(2237.75,22.57)(2237.75,0)
\z(2237.75,0)(2237.75,19.75)(2244.37,19.75)(2244.37,0)
\z(2244.37,0)(2244.37,16.93)(2250.99,16.93)(2250.99,0)
\z(2250.99,0)(2250.99,21.63)(2257.60,21.63)(2257.60,0)
\z(2257.60,0)(2257.60,12.23)(2264.22,12.23)(2264.22,0)
\z(2264.22,0)(2264.22,9.40)(2270.83,9.40)(2270.83,0)
\z(2270.83,0)(2270.83,12.23)(2277.45,12.23)(2277.45,0)
\z(2277.45,0)(2277.45,14.11)(2284.07,14.11)(2284.07,0)
\z(2284.07,0)(2284.07,4.70)(2290.68,4.70)(2290.68,0)
\z(2290.68,0)(2290.68,6.58)(2297.30,6.58)(2297.30,0)
\z(2297.30,0)(2297.30,1.88)(2303.91,1.88)(2303.91,0)
\z(2303.91,0)(2303.91,4.70)(2310.53,4.70)(2310.53,0)
\z(2310.53,0)(2310.53,2.82)(2317.15,2.82)(2317.15,0)
\z(2317.15,0)(2317.15,2.82)(2323.76,2.82)(2323.76,0)
\z(2323.76,0)(2323.76,0.94)(2330.38,0.94)(2330.38,0)
\z(2330.38,0)(2330.38,2.82)(2336.99,2.82)(2336.99,0)
\z(2343.61,0)(2343.61,0.94)(2350.23,0.94)(2350.23,0)
\color{black}
\dottedline{12}(1400,142.00)(2400,142.00)
\z(1400,0)(1400,640)(2400,640)(2400,0)(1400,0)
\end{picture}

%% file: g2_a1_20_22.tex
\hspace{-30pt}
\begin{picture}(2400,640)
\color[rgb]{0.0,0.2,0.8}
\z(62.42,0)(62.42,0.50)(65.92,0.50)(65.92,0)
\z(69.41,0)(69.41,0.50)(72.90,0.50)(72.90,0)
\z(79.89,0)(79.89,0.50)(83.38,0.50)(83.38,0)
\z(86.87,0)(86.87,0.50)(90.37,0.50)(90.37,0)
\z(93.86,0)(93.86,0.50)(97.35,0.50)(97.35,0)
\z(97.35,0)(97.35,1.99)(100.85,1.99)(100.85,0)
\z(100.85,0)(100.85,0.99)(104.34,0.99)(104.34,0)
\z(104.34,0)(104.34,1.99)(107.83,1.99)(107.83,0)
\z(107.83,0)(107.83,1.49)(111.33,1.49)(111.33,0)
\z(111.33,0)(111.33,2.48)(114.82,2.48)(114.82,0)
\z(114.82,0)(114.82,2.98)(118.31,2.98)(118.31,0)
\z(118.31,0)(118.31,1.99)(121.80,1.99)(121.80,0)
\z(121.80,0)(121.80,2.48)(125.30,2.48)(125.30,0)
\z(125.30,0)(125.30,2.48)(128.79,2.48)(128.79,0)
\z(128.79,0)(128.79,1.99)(132.28,1.99)(132.28,0)
\z(132.28,0)(132.28,3.97)(135.78,3.97)(135.78,0)
\z(135.78,0)(135.78,2.48)(139.27,2.48)(139.27,0)
\z(139.27,0)(139.27,4.47)(142.76,4.47)(142.76,0)
\z(142.76,0)(142.76,4.97)(146.26,4.97)(146.26,0)
\z(146.26,0)(146.26,5.96)(149.75,5.96)(149.75,0)
\z(149.75,0)(149.75,8.44)(153.24,8.44)(153.24,0)
\z(153.24,0)(153.24,3.48)(156.73,3.48)(156.73,0)
\z(156.73,0)(156.73,4.47)(160.23,4.47)(160.23,0)
\z(160.23,0)(160.23,8.44)(163.72,8.44)(163.72,0)
\z(163.72,0)(163.72,6.95)(167.21,6.95)(167.21,0)
\z(167.21,0)(167.21,8.94)(170.71,8.94)(170.71,0)
\z(170.71,0)(170.71,8.44)(174.20,8.44)(174.20,0)
\z(174.20,0)(174.20,7.94)(177.69,7.94)(177.69,0)
\z(177.69,0)(177.69,15.89)(181.19,15.89)(181.19,0)
\z(181.19,0)(181.19,13.41)(184.68,13.41)(184.68,0)
\z(184.68,0)(184.68,14.40)(188.17,14.40)(188.17,0)
\z(188.17,0)(188.17,16.38)(191.66,16.38)(191.66,0)
\z(191.66,0)(191.66,13.41)(195.16,13.41)(195.16,0)
\z(195.16,0)(195.16,20.36)(198.65,20.36)(198.65,0)
\z(198.65,0)(198.65,17.87)(202.14,17.87)(202.14,0)
\z(202.14,0)(202.14,16.88)(205.64,16.88)(205.64,0)
\z(205.64,0)(205.64,24.83)(209.13,24.83)(209.13,0)
\z(209.13,0)(209.13,18.87)(212.62,18.87)(212.62,0)
\z(212.62,0)(212.62,24.83)(216.12,24.83)(216.12,0)
\z(216.12,0)(216.12,23.34)(219.61,23.34)(219.61,0)
\z(219.61,0)(219.61,24.83)(223.10,24.83)(223.10,0)
\z(223.10,0)(223.10,34.76)(226.59,34.76)(226.59,0)
\z(226.59,0)(226.59,28.30)(230.09,28.30)(230.09,0)
\z(230.09,0)(230.09,33.76)(233.58,33.76)(233.58,0)
\z(233.58,0)(233.58,31.78)(237.07,31.78)(237.07,0)
\z(237.07,0)(237.07,35.75)(240.57,35.75)(240.57,0)
\z(240.57,0)(240.57,43.20)(244.06,43.20)(244.06,0)
\z(244.06,0)(244.06,34.26)(247.55,34.26)(247.55,0)
\z(247.55,0)(247.55,41.21)(251.05,41.21)(251.05,0)
\z(251.05,0)(251.05,55.11)(254.54,55.11)(254.54,0)
\z(254.54,0)(254.54,44.69)(258.03,44.69)(258.03,0)
\z(258.03,0)(258.03,43.20)(261.52,43.20)(261.52,0)
\z(261.52,0)(261.52,58.59)(265.02,58.59)(265.02,0)
\z(265.02,0)(265.02,54.12)(268.51,54.12)(268.51,0)
\z(268.51,0)(268.51,69.01)(272.00,69.01)(272.00,0)
\z(272.00,0)(272.00,57.10)(275.50,57.10)(275.50,0)
\z(275.50,0)(275.50,68.52)(278.99,68.52)(278.99,0)
\z(278.99,0)(278.99,65.54)(282.48,65.54)(282.48,0)
\z(282.48,0)(282.48,87.38)(285.98,87.38)(285.98,0)
\z(285.98,0)(285.98,75.47)(289.47,75.47)(289.47,0)
\z(289.47,0)(289.47,79.94)(292.96,79.94)(292.96,0)
\z(292.96,0)(292.96,100.79)(296.45,100.79)(296.45,0)
\z(296.45,0)(296.45,83.91)(299.95,83.91)(299.95,0)
\z(299.95,0)(299.95,103.27)(303.44,103.27)(303.44,0)
\z(303.44,0)(303.44,109.23)(306.93,109.23)(306.93,0)
\z(306.93,0)(306.93,120.65)(310.43,120.65)(310.43,0)
\z(310.43,0)(310.43,105.26)(313.92,105.26)(313.92,0)
\z(313.92,0)(313.92,114.69)(317.41,114.69)(317.41,0)
\z(317.41,0)(317.41,125.62)(320.91,125.62)(320.91,0)
\z(320.91,0)(320.91,114.69)(324.40,114.69)(324.40,0)
\z(324.40,0)(324.40,121.64)(327.89,121.64)(327.89,0)
\z(327.89,0)(327.89,131.08)(331.38,131.08)(331.38,0)
\z(331.38,0)(331.38,138.52)(334.88,138.52)(334.88,0)
\z(334.88,0)(334.88,140.01)(338.37,140.01)(338.37,0)
\z(338.37,0)(338.37,155.41)(341.86,155.41)(341.86,0)
\z(341.86,0)(341.86,158.38)(345.36,158.38)(345.36,0)
\z(345.36,0)(345.36,166.33)(348.85,166.33)(348.85,0)
\z(348.85,0)(348.85,157.89)(352.34,157.89)(352.34,0)
\z(352.34,0)(352.34,162.36)(355.84,162.36)(355.84,0)
\z(355.84,0)(355.84,186.19)(359.33,186.19)(359.33,0)
\z(359.33,0)(359.33,196.12)(362.82,196.12)(362.82,0)
\z(362.82,0)(362.82,195.62)(366.31,195.62)(366.31,0)
\z(366.31,0)(366.31,209.52)(369.81,209.52)(369.81,0)
\z(369.81,0)(369.81,202.08)(373.30,202.08)(373.30,0)
\z(373.30,0)(373.30,201.58)(376.79,201.58)(376.79,0)
\z(376.79,0)(376.79,198.60)(380.29,198.60)(380.29,0)
\z(380.29,0)(380.29,233.36)(383.78,233.36)(383.78,0)
\z(383.78,0)(383.78,228.39)(387.27,228.39)(387.27,0)
\z(387.27,0)(387.27,224.42)(390.77,224.42)(390.77,0)
\z(390.77,0)(390.77,251.23)(394.26,251.23)(394.26,0)
\z(394.26,0)(394.26,249.74)(397.75,249.74)(397.75,0)
\z(397.75,0)(397.75,268.61)(401.24,268.61)(401.24,0)
\z(401.24,0)(401.24,267.12)(404.74,267.12)(404.74,0)
\z(404.74,0)(404.74,258.18)(408.23,258.18)(408.23,0)
\z(408.23,0)(408.23,284.99)(411.72,284.99)(411.72,0)
\z(411.72,0)(411.72,287.97)(415.22,287.97)(415.22,0)
\z(415.22,0)(415.22,301.38)(418.71,301.38)(418.71,0)
\z(418.71,0)(418.71,295.92)(422.20,295.92)(422.20,0)
\z(422.20,0)(422.20,305.35)(425.70,305.35)(425.70,0)
\z(425.70,0)(425.70,304.85)(429.19,304.85)(429.19,0)
\z(429.19,0)(429.19,332.16)(432.68,332.16)(432.68,0)
\z(432.68,0)(432.68,327.69)(436.17,327.69)(436.17,0)
\z(436.17,0)(436.17,351.52)(439.67,351.52)(439.67,0)
\z(439.67,0)(439.67,346.56)(443.16,346.56)(443.16,0)
\z(443.16,0)(443.16,356.99)(446.65,356.99)(446.65,0)
\z(446.65,0)(446.65,361.45)(450.15,361.45)(450.15,0)
\z(450.15,0)(450.15,393.23)(453.64,393.23)(453.64,0)
\z(453.64,0)(453.64,391.74)(457.13,391.74)(457.13,0)
\z(457.13,0)(457.13,369.90)(460.63,369.90)(460.63,0)
\z(460.63,0)(460.63,420.54)(464.12,420.54)(464.12,0)
\z(464.12,0)(464.12,430.47)(467.61,430.47)(467.61,0)
\z(467.61,0)(467.61,421.53)(471.10,421.53)(471.10,0)
\z(471.10,0)(471.10,435.43)(474.60,435.43)(474.60,0)
\z(474.60,0)(474.60,422.03)(478.09,422.03)(478.09,0)
\z(478.09,0)(478.09,426.50)(481.58,426.50)(481.58,0)
\z(481.58,0)(481.58,455.29)(485.08,455.29)(485.08,0)
\z(485.08,0)(485.08,474.16)(488.57,474.16)(488.57,0)
\z(488.57,0)(488.57,459.76)(492.06,459.76)(492.06,0)
\z(492.06,0)(492.06,470.69)(495.56,470.69)(495.56,0)
\z(495.56,0)(495.56,465.72)(499.05,465.72)(499.05,0)
\z(499.05,0)(499.05,471.18)(502.54,471.18)(502.54,0)
\z(502.54,0)(502.54,490.05)(506.03,490.05)(506.03,0)
\z(506.03,0)(506.03,503.95)(509.53,503.95)(509.53,0)
\z(509.53,0)(509.53,467.71)(513.02,467.71)(513.02,0)
\z(513.02,0)(513.02,497.50)(516.51,497.50)(516.51,0)
\z(516.51,0)(516.51,521.83)(520.01,521.83)(520.01,0)
\z(520.01,0)(520.01,480.12)(523.50,480.12)(523.50,0)
\z(523.50,0)(523.50,524.31)(526.99,524.31)(526.99,0)
\z(526.99,0)(526.99,526.79)(530.49,526.79)(530.49,0)
\z(530.49,0)(530.49,476.64)(533.98,476.64)(533.98,0)
\z(533.98,0)(533.98,473.66)(537.47,473.66)(537.47,0)
\z(537.47,0)(537.47,488.06)(540.97,488.06)(540.97,0)
\z(540.97,0)(540.97,482.60)(544.46,482.60)(544.46,0)
\z(544.46,0)(544.46,472.17)(547.95,472.17)(547.95,0)
\z(547.95,0)(547.95,452.81)(551.44,452.81)(551.44,0)
\z(551.44,0)(551.44,459.27)(554.94,459.27)(554.94,0)
\z(554.94,0)(554.94,453.31)(558.43,453.31)(558.43,0)
\z(558.43,0)(558.43,423.02)(561.92,423.02)(561.92,0)
\z(561.92,0)(561.92,426.50)(565.42,426.50)(565.42,0)
\z(565.42,0)(565.42,441.89)(568.91,441.89)(568.91,0)
\z(568.91,0)(568.91,427.49)(572.40,427.49)(572.40,0)
\z(572.40,0)(572.40,439.41)(575.90,439.41)(575.90,0)
\z(575.90,0)(575.90,373.37)(579.39,373.37)(579.39,0)
\z(579.39,0)(579.39,410.61)(582.88,410.61)(582.88,0)
\z(582.88,0)(582.88,417.06)(586.37,417.06)(586.37,0)
\z(586.37,0)(586.37,388.76)(589.87,388.76)(589.87,0)
\z(589.87,0)(589.87,385.78)(593.36,385.78)(593.36,0)
\z(593.36,0)(593.36,376.35)(596.85,376.35)(596.85,0)
\z(596.85,0)(596.85,352.52)(600.35,352.52)(600.35,0)
\z(600.35,0)(600.35,375.85)(603.84,375.85)(603.84,0)
\z(603.84,0)(603.84,348.55)(607.33,348.55)(607.33,0)
\z(607.33,0)(607.33,347.06)(610.83,347.06)(610.83,0)
\z(610.83,0)(610.83,330.17)(614.32,330.17)(614.32,0)
\z(614.32,0)(614.32,314.29)(617.81,314.29)(617.81,0)
\z(617.81,0)(617.81,315.28)(621.30,315.28)(621.30,0)
\z(621.30,0)(621.30,306.34)(624.80,306.34)(624.80,0)
\z(624.80,0)(624.80,301.87)(628.29,301.87)(628.29,0)
\z(628.29,0)(628.29,293.93)(631.78,293.93)(631.78,0)
\z(631.78,0)(631.78,272.08)(635.28,272.08)(635.28,0)
\z(635.28,0)(635.28,299.39)(638.77,299.39)(638.77,0)
\z(638.77,0)(638.77,296.91)(642.26,296.91)(642.26,0)
\z(642.26,0)(642.26,250.24)(645.76,250.24)(645.76,0)
\z(645.76,0)(645.76,307.34)(649.25,307.34)(649.25,0)
\z(649.25,0)(649.25,222.93)(652.74,222.93)(652.74,0)
\z(652.74,0)(652.74,246.76)(656.23,246.76)(656.23,0)
\z(656.23,0)(656.23,231.37)(659.73,231.37)(659.73,0)
\z(659.73,0)(659.73,211.51)(663.22,211.51)(663.22,0)
\z(663.22,0)(663.22,239.81)(666.71,239.81)(666.71,0)
\z(666.71,0)(666.71,225.91)(670.21,225.91)(670.21,0)
\z(670.21,0)(670.21,226.90)(673.70,226.90)(673.70,0)
\z(673.70,0)(673.70,198.60)(677.19,198.60)(677.19,0)
\z(677.19,0)(677.19,194.63)(680.69,194.63)(680.69,0)
\z(680.69,0)(680.69,193.64)(684.18,193.64)(684.18,0)
\z(684.18,0)(684.18,193.64)(687.67,193.64)(687.67,0)
\z(687.67,0)(687.67,196.62)(691.16,196.62)(691.16,0)
\z(691.16,0)(691.16,189.17)(694.66,189.17)(694.66,0)
\z(694.66,0)(694.66,155.41)(698.15,155.41)(698.15,0)
\z(698.15,0)(698.15,155.41)(701.64,155.41)(701.64,0)
\z(701.64,0)(701.64,142.99)(705.14,142.99)(705.14,0)
\z(705.14,0)(705.14,151.43)(708.63,151.43)(708.63,0)
\z(708.63,0)(708.63,136.54)(712.12,136.54)(712.12,0)
\z(712.12,0)(712.12,136.54)(715.62,136.54)(715.62,0)
\z(715.62,0)(715.62,142.99)(719.11,142.99)(719.11,0)
\z(719.11,0)(719.11,143.49)(722.60,143.49)(722.60,0)
\z(722.60,0)(722.60,117.17)(726.09,117.17)(726.09,0)
\z(726.09,0)(726.09,116.68)(729.59,116.68)(729.59,0)
\z(729.59,0)(729.59,102.78)(733.08,102.78)(733.08,0)
\z(733.08,0)(733.08,119.16)(736.57,119.16)(736.57,0)
\z(736.57,0)(736.57,109.73)(740.07,109.73)(740.07,0)
\z(740.07,0)(740.07,93.84)(743.56,93.84)(743.56,0)
\z(743.56,0)(743.56,81.43)(747.05,81.43)(747.05,0)
\z(747.05,0)(747.05,92.85)(750.55,92.85)(750.55,0)
\z(750.55,0)(750.55,88.38)(754.04,88.38)(754.04,0)
\z(754.04,0)(754.04,82.92)(757.53,82.92)(757.53,0)
\z(757.53,0)(757.53,83.91)(761.02,83.91)(761.02,0)
\z(761.02,0)(761.02,84.41)(764.52,84.41)(764.52,0)
\z(764.52,0)(764.52,84.41)(768.01,84.41)(768.01,0)
\z(768.01,0)(768.01,73.98)(771.50,73.98)(771.50,0)
\z(771.50,0)(771.50,69.01)(775.00,69.01)(775.00,0)
\z(775.00,0)(775.00,57.10)(778.49,57.10)(778.49,0)
\z(778.49,0)(778.49,52.63)(781.98,52.63)(781.98,0)
\z(781.98,0)(781.98,54.12)(785.48,54.12)(785.48,0)
\z(785.48,0)(785.48,51.14)(788.97,51.14)(788.97,0)
\z(788.97,0)(788.97,54.12)(792.46,54.12)(792.46,0)
\z(792.46,0)(792.46,47.17)(795.95,47.17)(795.95,0)
\z(795.95,0)(795.95,50.15)(799.45,50.15)(799.45,0)
\z(799.45,0)(799.45,45.18)(802.94,45.18)(802.94,0)
\z(802.94,0)(802.94,33.76)(806.43,33.76)(806.43,0)
\z(806.43,0)(806.43,40.71)(809.93,40.71)(809.93,0)
\z(809.93,0)(809.93,34.26)(813.42,34.26)(813.42,0)
\z(813.42,0)(813.42,34.26)(816.91,34.26)(816.91,0)
\z(816.91,0)(816.91,38.73)(820.41,38.73)(820.41,0)
\z(820.41,0)(820.41,24.83)(823.90,24.83)(823.90,0)
\z(823.90,0)(823.90,25.82)(827.39,25.82)(827.39,0)
\z(827.39,0)(827.39,23.34)(830.88,23.34)(830.88,0)
\z(830.88,0)(830.88,23.34)(834.38,23.34)(834.38,0)
\z(834.38,0)(834.38,22.84)(837.87,22.84)(837.87,0)
\z(837.87,0)(837.87,22.84)(841.36,22.84)(841.36,0)
\z(841.36,0)(841.36,16.88)(844.86,16.88)(844.86,0)
\z(844.86,0)(844.86,20.85)(848.35,20.85)(848.35,0)
\z(848.35,0)(848.35,13.41)(851.84,13.41)(851.84,0)
\z(851.84,0)(851.84,18.87)(855.34,18.87)(855.34,0)
\z(855.34,0)(855.34,16.38)(858.83,16.38)(858.83,0)
\z(858.83,0)(858.83,11.92)(862.32,11.92)(862.32,0)
\z(862.32,0)(862.32,13.41)(865.81,13.41)(865.81,0)
\z(865.81,0)(865.81,10.43)(869.31,10.43)(869.31,0)
\z(869.31,0)(869.31,10.92)(872.80,10.92)(872.80,0)
\z(872.80,0)(872.80,11.42)(876.29,11.42)(876.29,0)
\z(876.29,0)(876.29,5.46)(879.79,5.46)(879.79,0)
\z(879.79,0)(879.79,10.92)(883.28,10.92)(883.28,0)
\z(883.28,0)(883.28,9.43)(886.77,9.43)(886.77,0)
\z(886.77,0)(886.77,7.45)(890.27,7.45)(890.27,0)
\z(890.27,0)(890.27,9.43)(893.76,9.43)(893.76,0)
\z(893.76,0)(893.76,5.96)(897.25,5.96)(897.25,0)
\z(897.25,0)(897.25,1.49)(900.74,1.49)(900.74,0)
\z(900.74,0)(900.74,2.48)(904.24,2.48)(904.24,0)
\z(904.24,0)(904.24,4.97)(907.73,4.97)(907.73,0)
\z(907.73,0)(907.73,3.48)(911.22,3.48)(911.22,0)
\z(911.22,0)(911.22,2.48)(914.72,2.48)(914.72,0)
\z(914.72,0)(914.72,2.48)(918.21,2.48)(918.21,0)
\z(918.21,0)(918.21,2.98)(921.70,2.98)(921.70,0)
\z(921.70,0)(921.70,1.99)(925.20,1.99)(925.20,0)
\z(928.69,0)(928.69,1.49)(932.18,1.49)(932.18,0)
\z(932.18,0)(932.18,1.49)(935.67,1.49)(935.67,0)
\z(935.67,0)(935.67,0.50)(939.17,0.50)(939.17,0)
\z(942.66,0)(942.66,1.49)(946.15,1.49)(946.15,0)
\z(946.15,0)(946.15,0.50)(949.65,0.50)(949.65,0)
\z(949.65,0)(949.65,0.99)(953.14,0.99)(953.14,0)
\z(953.14,0)(953.14,0.50)(956.63,0.50)(956.63,0)
\z(984.58,0)(984.58,0.50)(988.07,0.50)(988.07,0)
\z(1002.04,0)(1002.04,0.50)(1005.53,0.50)(1005.53,0)
\color{black}
\dottedline{12}(0,142.00)(1000,142.00)
\z(0,0)(0,640)(1000,640)(1000,0)(0,0)
\color[rgb]{0.0,0.2,0.8}
\z(1462.56,0)(1462.56,0.26)(1464.40,0.26)(1464.40,0)
\z(1469.91,0)(1469.91,0.26)(1471.75,0.26)(1471.75,0)
\z(1471.75,0)(1471.75,0.26)(1473.58,0.26)(1473.58,0)
\z(1473.58,0)(1473.58,0.52)(1475.42,0.52)(1475.42,0)
\z(1475.42,0)(1475.42,0.26)(1477.26,0.26)(1477.26,0)
\z(1480.93,0)(1480.93,0.26)(1482.76,0.26)(1482.76,0)
\z(1488.27,0)(1488.27,0.26)(1490.11,0.26)(1490.11,0)
\z(1490.11,0)(1490.11,0.78)(1491.95,0.78)(1491.95,0)
\z(1491.95,0)(1491.95,0.26)(1493.78,0.26)(1493.78,0)
\z(1493.78,0)(1493.78,0.52)(1495.62,0.52)(1495.62,0)
\z(1495.62,0)(1495.62,0.52)(1497.46,0.52)(1497.46,0)
\z(1497.46,0)(1497.46,0.78)(1499.29,0.78)(1499.29,0)
\z(1499.29,0)(1499.29,0.78)(1501.13,0.78)(1501.13,0)
\z(1501.13,0)(1501.13,0.52)(1502.97,0.52)(1502.97,0)
\z(1502.97,0)(1502.97,1.04)(1504.80,1.04)(1504.80,0)
\z(1504.80,0)(1504.80,0.78)(1506.64,0.78)(1506.64,0)
\z(1506.64,0)(1506.64,1.31)(1508.47,1.31)(1508.47,0)
\z(1508.47,0)(1508.47,0.78)(1510.31,0.78)(1510.31,0)
\z(1510.31,0)(1510.31,2.35)(1512.15,2.35)(1512.15,0)
\z(1512.15,0)(1512.15,2.35)(1513.98,2.35)(1513.98,0)
\z(1513.98,0)(1513.98,2.61)(1515.82,2.61)(1515.82,0)
\z(1515.82,0)(1515.82,1.04)(1517.66,1.04)(1517.66,0)
\z(1517.66,0)(1517.66,1.83)(1519.49,1.83)(1519.49,0)
\z(1519.49,0)(1519.49,3.13)(1521.33,3.13)(1521.33,0)
\z(1521.33,0)(1521.33,1.83)(1523.17,1.83)(1523.17,0)
\z(1523.17,0)(1523.17,1.57)(1525.00,1.57)(1525.00,0)
\z(1525.00,0)(1525.00,2.35)(1526.84,2.35)(1526.84,0)
\z(1526.84,0)(1526.84,3.13)(1528.67,3.13)(1528.67,0)
\z(1528.67,0)(1528.67,2.09)(1530.51,2.09)(1530.51,0)
\z(1530.51,0)(1530.51,2.35)(1532.35,2.35)(1532.35,0)
\z(1532.35,0)(1532.35,3.92)(1534.18,3.92)(1534.18,0)
\z(1534.18,0)(1534.18,3.39)(1536.02,3.39)(1536.02,0)
\z(1536.02,0)(1536.02,3.65)(1537.86,3.65)(1537.86,0)
\z(1537.86,0)(1537.86,5.22)(1539.69,5.22)(1539.69,0)
\z(1539.69,0)(1539.69,4.18)(1541.53,4.18)(1541.53,0)
\z(1541.53,0)(1541.53,2.61)(1543.37,2.61)(1543.37,0)
\z(1543.37,0)(1543.37,4.18)(1545.20,4.18)(1545.20,0)
\z(1545.20,0)(1545.20,3.39)(1547.04,3.39)(1547.04,0)
\z(1547.04,0)(1547.04,4.70)(1548.87,4.70)(1548.87,0)
\z(1548.87,0)(1548.87,7.05)(1550.71,7.05)(1550.71,0)
\z(1550.71,0)(1550.71,4.96)(1552.55,4.96)(1552.55,0)
\z(1552.55,0)(1552.55,3.92)(1554.38,3.92)(1554.38,0)
\z(1554.38,0)(1554.38,6.26)(1556.22,6.26)(1556.22,0)
\z(1556.22,0)(1556.22,3.65)(1558.06,3.65)(1558.06,0)
\z(1558.06,0)(1558.06,5.22)(1559.89,5.22)(1559.89,0)
\z(1559.89,0)(1559.89,6.26)(1561.73,6.26)(1561.73,0)
\z(1561.73,0)(1561.73,7.83)(1563.57,7.83)(1563.57,0)
\z(1563.57,0)(1563.57,8.61)(1565.40,8.61)(1565.40,0)
\z(1565.40,0)(1565.40,8.61)(1567.24,8.61)(1567.24,0)
\z(1567.24,0)(1567.24,7.83)(1569.08,7.83)(1569.08,0)
\z(1569.08,0)(1569.08,10.18)(1570.91,10.18)(1570.91,0)
\z(1570.91,0)(1570.91,9.92)(1572.75,9.92)(1572.75,0)
\z(1572.75,0)(1572.75,8.09)(1574.58,8.09)(1574.58,0)
\z(1574.58,0)(1574.58,10.18)(1576.42,10.18)(1576.42,0)
\z(1576.42,0)(1576.42,8.35)(1578.26,8.35)(1578.26,0)
\z(1578.26,0)(1578.26,13.05)(1580.09,13.05)(1580.09,0)
\z(1580.09,0)(1580.09,14.88)(1581.93,14.88)(1581.93,0)
\z(1581.93,0)(1581.93,11.49)(1583.77,11.49)(1583.77,0)
\z(1583.77,0)(1583.77,12.53)(1585.60,12.53)(1585.60,0)
\z(1585.60,0)(1585.60,11.75)(1587.44,11.75)(1587.44,0)
\z(1587.44,0)(1587.44,12.53)(1589.28,12.53)(1589.28,0)
\z(1589.28,0)(1589.28,13.83)(1591.11,13.83)(1591.11,0)
\z(1591.11,0)(1591.11,14.88)(1592.95,14.88)(1592.95,0)
\z(1592.95,0)(1592.95,14.10)(1594.78,14.10)(1594.78,0)
\z(1594.78,0)(1594.78,15.14)(1596.62,15.14)(1596.62,0)
\z(1596.62,0)(1596.62,19.58)(1598.46,19.58)(1598.46,0)
\z(1598.46,0)(1598.46,16.71)(1600.29,16.71)(1600.29,0)
\z(1600.29,0)(1600.29,22.71)(1602.13,22.71)(1602.13,0)
\z(1602.13,0)(1602.13,15.92)(1603.97,15.92)(1603.97,0)
\z(1603.97,0)(1603.97,17.49)(1605.80,17.49)(1605.80,0)
\z(1605.80,0)(1605.80,21.93)(1607.64,21.93)(1607.64,0)
\z(1607.64,0)(1607.64,22.19)(1609.48,22.19)(1609.48,0)
\z(1609.48,0)(1609.48,21.67)(1611.31,21.67)(1611.31,0)
\z(1611.31,0)(1611.31,26.89)(1613.15,26.89)(1613.15,0)
\z(1613.15,0)(1613.15,19.58)(1614.99,19.58)(1614.99,0)
\z(1614.99,0)(1614.99,24.80)(1616.82,24.80)(1616.82,0)
\z(1616.82,0)(1616.82,24.54)(1618.66,24.54)(1618.66,0)
\z(1618.66,0)(1618.66,28.45)(1620.49,28.45)(1620.49,0)
\z(1620.49,0)(1620.49,29.24)(1622.33,29.24)(1622.33,0)
\z(1622.33,0)(1622.33,26.89)(1624.17,26.89)(1624.17,0)
\z(1624.17,0)(1624.17,34.98)(1626.00,34.98)(1626.00,0)
\z(1626.00,0)(1626.00,27.41)(1627.84,27.41)(1627.84,0)
\z(1627.84,0)(1627.84,32.89)(1629.68,32.89)(1629.68,0)
\z(1629.68,0)(1629.68,31.32)(1631.51,31.32)(1631.51,0)
\z(1631.51,0)(1631.51,35.24)(1633.35,35.24)(1633.35,0)
\z(1633.35,0)(1633.35,33.41)(1635.19,33.41)(1635.19,0)
\z(1635.19,0)(1635.19,40.72)(1637.02,40.72)(1637.02,0)
\z(1637.02,0)(1637.02,40.20)(1638.86,40.20)(1638.86,0)
\z(1638.86,0)(1638.86,35.50)(1640.69,35.50)(1640.69,0)
\z(1640.69,0)(1640.69,45.94)(1642.53,45.94)(1642.53,0)
\z(1642.53,0)(1642.53,42.81)(1644.37,42.81)(1644.37,0)
\z(1644.37,0)(1644.37,40.46)(1646.20,40.46)(1646.20,0)
\z(1646.20,0)(1646.20,41.24)(1648.04,41.24)(1648.04,0)
\z(1648.04,0)(1648.04,42.81)(1649.88,42.81)(1649.88,0)
\z(1649.88,0)(1649.88,47.51)(1651.71,47.51)(1651.71,0)
\z(1651.71,0)(1651.71,45.42)(1653.55,45.42)(1653.55,0)
\z(1653.55,0)(1653.55,47.51)(1655.39,47.51)(1655.39,0)
\z(1655.39,0)(1655.39,46.46)(1657.22,46.46)(1657.22,0)
\z(1657.22,0)(1657.22,50.64)(1659.06,50.64)(1659.06,0)
\z(1659.06,0)(1659.06,51.16)(1660.90,51.16)(1660.90,0)
\z(1660.90,0)(1660.90,59.25)(1662.73,59.25)(1662.73,0)
\z(1662.73,0)(1662.73,56.64)(1664.57,56.64)(1664.57,0)
\z(1664.57,0)(1664.57,63.17)(1666.40,63.17)(1666.40,0)
\z(1666.40,0)(1666.40,61.86)(1668.24,61.86)(1668.24,0)
\z(1668.24,0)(1668.24,62.39)(1670.08,62.39)(1670.08,0)
\z(1670.08,0)(1670.08,65.00)(1671.91,65.00)(1671.91,0)
\z(1671.91,0)(1671.91,56.64)(1673.75,56.64)(1673.75,0)
\z(1673.75,0)(1673.75,60.82)(1675.59,60.82)(1675.59,0)
\z(1675.59,0)(1675.59,63.17)(1677.42,63.17)(1677.42,0)
\z(1677.42,0)(1677.42,74.39)(1679.26,74.39)(1679.26,0)
\z(1679.26,0)(1679.26,68.91)(1681.10,68.91)(1681.10,0)
\z(1681.10,0)(1681.10,74.13)(1682.93,74.13)(1682.93,0)
\z(1682.93,0)(1682.93,69.96)(1684.77,69.96)(1684.77,0)
\z(1684.77,0)(1684.77,88.23)(1686.60,88.23)(1686.60,0)
\z(1686.60,0)(1686.60,79.61)(1688.44,79.61)(1688.44,0)
\z(1688.44,0)(1688.44,84.31)(1690.28,84.31)(1690.28,0)
\z(1690.28,0)(1690.28,86.92)(1692.11,86.92)(1692.11,0)
\z(1692.11,0)(1692.11,88.75)(1693.95,88.75)(1693.95,0)
\z(1693.95,0)(1693.95,87.97)(1695.79,87.97)(1695.79,0)
\z(1695.79,0)(1695.79,90.06)(1697.62,90.06)(1697.62,0)
\z(1697.62,0)(1697.62,91.88)(1699.46,91.88)(1699.46,0)
\z(1699.46,0)(1699.46,95.54)(1701.30,95.54)(1701.30,0)
\z(1701.30,0)(1701.30,96.06)(1703.13,96.06)(1703.13,0)
\z(1703.13,0)(1703.13,101.80)(1704.97,101.80)(1704.97,0)
\z(1704.97,0)(1704.97,108.85)(1706.81,108.85)(1706.81,0)
\z(1706.81,0)(1706.81,110.42)(1708.64,110.42)(1708.64,0)
\z(1708.64,0)(1708.64,107.02)(1710.48,107.02)(1710.48,0)
\z(1710.48,0)(1710.48,118.51)(1712.31,118.51)(1712.31,0)
\z(1712.31,0)(1712.31,114.59)(1714.15,114.59)(1714.15,0)
\z(1714.15,0)(1714.15,109.11)(1715.99,109.11)(1715.99,0)
\z(1715.99,0)(1715.99,115.64)(1717.82,115.64)(1717.82,0)
\z(1717.82,0)(1717.82,113.29)(1719.66,113.29)(1719.66,0)
\z(1719.66,0)(1719.66,124.51)(1721.50,124.51)(1721.50,0)
\z(1721.50,0)(1721.50,118.51)(1723.33,118.51)(1723.33,0)
\z(1723.33,0)(1723.33,124.51)(1725.17,124.51)(1725.17,0)
\z(1725.17,0)(1725.17,131.82)(1727.01,131.82)(1727.01,0)
\z(1727.01,0)(1727.01,127.64)(1728.84,127.64)(1728.84,0)
\z(1728.84,0)(1728.84,123.99)(1730.68,123.99)(1730.68,0)
\z(1730.68,0)(1730.68,132.34)(1732.51,132.34)(1732.51,0)
\z(1732.51,0)(1732.51,136.52)(1734.35,136.52)(1734.35,0)
\z(1734.35,0)(1734.35,149.31)(1736.19,149.31)(1736.19,0)
\z(1736.19,0)(1736.19,144.35)(1738.02,144.35)(1738.02,0)
\z(1738.02,0)(1738.02,149.05)(1739.86,149.05)(1739.86,0)
\z(1739.86,0)(1739.86,154.79)(1741.70,154.79)(1741.70,0)
\z(1741.70,0)(1741.70,152.70)(1743.53,152.70)(1743.53,0)
\z(1743.53,0)(1743.53,155.83)(1745.37,155.83)(1745.37,0)
\z(1745.37,0)(1745.37,155.31)(1747.21,155.31)(1747.21,0)
\z(1747.21,0)(1747.21,167.06)(1749.04,167.06)(1749.04,0)
\z(1749.04,0)(1749.04,167.58)(1750.88,167.58)(1750.88,0)
\z(1750.88,0)(1750.88,157.14)(1752.72,157.14)(1752.72,0)
\z(1752.72,0)(1752.72,163.40)(1754.55,163.40)(1754.55,0)
\z(1754.55,0)(1754.55,179.07)(1756.39,179.07)(1756.39,0)
\z(1756.39,0)(1756.39,173.58)(1758.22,173.58)(1758.22,0)
\z(1758.22,0)(1758.22,178.02)(1760.06,178.02)(1760.06,0)
\z(1760.06,0)(1760.06,193.68)(1761.90,193.68)(1761.90,0)
\z(1761.90,0)(1761.90,191.60)(1763.73,191.60)(1763.73,0)
\z(1763.73,0)(1763.73,194.73)(1765.57,194.73)(1765.57,0)
\z(1765.57,0)(1765.57,197.60)(1767.41,197.60)(1767.41,0)
\z(1767.41,0)(1767.41,198.90)(1769.24,198.90)(1769.24,0)
\z(1769.24,0)(1769.24,205.43)(1771.08,205.43)(1771.08,0)
\z(1771.08,0)(1771.08,197.86)(1772.92,197.86)(1772.92,0)
\z(1772.92,0)(1772.92,205.43)(1774.75,205.43)(1774.75,0)
\z(1774.75,0)(1774.75,201.51)(1776.59,201.51)(1776.59,0)
\z(1776.59,0)(1776.59,213.52)(1778.42,213.52)(1778.42,0)
\z(1778.42,0)(1778.42,209.87)(1780.26,209.87)(1780.26,0)
\z(1780.26,0)(1780.26,225.01)(1782.10,225.01)(1782.10,0)
\z(1782.10,0)(1782.10,231.53)(1783.93,231.53)(1783.93,0)
\z(1783.93,0)(1783.93,238.06)(1785.77,238.06)(1785.77,0)
\z(1785.77,0)(1785.77,219.53)(1787.61,219.53)(1787.61,0)
\z(1787.61,0)(1787.61,229.71)(1789.44,229.71)(1789.44,0)
\z(1789.44,0)(1789.44,245.37)(1791.28,245.37)(1791.28,0)
\z(1791.28,0)(1791.28,242.24)(1793.12,242.24)(1793.12,0)
\z(1793.12,0)(1793.12,239.10)(1794.95,239.10)(1794.95,0)
\z(1794.95,0)(1794.95,261.03)(1796.79,261.03)(1796.79,0)
\z(1796.79,0)(1796.79,261.03)(1798.62,261.03)(1798.62,0)
\z(1798.62,0)(1798.62,268.08)(1800.46,268.08)(1800.46,0)
\z(1800.46,0)(1800.46,257.11)(1802.30,257.11)(1802.30,0)
\z(1802.30,0)(1802.30,267.03)(1804.13,267.03)(1804.13,0)
\z(1804.13,0)(1804.13,262.60)(1805.97,262.60)(1805.97,0)
\z(1805.97,0)(1805.97,270.69)(1807.81,270.69)(1807.81,0)
\z(1807.81,0)(1807.81,276.95)(1809.64,276.95)(1809.64,0)
\z(1809.64,0)(1809.64,290)(1811.48,290)(1811.48,0)
\z(1811.48,0)(1811.48,280.61)(1813.32,280.61)(1813.32,0)
\z(1813.32,0)(1813.32,299.14)(1815.15,299.14)(1815.15,0)
\z(1815.15,0)(1815.15,297.31)(1816.99,297.31)(1816.99,0)
\z(1816.99,0)(1816.99,289.48)(1818.83,289.48)(1818.83,0)
\z(1818.83,0)(1818.83,298.88)(1820.66,298.88)(1820.66,0)
\z(1820.66,0)(1820.66,313.50)(1822.50,313.50)(1822.50,0)
\z(1822.50,0)(1822.50,312.45)(1824.33,312.45)(1824.33,0)
\z(1824.33,0)(1824.33,303.06)(1826.17,303.06)(1826.17,0)
\z(1826.17,0)(1826.17,328.37)(1828.01,328.37)(1828.01,0)
\z(1828.01,0)(1828.01,332.03)(1829.84,332.03)(1829.84,0)
\z(1829.84,0)(1829.84,329.94)(1831.68,329.94)(1831.68,0)
\z(1831.68,0)(1831.68,316.11)(1833.52,316.11)(1833.52,0)
\z(1833.52,0)(1833.52,320.81)(1835.35,320.81)(1835.35,0)
\z(1835.35,0)(1835.35,339.60)(1837.19,339.60)(1837.19,0)
\z(1837.19,0)(1837.19,348.21)(1839.03,348.21)(1839.03,0)
\z(1839.03,0)(1839.03,352.39)(1840.86,352.39)(1840.86,0)
\z(1840.86,0)(1840.86,336.99)(1842.70,336.99)(1842.70,0)
\z(1842.70,0)(1842.70,340.90)(1844.53,340.90)(1844.53,0)
\z(1844.53,0)(1844.53,364.92)(1846.37,364.92)(1846.37,0)
\z(1846.37,0)(1846.37,368.57)(1848.21,368.57)(1848.21,0)
\z(1848.21,0)(1848.21,373.53)(1850.04,373.53)(1850.04,0)
\z(1850.04,0)(1850.04,367.53)(1851.88,367.53)(1851.88,0)
\z(1851.88,0)(1851.88,373.27)(1853.72,373.27)(1853.72,0)
\z(1853.72,0)(1853.72,375.62)(1855.55,375.62)(1855.55,0)
\z(1855.55,0)(1855.55,382.93)(1857.39,382.93)(1857.39,0)
\z(1857.39,0)(1857.39,394.94)(1859.23,394.94)(1859.23,0)
\z(1859.23,0)(1859.23,398.59)(1861.06,398.59)(1861.06,0)
\z(1861.06,0)(1861.06,406.68)(1862.90,406.68)(1862.90,0)
\z(1862.90,0)(1862.90,401.99)(1864.74,401.99)(1864.74,0)
\z(1864.74,0)(1864.74,403.81)(1866.57,403.81)(1866.57,0)
\z(1866.57,0)(1866.57,402.51)(1868.41,402.51)(1868.41,0)
\z(1868.41,0)(1868.41,418.43)(1870.24,418.43)(1870.24,0)
\z(1870.24,0)(1870.24,421.82)(1872.08,421.82)(1872.08,0)
\z(1872.08,0)(1872.08,435.40)(1873.92,435.40)(1873.92,0)
\z(1873.92,0)(1873.92,429.39)(1875.75,429.39)(1875.75,0)
\z(1875.75,0)(1875.75,418.95)(1877.59,418.95)(1877.59,0)
\z(1877.59,0)(1877.59,440.88)(1879.43,440.88)(1879.43,0)
\z(1879.43,0)(1879.43,430.96)(1881.26,430.96)(1881.26,0)
\z(1881.26,0)(1881.26,444.53)(1883.10,444.53)(1883.10,0)
\z(1883.10,0)(1883.10,455.50)(1884.94,455.50)(1884.94,0)
\z(1884.94,0)(1884.94,434.09)(1886.77,434.09)(1886.77,0)
\z(1886.77,0)(1886.77,446.88)(1888.61,446.88)(1888.61,0)
\z(1888.61,0)(1888.61,451.58)(1890.44,451.58)(1890.44,0)
\z(1890.44,0)(1890.44,456.54)(1892.28,456.54)(1892.28,0)
\z(1892.28,0)(1892.28,466.46)(1894.12,466.46)(1894.12,0)
\z(1894.12,0)(1894.12,469.07)(1895.95,469.07)(1895.95,0)
\z(1895.95,0)(1895.95,462.02)(1897.79,462.02)(1897.79,0)
\z(1897.79,0)(1897.79,482.12)(1899.63,482.12)(1899.63,0)
\z(1899.63,0)(1899.63,464.63)(1901.46,464.63)(1901.46,0)
\z(1901.46,0)(1901.46,464.89)(1903.30,464.89)(1903.30,0)
\z(1903.30,0)(1903.30,454.45)(1905.14,454.45)(1905.14,0)
\z(1905.14,0)(1905.14,505.61)(1906.97,505.61)(1906.97,0)
\z(1906.97,0)(1906.97,476.12)(1908.81,476.12)(1908.81,0)
\z(1908.81,0)(1908.81,469.33)(1910.65,469.33)(1910.65,0)
\z(1910.65,0)(1910.65,490.74)(1912.48,490.74)(1912.48,0)
\z(1912.48,0)(1912.48,503.00)(1914.32,503.00)(1914.32,0)
\z(1914.32,0)(1914.32,480.56)(1916.15,480.56)(1916.15,0)
\z(1916.15,0)(1916.15,497.26)(1917.99,497.26)(1917.99,0)
\z(1917.99,0)(1917.99,497.78)(1919.83,497.78)(1919.83,0)
\z(1919.83,0)(1919.83,500.39)(1921.66,500.39)(1921.66,0)
\z(1921.66,0)(1921.66,472.46)(1923.50,472.46)(1923.50,0)
\z(1923.50,0)(1923.50,509.53)(1925.34,509.53)(1925.34,0)
\z(1925.34,0)(1925.34,491.00)(1927.17,491.00)(1927.17,0)
\z(1927.17,0)(1927.17,498.31)(1929.01,498.31)(1929.01,0)
\z(1929.01,0)(1929.01,522.84)(1930.85,522.84)(1930.85,0)
\z(1930.85,0)(1930.85,487.60)(1932.68,487.60)(1932.68,0)
\z(1932.68,0)(1932.68,469.59)(1934.52,469.59)(1934.52,0)
\z(1934.52,0)(1934.52,492.56)(1936.35,492.56)(1936.35,0)
\z(1936.35,0)(1936.35,479.51)(1938.19,479.51)(1938.19,0)
\z(1938.19,0)(1938.19,475.86)(1940.03,475.86)(1940.03,0)
\z(1940.03,0)(1940.03,492.30)(1941.86,492.30)(1941.86,0)
\z(1941.86,0)(1941.86,484.47)(1943.70,484.47)(1943.70,0)
\z(1943.70,0)(1943.70,472.46)(1945.54,472.46)(1945.54,0)
\z(1945.54,0)(1945.54,456.80)(1947.37,456.80)(1947.37,0)
\z(1947.37,0)(1947.37,475.07)(1949.21,475.07)(1949.21,0)
\z(1949.21,0)(1949.21,465.15)(1951.05,465.15)(1951.05,0)
\z(1951.05,0)(1951.05,458.11)(1952.88,458.11)(1952.88,0)
\z(1952.88,0)(1952.88,458.37)(1954.72,458.37)(1954.72,0)
\z(1954.72,0)(1954.72,463.07)(1956.56,463.07)(1956.56,0)
\z(1956.56,0)(1956.56,447.93)(1958.39,447.93)(1958.39,0)
\z(1958.39,0)(1958.39,452.62)(1960.23,452.62)(1960.23,0)
\z(1960.23,0)(1960.23,427.31)(1962.06,427.31)(1962.06,0)
\z(1962.06,0)(1962.06,454.97)(1963.90,454.97)(1963.90,0)
\z(1963.90,0)(1963.90,426.00)(1965.74,426.00)(1965.74,0)
\z(1965.74,0)(1965.74,441.66)(1967.57,441.66)(1967.57,0)
\z(1967.57,0)(1967.57,427.31)(1969.41,427.31)(1969.41,0)
\z(1969.41,0)(1969.41,448.71)(1971.25,448.71)(1971.25,0)
\z(1971.25,0)(1971.25,434.61)(1973.08,434.61)(1973.08,0)
\z(1973.08,0)(1973.08,424.96)(1974.92,424.96)(1974.92,0)
\z(1974.92,0)(1974.92,400.16)(1976.76,400.16)(1976.76,0)
\z(1976.76,0)(1976.76,406.68)(1978.59,406.68)(1978.59,0)
\z(1978.59,0)(1978.59,409.03)(1980.43,409.03)(1980.43,0)
\z(1980.43,0)(1980.43,421.56)(1982.26,421.56)(1982.26,0)
\z(1982.26,0)(1982.26,408.77)(1984.10,408.77)(1984.10,0)
\z(1984.10,0)(1984.10,413.99)(1985.94,413.99)(1985.94,0)
\z(1985.94,0)(1985.94,404.86)(1987.77,404.86)(1987.77,0)
\z(1987.77,0)(1987.77,388.93)(1989.61,388.93)(1989.61,0)
\z(1989.61,0)(1989.61,381.10)(1991.45,381.10)(1991.45,0)
\z(1991.45,0)(1991.45,374.58)(1993.28,374.58)(1993.28,0)
\z(1993.28,0)(1993.28,370.92)(1995.12,370.92)(1995.12,0)
\z(1995.12,0)(1995.12,389.72)(1996.96,389.72)(1996.96,0)
\z(1996.96,0)(1996.96,368.57)(1998.79,368.57)(1998.79,0)
\z(1998.79,0)(1998.79,354.74)(2000.63,354.74)(2000.63,0)
\z(2000.63,0)(2000.63,343.51)(2002.47,343.51)(2002.47,0)
\z(2002.47,0)(2002.47,371.71)(2004.30,371.71)(2004.30,0)
\z(2004.30,0)(2004.30,364.40)(2006.14,364.40)(2006.14,0)
\z(2006.14,0)(2006.14,347.69)(2007.97,347.69)(2007.97,0)
\z(2007.97,0)(2007.97,345.08)(2009.81,345.08)(2009.81,0)
\z(2009.81,0)(2009.81,342.99)(2011.65,342.99)(2011.65,0)
\z(2011.65,0)(2011.65,319.24)(2013.48,319.24)(2013.48,0)
\z(2013.48,0)(2013.48,326.29)(2015.32,326.29)(2015.32,0)
\z(2015.32,0)(2015.32,314.54)(2017.16,314.54)(2017.16,0)
\z(2017.16,0)(2017.16,325.76)(2018.99,325.76)(2018.99,0)
\z(2018.99,0)(2018.99,308.54)(2020.83,308.54)(2020.83,0)
\z(2020.83,0)(2020.83,310.10)(2022.67,310.10)(2022.67,0)
\z(2022.67,0)(2022.67,312.45)(2024.50,312.45)(2024.50,0)
\z(2024.50,0)(2024.50,308.54)(2026.34,308.54)(2026.34,0)
\z(2026.34,0)(2026.34,293.92)(2028.17,293.92)(2028.17,0)
\z(2028.17,0)(2028.17,283.74)(2030.01,283.74)(2030.01,0)
\z(2030.01,0)(2030.01,308.01)(2031.85,308.01)(2031.85,0)
\z(2031.85,0)(2031.85,288.96)(2033.68,288.96)(2033.68,0)
\z(2033.68,0)(2033.68,276.69)(2035.52,276.69)(2035.52,0)
\z(2035.52,0)(2035.52,281.13)(2037.36,281.13)(2037.36,0)
\z(2037.36,0)(2037.36,288.18)(2039.19,288.18)(2039.19,0)
\z(2039.19,0)(2039.19,283.22)(2041.03,283.22)(2041.03,0)
\z(2041.03,0)(2041.03,281.13)(2042.87,281.13)(2042.87,0)
\z(2042.87,0)(2042.87,253.20)(2044.70,253.20)(2044.70,0)
\z(2044.70,0)(2044.70,265.47)(2046.54,265.47)(2046.54,0)
\z(2046.54,0)(2046.54,270.43)(2048.37,270.43)(2048.37,0)
\z(2048.37,0)(2048.37,248.24)(2050.21,248.24)(2050.21,0)
\z(2050.21,0)(2050.21,249.28)(2052.05,249.28)(2052.05,0)
\z(2052.05,0)(2052.05,231.53)(2053.88,231.53)(2053.88,0)
\z(2053.88,0)(2053.88,251.89)(2055.72,251.89)(2055.72,0)
\z(2055.72,0)(2055.72,240.41)(2057.56,240.41)(2057.56,0)
\z(2057.56,0)(2057.56,238.32)(2059.39,238.32)(2059.39,0)
\z(2059.39,0)(2059.39,233.36)(2061.23,233.36)(2061.23,0)
\z(2061.23,0)(2061.23,219.26)(2063.07,219.26)(2063.07,0)
\z(2063.07,0)(2063.07,234.14)(2064.90,234.14)(2064.90,0)
\z(2064.90,0)(2064.90,225.27)(2066.74,225.27)(2066.74,0)
\z(2066.74,0)(2066.74,232.06)(2068.58,232.06)(2068.58,0)
\z(2068.58,0)(2068.58,214.31)(2070.41,214.31)(2070.41,0)
\z(2070.41,0)(2070.41,210.91)(2072.25,210.91)(2072.25,0)
\z(2072.25,0)(2072.25,215.61)(2074.08,215.61)(2074.08,0)
\z(2074.08,0)(2074.08,203.60)(2075.92,203.60)(2075.92,0)
\z(2075.92,0)(2075.92,197.34)(2077.76,197.34)(2077.76,0)
\z(2077.76,0)(2077.76,198.38)(2079.59,198.38)(2079.59,0)
\z(2079.59,0)(2079.59,188.46)(2081.43,188.46)(2081.43,0)
\z(2081.43,0)(2081.43,204.39)(2083.27,204.39)(2083.27,0)
\z(2083.27,0)(2083.27,180.63)(2085.10,180.63)(2085.10,0)
\z(2085.10,0)(2085.10,192.64)(2086.94,192.64)(2086.94,0)
\z(2086.94,0)(2086.94,182.46)(2088.78,182.46)(2088.78,0)
\z(2088.78,0)(2088.78,195.51)(2090.61,195.51)(2090.61,0)
\z(2090.61,0)(2090.61,176.19)(2092.45,176.19)(2092.45,0)
\z(2092.45,0)(2092.45,175.93)(2094.28,175.93)(2094.28,0)
\z(2094.28,0)(2094.28,161.58)(2096.12,161.58)(2096.12,0)
\z(2096.12,0)(2096.12,160.79)(2097.96,160.79)(2097.96,0)
\z(2097.96,0)(2097.96,152.44)(2099.79,152.44)(2099.79,0)
\z(2099.79,0)(2099.79,155.31)(2101.63,155.31)(2101.63,0)
\z(2101.63,0)(2101.63,155.83)(2103.47,155.83)(2103.47,0)
\z(2103.47,0)(2103.47,145.92)(2105.30,145.92)(2105.30,0)
\z(2105.30,0)(2105.30,155.57)(2107.14,155.57)(2107.14,0)
\z(2107.14,0)(2107.14,144.09)(2108.98,144.09)(2108.98,0)
\z(2108.98,0)(2108.98,146.18)(2110.81,146.18)(2110.81,0)
\z(2110.81,0)(2110.81,140.43)(2112.65,140.43)(2112.65,0)
\z(2112.65,0)(2112.65,140.43)(2114.49,140.43)(2114.49,0)
\z(2114.49,0)(2114.49,132.08)(2116.32,132.08)(2116.32,0)
\z(2116.32,0)(2116.32,136.78)(2118.16,136.78)(2118.16,0)
\z(2118.16,0)(2118.16,127.90)(2119.99,127.90)(2119.99,0)
\z(2119.99,0)(2119.99,140.96)(2121.83,140.96)(2121.83,0)
\z(2121.83,0)(2121.83,123.21)(2123.67,123.21)(2123.67,0)
\z(2123.67,0)(2123.67,119.55)(2125.50,119.55)(2125.50,0)
\z(2125.50,0)(2125.50,126.60)(2127.34,126.60)(2127.34,0)
\z(2127.34,0)(2127.34,120.86)(2129.18,120.86)(2129.18,0)
\z(2129.18,0)(2129.18,112.76)(2131.01,112.76)(2131.01,0)
\z(2131.01,0)(2131.01,100.50)(2132.85,100.50)(2132.85,0)
\z(2132.85,0)(2132.85,98.93)(2134.69,98.93)(2134.69,0)
\z(2134.69,0)(2134.69,112.24)(2136.52,112.24)(2136.52,0)
\z(2136.52,0)(2136.52,109.89)(2138.36,109.89)(2138.36,0)
\z(2138.36,0)(2138.36,109.89)(2140.19,109.89)(2140.19,0)
\z(2140.19,0)(2140.19,99.19)(2142.03,99.19)(2142.03,0)
\z(2142.03,0)(2142.03,105.19)(2143.87,105.19)(2143.87,0)
\z(2143.87,0)(2143.87,87.97)(2145.70,87.97)(2145.70,0)
\z(2145.70,0)(2145.70,90.32)(2147.54,90.32)(2147.54,0)
\z(2147.54,0)(2147.54,94.75)(2149.38,94.75)(2149.38,0)
\z(2149.38,0)(2149.38,99.45)(2151.21,99.45)(2151.21,0)
\z(2151.21,0)(2151.21,85.36)(2153.05,85.36)(2153.05,0)
\z(2153.05,0)(2153.05,91.36)(2154.89,91.36)(2154.89,0)
\z(2154.89,0)(2154.89,85.10)(2156.72,85.10)(2156.72,0)
\z(2156.72,0)(2156.72,85.62)(2158.56,85.62)(2158.56,0)
\z(2158.56,0)(2158.56,82.49)(2160.40,82.49)(2160.40,0)
\z(2160.40,0)(2160.40,81.96)(2162.23,81.96)(2162.23,0)
\z(2162.23,0)(2162.23,83.79)(2164.07,83.79)(2164.07,0)
\z(2164.07,0)(2164.07,77.79)(2165.90,77.79)(2165.90,0)
\z(2165.90,0)(2165.90,71.26)(2167.74,71.26)(2167.74,0)
\z(2167.74,0)(2167.74,75.44)(2169.58,75.44)(2169.58,0)
\z(2169.58,0)(2169.58,68.65)(2171.41,68.65)(2171.41,0)
\z(2171.41,0)(2171.41,65.52)(2173.25,65.52)(2173.25,0)
\z(2173.25,0)(2173.25,65.26)(2175.09,65.26)(2175.09,0)
\z(2175.09,0)(2175.09,62.91)(2176.92,62.91)(2176.92,0)
\z(2176.92,0)(2176.92,62.39)(2178.76,62.39)(2178.76,0)
\z(2178.76,0)(2178.76,57.69)(2180.60,57.69)(2180.60,0)
\z(2180.60,0)(2180.60,59.51)(2182.43,59.51)(2182.43,0)
\z(2182.43,0)(2182.43,52.99)(2184.27,52.99)(2184.27,0)
\z(2184.27,0)(2184.27,48.29)(2186.10,48.29)(2186.10,0)
\z(2186.10,0)(2186.10,53.25)(2187.94,53.25)(2187.94,0)
\z(2187.94,0)(2187.94,54.56)(2189.78,54.56)(2189.78,0)
\z(2189.78,0)(2189.78,57.43)(2191.61,57.43)(2191.61,0)
\z(2191.61,0)(2191.61,54.03)(2193.45,54.03)(2193.45,0)
\z(2193.45,0)(2193.45,43.85)(2195.29,43.85)(2195.29,0)
\z(2195.29,0)(2195.29,43.07)(2197.12,43.07)(2197.12,0)
\z(2197.12,0)(2197.12,52.47)(2198.96,52.47)(2198.96,0)
\z(2198.96,0)(2198.96,43.85)(2200.80,43.85)(2200.80,0)
\z(2200.80,0)(2200.80,46.99)(2202.63,46.99)(2202.63,0)
\z(2202.63,0)(2202.63,40.20)(2204.47,40.20)(2204.47,0)
\z(2204.47,0)(2204.47,40.46)(2206.31,40.46)(2206.31,0)
\z(2206.31,0)(2206.31,36.28)(2208.14,36.28)(2208.14,0)
\z(2208.14,0)(2208.14,44.64)(2209.98,44.64)(2209.98,0)
\z(2209.98,0)(2209.98,38.11)(2211.81,38.11)(2211.81,0)
\z(2211.81,0)(2211.81,31.32)(2213.65,31.32)(2213.65,0)
\z(2213.65,0)(2213.65,33.15)(2215.49,33.15)(2215.49,0)
\z(2215.49,0)(2215.49,35.24)(2217.32,35.24)(2217.32,0)
\z(2217.32,0)(2217.32,34.46)(2219.16,34.46)(2219.16,0)
\z(2219.16,0)(2219.16,31.85)(2221.00,31.85)(2221.00,0)
\z(2221.00,0)(2221.00,27.93)(2222.83,27.93)(2222.83,0)
\z(2222.83,0)(2222.83,29.50)(2224.67,29.50)(2224.67,0)
\z(2224.67,0)(2224.67,33.15)(2226.51,33.15)(2226.51,0)
\z(2226.51,0)(2226.51,31.58)(2228.34,31.58)(2228.34,0)
\z(2228.34,0)(2228.34,27.93)(2230.18,27.93)(2230.18,0)
\z(2230.18,0)(2230.18,23.49)(2232.01,23.49)(2232.01,0)
\z(2232.01,0)(2232.01,23.49)(2233.85,23.49)(2233.85,0)
\z(2233.85,0)(2233.85,26.10)(2235.69,26.10)(2235.69,0)
\z(2235.69,0)(2235.69,24.01)(2237.52,24.01)(2237.52,0)
\z(2237.52,0)(2237.52,21.40)(2239.36,21.40)(2239.36,0)
\z(2239.36,0)(2239.36,19.06)(2241.20,19.06)(2241.20,0)
\z(2241.20,0)(2241.20,19.58)(2243.03,19.58)(2243.03,0)
\z(2243.03,0)(2243.03,20.62)(2244.87,20.62)(2244.87,0)
\z(2244.87,0)(2244.87,21.40)(2246.71,21.40)(2246.71,0)
\z(2246.71,0)(2246.71,18.53)(2248.54,18.53)(2248.54,0)
\z(2248.54,0)(2248.54,16.71)(2250.38,16.71)(2250.38,0)
\z(2250.38,0)(2250.38,15.66)(2252.22,15.66)(2252.22,0)
\z(2252.22,0)(2252.22,16.71)(2254.05,16.71)(2254.05,0)
\z(2254.05,0)(2254.05,13.83)(2255.89,13.83)(2255.89,0)
\z(2255.89,0)(2255.89,19.06)(2257.72,19.06)(2257.72,0)
\z(2257.72,0)(2257.72,17.49)(2259.56,17.49)(2259.56,0)
\z(2259.56,0)(2259.56,13.57)(2261.40,13.57)(2261.40,0)
\z(2261.40,0)(2261.40,13.57)(2263.23,13.57)(2263.23,0)
\z(2263.23,0)(2263.23,12.79)(2265.07,12.79)(2265.07,0)
\z(2265.07,0)(2265.07,12.27)(2266.91,12.27)(2266.91,0)
\z(2266.91,0)(2266.91,10.18)(2268.74,10.18)(2268.74,0)
\z(2268.74,0)(2268.74,10.18)(2270.58,10.18)(2270.58,0)
\z(2270.58,0)(2270.58,10.96)(2272.42,10.96)(2272.42,0)
\z(2272.42,0)(2272.42,11.75)(2274.25,11.75)(2274.25,0)
\z(2274.25,0)(2274.25,12.01)(2276.09,12.01)(2276.09,0)
\z(2276.09,0)(2276.09,7.05)(2277.92,7.05)(2277.92,0)
\z(2277.92,0)(2277.92,5.74)(2279.76,5.74)(2279.76,0)
\z(2279.76,0)(2279.76,8.35)(2281.60,8.35)(2281.60,0)
\z(2281.60,0)(2281.60,9.14)(2283.43,9.14)(2283.43,0)
\z(2283.43,0)(2283.43,6.53)(2285.27,6.53)(2285.27,0)
\z(2285.27,0)(2285.27,6.26)(2287.11,6.26)(2287.11,0)
\z(2287.11,0)(2287.11,5.48)(2288.94,5.48)(2288.94,0)
\z(2288.94,0)(2288.94,5.48)(2290.78,5.48)(2290.78,0)
\z(2290.78,0)(2290.78,5.48)(2292.62,5.48)(2292.62,0)
\z(2292.62,0)(2292.62,7.83)(2294.45,7.83)(2294.45,0)
\z(2294.45,0)(2294.45,4.70)(2296.29,4.70)(2296.29,0)
\z(2296.29,0)(2296.29,5.22)(2298.12,5.22)(2298.12,0)
\z(2298.12,0)(2298.12,3.39)(2299.96,3.39)(2299.96,0)
\z(2299.96,0)(2299.96,4.18)(2301.80,4.18)(2301.80,0)
\z(2301.80,0)(2301.80,4.70)(2303.63,4.70)(2303.63,0)
\z(2303.63,0)(2303.63,3.92)(2305.47,3.92)(2305.47,0)
\z(2305.47,0)(2305.47,3.92)(2307.31,3.92)(2307.31,0)
\z(2307.31,0)(2307.31,3.13)(2309.14,3.13)(2309.14,0)
\z(2309.14,0)(2309.14,2.87)(2310.98,2.87)(2310.98,0)
\z(2310.98,0)(2310.98,2.09)(2312.82,2.09)(2312.82,0)
\z(2312.82,0)(2312.82,3.65)(2314.65,3.65)(2314.65,0)
\z(2314.65,0)(2314.65,3.92)(2316.49,3.92)(2316.49,0)
\z(2316.49,0)(2316.49,2.87)(2318.33,2.87)(2318.33,0)
\z(2318.33,0)(2318.33,3.92)(2320.16,3.92)(2320.16,0)
\z(2320.16,0)(2320.16,2.87)(2322.00,2.87)(2322.00,0)
\z(2322.00,0)(2322.00,1.31)(2323.83,1.31)(2323.83,0)
\z(2323.83,0)(2323.83,1.83)(2325.67,1.83)(2325.67,0)
\z(2325.67,0)(2325.67,2.35)(2327.51,2.35)(2327.51,0)
\z(2327.51,0)(2327.51,1.04)(2329.34,1.04)(2329.34,0)
\z(2329.34,0)(2329.34,1.57)(2331.18,1.57)(2331.18,0)
\z(2331.18,0)(2331.18,1.57)(2333.02,1.57)(2333.02,0)
\z(2333.02,0)(2333.02,1.83)(2334.85,1.83)(2334.85,0)
\z(2334.85,0)(2334.85,0.52)(2336.69,0.52)(2336.69,0)
\z(2336.69,0)(2336.69,0.26)(2338.53,0.26)(2338.53,0)
\z(2338.53,0)(2338.53,0.78)(2340.36,0.78)(2340.36,0)
\z(2340.36,0)(2340.36,0.52)(2342.20,0.52)(2342.20,0)
\z(2342.20,0)(2342.20,1.57)(2344.03,1.57)(2344.03,0)
\z(2344.03,0)(2344.03,0.78)(2345.87,0.78)(2345.87,0)
\z(2347.71,0)(2347.71,0.78)(2349.54,0.78)(2349.54,0)
\z(2349.54,0)(2349.54,0.52)(2351.38,0.52)(2351.38,0)
\z(2353.22,0)(2353.22,0.52)(2355.05,0.52)(2355.05,0)
\z(2356.89,0)(2356.89,0.78)(2358.73,0.78)(2358.73,0)
\z(2358.73,0)(2358.73,0.26)(2360.56,0.26)(2360.56,0)
\z(2371.58,0)(2371.58,0.26)(2373.42,0.26)(2373.42,0)
\z(2373.42,0)(2373.42,0.78)(2375.25,0.78)(2375.25,0)
\z(2384.44,0)(2384.44,0.26)(2386.27,0.26)(2386.27,0)
\z(2402.80,0)(2402.80,0.26)(2404.64,0.26)(2404.64,0)
\color{black}
\dottedline{12}(1400,142.00)(2400,142.00)
\z(1400,0)(1400,640)(2400,640)(2400,0)(1400,0)
\end{picture}

%% file: g3_a1_11_13.tex
\hspace{-30pt}
\begin{picture}(2400,640)
\color[rgb]{0.0,0.2,0.8}
\z(259.06,0)(259.06,12.53)(317.82,12.53)(317.82,0)
\z(317.82,0)(317.82,43.85)(376.59,43.85)(376.59,0)
\z(376.59,0)(376.59,194.21)(435.35,194.21)(435.35,0)
\z(435.35,0)(435.35,419.74)(494.12,419.74)(494.12,0)
\z(494.12,0)(494.12,494.91)(552.88,494.91)(552.88,0)
\z(552.88,0)(552.88,451.06)(611.65,451.06)(611.65,0)
\z(611.65,0)(611.65,175.41)(670.41,175.41)(670.41,0)
\z(670.41,0)(670.41,68.91)(729.18,68.91)(729.18,0)
\z(729.18,0)(729.18,12.53)(787.94,12.53)(787.94,0)
\color{black}
\dottedline{12}(0,106.50)(1000,106.50)
\z(0,0)(0,640)(1000,640)(1000,0)(0,0)
\color[rgb]{0.0,0.2,0.8}
\z(1642.53,0)(1642.53,3.33)(1673.75,3.33)(1673.75,0)
\z(1673.75,0)(1673.75,6.66)(1704.97,6.66)(1704.97,0)
\z(1704.97,0)(1704.97,9.98)(1736.19,9.98)(1736.19,0)
\z(1736.19,0)(1736.19,66.56)(1767.41,66.56)(1767.41,0)
\z(1767.41,0)(1767.41,123.14)(1798.62,123.14)(1798.62,0)
\z(1798.62,0)(1798.62,203.02)(1829.84,203.02)(1829.84,0)
\z(1829.84,0)(1829.84,336.14)(1861.06,336.14)(1861.06,0)
\z(1861.06,0)(1861.06,429.33)(1892.28,429.33)(1892.28,0)
\z(1892.28,0)(1892.28,535.83)(1923.50,535.83)(1923.50,0)
\z(1923.50,0)(1923.50,419.34)(1954.72,419.34)(1954.72,0)
\z(1954.72,0)(1954.72,445.97)(1985.94,445.97)(1985.94,0)
\z(1985.94,0)(1985.94,389.39)(2017.16,389.39)(2017.16,0)
\z(2017.16,0)(2017.16,216.33)(2048.38,216.33)(2048.38,0)
\z(2048.38,0)(2048.38,106.50)(2079.59,106.50)(2079.59,0)
\z(2079.59,0)(2079.59,43.27)(2110.81,43.27)(2110.81,0)
\z(2110.81,0)(2110.81,46.59)(2142.03,46.59)(2142.03,0)
\z(2142.03,0)(2142.03,3.33)(2173.25,3.33)(2173.25,0)
\z(2173.25,0)(2173.25,3.33)(2204.47,3.33)(2204.47,0)
\color{black}
\dottedline{12}(1400,106.50)(2400,106.50)
\z(1400,0)(1400,640)(2400,640)(2400,0)(1400,0)
\end{picture}

%% file: g3_a1_15_17.tex
\hspace{-30pt}
\begin{picture}(2400,640)
\color[rgb]{0.0,0.2,0.8}
\z(227.19,0)(227.19,3.61)(244.12,3.61)(244.12,0)
\z(261.05,0)(261.05,12.64)(277.98,12.64)(277.98,0)
\z(277.98,0)(277.98,7.22)(294.92,7.22)(294.92,0)
\z(294.92,0)(294.92,16.25)(311.85,16.25)(311.85,0)
\z(311.85,0)(311.85,36.10)(328.78,36.10)(328.78,0)
\z(328.78,0)(328.78,39.71)(345.71,39.71)(345.71,0)
\z(345.71,0)(345.71,83.03)(362.64,83.03)(362.64,0)
\z(362.64,0)(362.64,74.01)(379.58,74.01)(379.58,0)
\z(379.58,0)(379.58,138.99)(396.51,138.99)(396.51,0)
\z(396.51,0)(396.51,171.48)(413.44,171.48)(413.44,0)
\z(413.44,0)(413.44,223.83)(430.37,223.83)(430.37,0)
\z(430.37,0)(430.37,317.69)(447.31,317.69)(447.31,0)
\z(447.31,0)(447.31,409.75)(464.24,409.75)(464.24,0)
\z(464.24,0)(464.24,426.00)(481.17,426.00)(481.17,0)
\z(481.17,0)(481.17,483.76)(498.10,483.76)(498.10,0)
\z(498.10,0)(498.10,476.54)(515.03,476.54)(515.03,0)
\z(515.03,0)(515.03,472.93)(531.97,472.93)(531.97,0)
\z(531.97,0)(531.97,507.23)(548.90,507.23)(548.90,0)
\z(548.90,0)(548.90,481.96)(565.83,481.96)(565.83,0)
\z(565.83,0)(565.83,436.83)(582.76,436.83)(582.76,0)
\z(582.76,0)(582.76,404.34)(599.69,404.34)(599.69,0)
\z(599.69,0)(599.69,303.25)(616.63,303.25)(616.63,0)
\z(616.63,0)(616.63,216.61)(633.56,216.61)(633.56,0)
\z(633.56,0)(633.56,166.07)(650.49,166.07)(650.49,0)
\z(650.49,0)(650.49,137.19)(667.42,137.19)(667.42,0)
\z(667.42,0)(667.42,104.69)(684.36,104.69)(684.36,0)
\z(684.36,0)(684.36,64.98)(701.29,64.98)(701.29,0)
\z(701.29,0)(701.29,36.10)(718.22,36.10)(718.22,0)
\z(718.22,0)(718.22,32.49)(735.15,32.49)(735.15,0)
\z(735.15,0)(735.15,21.66)(752.08,21.66)(752.08,0)
\z(752.08,0)(752.08,1.81)(769.02,1.81)(769.02,0)
\z(769.02,0)(769.02,5.42)(785.95,5.42)(785.95,0)
\z(785.95,0)(785.95,3.61)(802.88,3.61)(802.88,0)
\z(819.81,0)(819.81,1.81)(836.75,1.81)(836.75,0)
\z(836.75,0)(836.75,1.81)(853.68,1.81)(853.68,0)
\color{black}
\dottedline{12}(0,106.50)(1000,106.50)
\z(0,0)(0,640)(1000,640)(1000,0)(0,0)
\color[rgb]{0.0,0.2,0.8}
\z(1605.64,0)(1605.64,0.97)(1614.72,0.97)(1614.72,0)
\z(1623.80,0)(1623.80,0.97)(1632.88,0.97)(1632.88,0)
\z(1632.88,0)(1632.88,2.90)(1641.96,2.90)(1641.96,0)
\z(1651.05,0)(1651.05,2.90)(1660.13,2.90)(1660.13,0)
\z(1660.13,0)(1660.13,3.87)(1669.21,3.87)(1669.21,0)
\z(1669.21,0)(1669.21,12.59)(1678.29,12.59)(1678.29,0)
\z(1678.29,0)(1678.29,1.94)(1687.37,1.94)(1687.37,0)
\z(1687.37,0)(1687.37,14.52)(1696.45,14.52)(1696.45,0)
\z(1696.45,0)(1696.45,17.43)(1705.54,17.43)(1705.54,0)
\z(1705.54,0)(1705.54,22.27)(1714.62,22.27)(1714.62,0)
\z(1714.62,0)(1714.62,31.95)(1723.70,31.95)(1723.70,0)
\z(1723.70,0)(1723.70,39.70)(1732.78,39.70)(1732.78,0)
\z(1732.78,0)(1732.78,39.70)(1741.86,39.70)(1741.86,0)
\z(1741.86,0)(1741.86,55.19)(1750.95,55.19)(1750.95,0)
\z(1750.95,0)(1750.95,73.58)(1760.03,73.58)(1760.03,0)
\z(1760.03,0)(1760.03,93.91)(1769.11,93.91)(1769.11,0)
\z(1769.11,0)(1769.11,73.58)(1778.19,73.58)(1778.19,0)
\z(1778.19,0)(1778.19,115.21)(1787.27,115.21)(1787.27,0)
\z(1787.27,0)(1787.27,130.70)(1796.35,130.70)(1796.35,0)
\z(1796.35,0)(1796.35,175.24)(1805.44,175.24)(1805.44,0)
\z(1805.44,0)(1805.44,211.06)(1814.52,211.06)(1814.52,0)
\z(1814.52,0)(1814.52,195.57)(1823.60,195.57)(1823.60,0)
\z(1823.60,0)(1823.60,287.55)(1832.68,287.55)(1832.68,0)
\z(1832.68,0)(1832.68,290.45)(1841.76,290.45)(1841.76,0)
\z(1841.76,0)(1841.76,364.04)(1850.85,364.04)(1850.85,0)
\z(1850.85,0)(1850.85,383.40)(1859.93,383.40)(1859.93,0)
\z(1859.93,0)(1859.93,411.48)(1869.01,411.48)(1869.01,0)
\z(1869.01,0)(1869.01,410.51)(1878.09,410.51)(1878.09,0)
\z(1878.09,0)(1878.09,437.62)(1887.17,437.62)(1887.17,0)
\z(1887.17,0)(1887.17,496.68)(1896.25,496.68)(1896.25,0)
\z(1896.25,0)(1896.25,494.74)(1905.34,494.74)(1905.34,0)
\z(1905.34,0)(1905.34,531.53)(1914.42,531.53)(1914.42,0)
\z(1914.42,0)(1914.42,529.60)(1923.50,529.60)(1923.50,0)
\z(1923.50,0)(1923.50,521.85)(1932.58,521.85)(1932.58,0)
\z(1932.58,0)(1932.58,504.42)(1941.66,504.42)(1941.66,0)
\z(1941.66,0)(1941.66,494.74)(1950.75,494.74)(1950.75,0)
\z(1950.75,0)(1950.75,501.52)(1959.83,501.52)(1959.83,0)
\z(1959.83,0)(1959.83,456.01)(1968.91,456.01)(1968.91,0)
\z(1968.91,0)(1968.91,425.03)(1977.99,425.03)(1977.99,0)
\z(1977.99,0)(1977.99,426.00)(1987.07,426.00)(1987.07,0)
\z(1987.07,0)(1987.07,356.29)(1996.15,356.29)(1996.15,0)
\z(1996.15,0)(1996.15,337.90)(2005.24,337.90)(2005.24,0)
\z(2005.24,0)(2005.24,290.45)(2014.32,290.45)(2014.32,0)
\z(2014.32,0)(2014.32,276.90)(2023.40,276.90)(2023.40,0)
\z(2023.40,0)(2023.40,217.84)(2032.48,217.84)(2032.48,0)
\z(2032.48,0)(2032.48,186.86)(2041.56,186.86)(2041.56,0)
\z(2041.56,0)(2041.56,158.78)(2050.65,158.78)(2050.65,0)
\z(2050.65,0)(2050.65,133.61)(2059.73,133.61)(2059.73,0)
\z(2059.73,0)(2059.73,131.67)(2068.81,131.67)(2068.81,0)
\z(2068.81,0)(2068.81,112.31)(2077.89,112.31)(2077.89,0)
\z(2077.89,0)(2077.89,90.04)(2086.97,90.04)(2086.97,0)
\z(2086.97,0)(2086.97,48.41)(2096.05,48.41)(2096.05,0)
\z(2096.05,0)(2096.05,56.15)(2105.14,56.15)(2105.14,0)
\z(2105.14,0)(2105.14,40.66)(2114.22,40.66)(2114.22,0)
\z(2114.22,0)(2114.22,30.01)(2123.30,30.01)(2123.30,0)
\z(2123.30,0)(2123.30,29.05)(2132.38,29.05)(2132.38,0)
\z(2132.38,0)(2132.38,25.17)(2141.46,25.17)(2141.46,0)
\z(2141.46,0)(2141.46,11.62)(2150.55,11.62)(2150.55,0)
\z(2150.55,0)(2150.55,6.78)(2159.63,6.78)(2159.63,0)
\z(2159.63,0)(2159.63,6.78)(2168.71,6.78)(2168.71,0)
\z(2168.71,0)(2168.71,7.75)(2177.79,7.75)(2177.79,0)
\z(2177.79,0)(2177.79,3.87)(2186.87,3.87)(2186.87,0)
\z(2186.87,0)(2186.87,2.90)(2195.95,2.90)(2195.95,0)
\z(2195.95,0)(2195.95,0.97)(2205.04,0.97)(2205.04,0)
\z(2205.04,0)(2205.04,0.97)(2214.12,0.97)(2214.12,0)
\z(2223.20,0)(2223.20,2.90)(2232.28,2.90)(2232.28,0)
\z(2232.28,0)(2232.28,0.97)(2241.36,0.97)(2241.36,0)
\z(2268.61,0)(2268.61,0.97)(2277.69,0.97)(2277.69,0)
\color{black}
\dottedline{12}(1400,106.50)(2400,106.50)
\z(1400,0)(1400,640)(2400,640)(2400,0)(1400,0)
\end{picture}

%% file: g3_a1_19_21.tex
\hspace{-30pt}
\begin{picture}(2400,640)
\color[rgb]{0.0,0.2,0.8}
\z(196.90,0)(196.90,0.51)(201.71,0.51)(201.71,0)
\z(201.71,0)(201.71,1.54)(206.51,1.54)(206.51,0)
\z(206.51,0)(206.51,1.02)(211.31,1.02)(211.31,0)
\z(220.92,0)(220.92,1.02)(225.72,1.02)(225.72,0)
\z(225.72,0)(225.72,0.51)(230.52,0.51)(230.52,0)
\z(230.52,0)(230.52,0.51)(235.33,0.51)(235.33,0)
\z(235.33,0)(235.33,2.05)(240.13,2.05)(240.13,0)
\z(240.13,0)(240.13,1.02)(244.93,1.02)(244.93,0)
\z(244.93,0)(244.93,2.56)(249.74,2.56)(249.74,0)
\z(249.74,0)(249.74,1.54)(254.54,1.54)(254.54,0)
\z(254.54,0)(254.54,3.58)(259.34,3.58)(259.34,0)
\z(259.34,0)(259.34,3.07)(264.14,3.07)(264.14,0)
\z(264.14,0)(264.14,6.14)(268.95,6.14)(268.95,0)
\z(268.95,0)(268.95,6.14)(273.75,6.14)(273.75,0)
\z(273.75,0)(273.75,10.24)(278.55,10.24)(278.55,0)
\z(278.55,0)(278.55,3.07)(283.36,3.07)(283.36,0)
\z(283.36,0)(283.36,8.19)(288.16,8.19)(288.16,0)
\z(288.16,0)(288.16,8.70)(292.96,8.70)(292.96,0)
\z(292.96,0)(292.96,14.34)(297.76,14.34)(297.76,0)
\z(297.76,0)(297.76,16.38)(302.57,16.38)(302.57,0)
\z(302.57,0)(302.57,16.38)(307.37,16.38)(307.37,0)
\z(307.37,0)(307.37,15.87)(312.17,15.87)(312.17,0)
\z(312.17,0)(312.17,17.92)(316.98,17.92)(316.98,0)
\z(316.98,0)(316.98,25.60)(321.78,25.60)(321.78,0)
\z(321.78,0)(321.78,34.82)(326.58,34.82)(326.58,0)
\z(326.58,0)(326.58,41.47)(331.38,41.47)(331.38,0)
\z(331.38,0)(331.38,37.38)(336.19,37.38)(336.19,0)
\z(336.19,0)(336.19,37.89)(340.99,37.89)(340.99,0)
\z(340.99,0)(340.99,50.69)(345.79,50.69)(345.79,0)
\z(345.79,0)(345.79,55.81)(350.60,55.81)(350.60,0)
\z(350.60,0)(350.60,59.91)(355.40,59.91)(355.40,0)
\z(355.40,0)(355.40,78.85)(360.20,78.85)(360.20,0)
\z(360.20,0)(360.20,89.09)(365.00,89.09)(365.00,0)
\z(365.00,0)(365.00,82.95)(369.81,82.95)(369.81,0)
\z(369.81,0)(369.81,92.68)(374.61,92.68)(374.61,0)
\z(374.61,0)(374.61,112.64)(379.41,112.64)(379.41,0)
\z(379.41,0)(379.41,128.00)(384.22,128.00)(384.22,0)
\z(384.22,0)(384.22,123.40)(389.02,123.40)(389.02,0)
\z(389.02,0)(389.02,138.76)(393.82,138.76)(393.82,0)
\z(393.82,0)(393.82,153.09)(398.63,153.09)(398.63,0)
\z(398.63,0)(398.63,177.67)(403.43,177.67)(403.43,0)
\z(403.43,0)(403.43,187.91)(408.23,187.91)(408.23,0)
\z(408.23,0)(408.23,205.32)(413.03,205.32)(413.03,0)
\z(413.03,0)(413.03,213.00)(417.84,213.00)(417.84,0)
\z(417.84,0)(417.84,232.97)(422.64,232.97)(422.64,0)
\z(422.64,0)(422.64,252.94)(427.44,252.94)(427.44,0)
\z(427.44,0)(427.44,271.37)(432.25,271.37)(432.25,0)
\z(432.25,0)(432.25,284.17)(437.05,284.17)(437.05,0)
\z(437.05,0)(437.05,302.60)(441.85,302.60)(441.85,0)
\z(441.85,0)(441.85,350.22)(446.65,350.22)(446.65,0)
\z(446.65,0)(446.65,346.64)(451.46,346.64)(451.46,0)
\z(451.46,0)(451.46,377.87)(456.26,377.87)(456.26,0)
\z(456.26,0)(456.26,363.02)(461.06,363.02)(461.06,0)
\z(461.06,0)(461.06,403.98)(465.87,403.98)(465.87,0)
\z(465.87,0)(465.87,392.72)(470.67,392.72)(470.67,0)
\z(470.67,0)(470.67,416.27)(475.47,416.27)(475.47,0)
\z(475.47,0)(475.47,412.69)(480.27,412.69)(480.27,0)
\z(480.27,0)(480.27,470.55)(485.08,470.55)(485.08,0)
\z(485.08,0)(485.08,468.50)(489.88,468.50)(489.88,0)
\z(489.88,0)(489.88,484.88)(494.68,484.88)(494.68,0)
\z(494.68,0)(494.68,486.93)(499.49,486.93)(499.49,0)
\z(499.49,0)(499.49,507.41)(504.29,507.41)(504.29,0)
\z(504.29,0)(504.29,486.93)(509.09,486.93)(509.09,0)
\z(509.09,0)(509.09,509.46)(513.89,509.46)(513.89,0)
\z(513.89,0)(513.89,521.24)(518.70,521.24)(518.70,0)
\z(518.70,0)(518.70,514.07)(523.50,514.07)(523.50,0)
\z(523.50,0)(523.50,486.42)(528.30,486.42)(528.30,0)
\z(528.30,0)(528.30,498.19)(533.11,498.19)(533.11,0)
\z(533.11,0)(533.11,490)(537.91,490)(537.91,0)
\z(537.91,0)(537.91,504.34)(542.71,504.34)(542.71,0)
\z(542.71,0)(542.71,492.56)(547.51,492.56)(547.51,0)
\z(547.51,0)(547.51,474.13)(552.32,474.13)(552.32,0)
\z(552.32,0)(552.32,472.59)(557.12,472.59)(557.12,0)
\z(557.12,0)(557.12,459.79)(561.92,459.79)(561.92,0)
\z(561.92,0)(561.92,466.45)(566.73,466.45)(566.73,0)
\z(566.73,0)(566.73,460.31)(571.53,460.31)(571.53,0)
\z(571.53,0)(571.53,411.15)(576.33,411.15)(576.33,0)
\z(576.33,0)(576.33,431.12)(581.13,431.12)(581.13,0)
\z(581.13,0)(581.13,387.60)(585.94,387.60)(585.94,0)
\z(585.94,0)(585.94,402.96)(590.74,402.96)(590.74,0)
\z(590.74,0)(590.74,366.09)(595.54,366.09)(595.54,0)
\z(595.54,0)(595.54,343.56)(600.35,343.56)(600.35,0)
\z(600.35,0)(600.35,332.30)(605.15,332.30)(605.15,0)
\z(605.15,0)(605.15,307.21)(609.95,307.21)(609.95,0)
\z(609.95,0)(609.95,287.75)(614.75,287.75)(614.75,0)
\z(614.75,0)(614.75,280.59)(619.56,280.59)(619.56,0)
\z(619.56,0)(619.56,259.08)(624.36,259.08)(624.36,0)
\z(624.36,0)(624.36,223.75)(629.16,223.75)(629.16,0)
\z(629.16,0)(629.16,214.54)(633.97,214.54)(633.97,0)
\z(633.97,0)(633.97,183.81)(638.77,183.81)(638.77,0)
\z(638.77,0)(638.77,189.96)(643.57,189.96)(643.57,0)
\z(643.57,0)(643.57,164.87)(648.38,164.87)(648.38,0)
\z(648.38,0)(648.38,159.75)(653.18,159.75)(653.18,0)
\z(653.18,0)(653.18,155.14)(657.98,155.14)(657.98,0)
\z(657.98,0)(657.98,142.85)(662.78,142.85)(662.78,0)
\z(662.78,0)(662.78,123.40)(667.59,123.40)(667.59,0)
\z(667.59,0)(667.59,103.94)(672.39,103.94)(672.39,0)
\z(672.39,0)(672.39,107.01)(677.19,107.01)(677.19,0)
\z(677.19,0)(677.19,97.80)(682.00,97.80)(682.00,0)
\z(682.00,0)(682.00,79.88)(686.80,79.88)(686.80,0)
\z(686.80,0)(686.80,68.10)(691.60,68.10)(691.60,0)
\z(691.60,0)(691.60,50.69)(696.40,50.69)(696.40,0)
\z(696.40,0)(696.40,48.13)(701.21,48.13)(701.21,0)
\z(701.21,0)(701.21,44.03)(706.01,44.03)(706.01,0)
\z(706.01,0)(706.01,37.38)(710.81,37.38)(710.81,0)
\z(710.81,0)(710.81,40.96)(715.62,40.96)(715.62,0)
\z(715.62,0)(715.62,37.38)(720.42,37.38)(720.42,0)
\z(720.42,0)(720.42,23.55)(725.22,23.55)(725.22,0)
\z(725.22,0)(725.22,26.62)(730.02,26.62)(730.02,0)
\z(730.02,0)(730.02,20.48)(734.83,20.48)(734.83,0)
\z(734.83,0)(734.83,27.14)(739.63,27.14)(739.63,0)
\z(739.63,0)(739.63,19.97)(744.43,19.97)(744.43,0)
\z(744.43,0)(744.43,11.78)(749.24,11.78)(749.24,0)
\z(749.24,0)(749.24,13.82)(754.04,13.82)(754.04,0)
\z(754.04,0)(754.04,8.70)(758.84,8.70)(758.84,0)
\z(758.84,0)(758.84,5.12)(763.64,5.12)(763.64,0)
\z(763.64,0)(763.64,8.19)(768.45,8.19)(768.45,0)
\z(768.45,0)(768.45,7.17)(773.25,7.17)(773.25,0)
\z(773.25,0)(773.25,2.56)(778.05,2.56)(778.05,0)
\z(778.05,0)(778.05,2.56)(782.86,2.56)(782.86,0)
\z(782.86,0)(782.86,4.61)(787.66,4.61)(787.66,0)
\z(787.66,0)(787.66,2.05)(792.46,2.05)(792.46,0)
\z(792.46,0)(792.46,3.07)(797.26,3.07)(797.26,0)
\z(797.26,0)(797.26,1.54)(802.07,1.54)(802.07,0)
\z(802.07,0)(802.07,1.02)(806.87,1.02)(806.87,0)
\z(806.87,0)(806.87,0.51)(811.67,0.51)(811.67,0)
\z(811.67,0)(811.67,0.51)(816.48,0.51)(816.48,0)
\z(821.28,0)(821.28,2.56)(826.08,2.56)(826.08,0)
\z(826.08,0)(826.08,1.02)(830.88,1.02)(830.88,0)
\z(830.88,0)(830.88,0.51)(835.69,0.51)(835.69,0)
\z(835.69,0)(835.69,0.51)(840.49,0.51)(840.49,0)
\z(840.49,0)(840.49,0.51)(845.29,0.51)(845.29,0)
\z(854.90,0)(854.90,0.51)(859.70,0.51)(859.70,0)
\z(874.11,0)(874.11,0.51)(878.91,0.51)(878.91,0)
\color{black}
\dottedline{12}(0,106.50)(1000,106.50)
\z(0,0)(0,640)(1000,640)(1000,0)(0,0)
\color[rgb]{0.0,0.2,0.8}
\z(1543.17,0)(1543.17,0.27)(1545.71,0.27)(1545.71,0)
\z(1568.53,0)(1568.53,0.27)(1571.06,0.27)(1571.06,0)
\z(1581.20,0)(1581.20,0.27)(1583.74,0.27)(1583.74,0)
\z(1593.88,0)(1593.88,0.27)(1596.42,0.27)(1596.42,0)
\z(1596.42,0)(1596.42,1.08)(1598.95,1.08)(1598.95,0)
\z(1598.95,0)(1598.95,0.54)(1601.49,0.54)(1601.49,0)
\z(1601.49,0)(1601.49,0.27)(1604.02,0.27)(1604.02,0)
\z(1604.02,0)(1604.02,0.54)(1606.56,0.54)(1606.56,0)
\z(1606.56,0)(1606.56,0.27)(1609.09,0.27)(1609.09,0)
\z(1609.09,0)(1609.09,0.54)(1611.63,0.54)(1611.63,0)
\z(1611.63,0)(1611.63,0.27)(1614.16,0.27)(1614.16,0)
\z(1616.70,0)(1616.70,0.27)(1619.24,0.27)(1619.24,0)
\z(1619.24,0)(1619.24,0.54)(1621.77,0.54)(1621.77,0)
\z(1621.77,0)(1621.77,0.81)(1624.31,0.81)(1624.31,0)
\z(1624.31,0)(1624.31,0.54)(1626.84,0.54)(1626.84,0)
\z(1626.84,0)(1626.84,1.08)(1629.38,1.08)(1629.38,0)
\z(1629.38,0)(1629.38,1.08)(1631.91,1.08)(1631.91,0)
\z(1631.91,0)(1631.91,1.35)(1634.45,1.35)(1634.45,0)
\z(1634.45,0)(1634.45,0.81)(1636.98,0.81)(1636.98,0)
\z(1636.98,0)(1636.98,1.62)(1639.52,1.62)(1639.52,0)
\z(1639.52,0)(1639.52,1.89)(1642.06,1.89)(1642.06,0)
\z(1642.06,0)(1642.06,0.54)(1644.59,0.54)(1644.59,0)
\z(1644.59,0)(1644.59,2.16)(1647.13,2.16)(1647.13,0)
\z(1647.13,0)(1647.13,3.24)(1649.66,3.24)(1649.66,0)
\z(1649.66,0)(1649.66,2.43)(1652.20,2.43)(1652.20,0)
\z(1652.20,0)(1652.20,2.43)(1654.73,2.43)(1654.73,0)
\z(1654.73,0)(1654.73,4.05)(1657.27,4.05)(1657.27,0)
\z(1657.27,0)(1657.27,3.51)(1659.80,3.51)(1659.80,0)
\z(1659.80,0)(1659.80,3.51)(1662.34,3.51)(1662.34,0)
\z(1662.34,0)(1662.34,4.60)(1664.88,4.60)(1664.88,0)
\z(1664.88,0)(1664.88,3.78)(1667.41,3.78)(1667.41,0)
\z(1667.41,0)(1667.41,5.95)(1669.95,5.95)(1669.95,0)
\z(1669.95,0)(1669.95,5.95)(1672.48,5.95)(1672.48,0)
\z(1672.48,0)(1672.48,4.32)(1675.02,4.32)(1675.02,0)
\z(1675.02,0)(1675.02,7.84)(1677.55,7.84)(1677.55,0)
\z(1677.55,0)(1677.55,8.38)(1680.09,8.38)(1680.09,0)
\z(1680.09,0)(1680.09,6.76)(1682.62,6.76)(1682.62,0)
\z(1682.62,0)(1682.62,7.03)(1685.16,7.03)(1685.16,0)
\z(1685.16,0)(1685.16,9.73)(1687.70,9.73)(1687.70,0)
\z(1687.70,0)(1687.70,10)(1690.23,10)(1690.23,0)
\z(1690.23,0)(1690.23,11.62)(1692.77,11.62)(1692.77,0)
\z(1692.77,0)(1692.77,8.65)(1695.30,8.65)(1695.30,0)
\z(1695.30,0)(1695.30,11.35)(1697.84,11.35)(1697.84,0)
\z(1697.84,0)(1697.84,12.70)(1700.37,12.70)(1700.37,0)
\z(1700.37,0)(1700.37,18.11)(1702.91,18.11)(1702.91,0)
\z(1702.91,0)(1702.91,12.97)(1705.44,12.97)(1705.44,0)
\z(1705.44,0)(1705.44,15.68)(1707.98,15.68)(1707.98,0)
\z(1707.98,0)(1707.98,19.46)(1710.52,19.46)(1710.52,0)
\z(1710.52,0)(1710.52,20.81)(1713.05,20.81)(1713.05,0)
\z(1713.05,0)(1713.05,23.79)(1715.59,23.79)(1715.59,0)
\z(1715.59,0)(1715.59,21.89)(1718.12,21.89)(1718.12,0)
\z(1718.12,0)(1718.12,26.22)(1720.66,26.22)(1720.66,0)
\z(1720.66,0)(1720.66,26.22)(1723.19,26.22)(1723.19,0)
\z(1723.19,0)(1723.19,32.44)(1725.73,32.44)(1725.73,0)
\z(1725.73,0)(1725.73,30.27)(1728.26,30.27)(1728.26,0)
\z(1728.26,0)(1728.26,31.36)(1730.80,31.36)(1730.80,0)
\z(1730.80,0)(1730.80,35.41)(1733.34,35.41)(1733.34,0)
\z(1733.34,0)(1733.34,37.57)(1735.87,37.57)(1735.87,0)
\z(1735.87,0)(1735.87,43.52)(1738.41,43.52)(1738.41,0)
\z(1738.41,0)(1738.41,36.76)(1740.94,36.76)(1740.94,0)
\z(1740.94,0)(1740.94,45.95)(1743.48,45.95)(1743.48,0)
\z(1743.48,0)(1743.48,55.41)(1746.01,55.41)(1746.01,0)
\z(1746.01,0)(1746.01,50.28)(1748.55,50.28)(1748.55,0)
\z(1748.55,0)(1748.55,52.17)(1751.08,52.17)(1751.08,0)
\z(1751.08,0)(1751.08,54.60)(1753.62,54.60)(1753.62,0)
\z(1753.62,0)(1753.62,69.20)(1756.15,69.20)(1756.15,0)
\z(1756.15,0)(1756.15,70.28)(1758.69,70.28)(1758.69,0)
\z(1758.69,0)(1758.69,75.14)(1761.23,75.14)(1761.23,0)
\z(1761.23,0)(1761.23,75.14)(1763.76,75.14)(1763.76,0)
\z(1763.76,0)(1763.76,93.26)(1766.30,93.26)(1766.30,0)
\z(1766.30,0)(1766.30,90.55)(1768.83,90.55)(1768.83,0)
\z(1768.83,0)(1768.83,82.71)(1771.37,82.71)(1771.37,0)
\z(1771.37,0)(1771.37,95.69)(1773.90,95.69)(1773.90,0)
\z(1773.90,0)(1773.90,101.09)(1776.44,101.09)(1776.44,0)
\z(1776.44,0)(1776.44,116.50)(1778.97,116.50)(1778.97,0)
\z(1778.97,0)(1778.97,109.74)(1781.51,109.74)(1781.51,0)
\z(1781.51,0)(1781.51,117.85)(1784.05,117.85)(1784.05,0)
\z(1784.05,0)(1784.05,116.50)(1786.58,116.50)(1786.58,0)
\z(1786.58,0)(1786.58,141.10)(1789.12,141.10)(1789.12,0)
\z(1789.12,0)(1789.12,127.31)(1791.65,127.31)(1791.65,0)
\z(1791.65,0)(1791.65,148.40)(1794.19,148.40)(1794.19,0)
\z(1794.19,0)(1794.19,164.35)(1796.72,164.35)(1796.72,0)
\z(1796.72,0)(1796.72,165.43)(1799.26,165.43)(1799.26,0)
\z(1799.26,0)(1799.26,167.32)(1801.79,167.32)(1801.79,0)
\z(1801.79,0)(1801.79,172.18)(1804.33,172.18)(1804.33,0)
\z(1804.33,0)(1804.33,190.56)(1806.87,190.56)(1806.87,0)
\z(1806.87,0)(1806.87,190.56)(1809.40,190.56)(1809.40,0)
\z(1809.40,0)(1809.40,207.59)(1811.94,207.59)(1811.94,0)
\z(1811.94,0)(1811.94,216.78)(1814.47,216.78)(1814.47,0)
\z(1814.47,0)(1814.47,215.70)(1817.01,215.70)(1817.01,0)
\z(1817.01,0)(1817.01,238.41)(1819.54,238.41)(1819.54,0)
\z(1819.54,0)(1819.54,243.54)(1822.08,243.54)(1822.08,0)
\z(1822.08,0)(1822.08,253.82)(1824.61,253.82)(1824.61,0)
\z(1824.61,0)(1824.61,265.44)(1827.15,265.44)(1827.15,0)
\z(1827.15,0)(1827.15,267.87)(1829.69,267.87)(1829.69,0)
\z(1829.69,0)(1829.69,279.49)(1832.22,279.49)(1832.22,0)
\z(1832.22,0)(1832.22,283.82)(1834.76,283.82)(1834.76,0)
\z(1834.76,0)(1834.76,291.66)(1837.29,291.66)(1837.29,0)
\z(1837.29,0)(1837.29,305.98)(1839.83,305.98)(1839.83,0)
\z(1839.83,0)(1839.83,309.50)(1842.36,309.50)(1842.36,0)
\z(1842.36,0)(1842.36,322.74)(1844.90,322.74)(1844.90,0)
\z(1844.90,0)(1844.90,347.61)(1847.43,347.61)(1847.43,0)
\z(1847.43,0)(1847.43,350.31)(1849.97,350.31)(1849.97,0)
\z(1849.97,0)(1849.97,342.48)(1852.51,342.48)(1852.51,0)
\z(1852.51,0)(1852.51,385.99)(1855.04,385.99)(1855.04,0)
\z(1855.04,0)(1855.04,358.69)(1857.58,358.69)(1857.58,0)
\z(1857.58,0)(1857.58,364.64)(1860.11,364.64)(1860.11,0)
\z(1860.11,0)(1860.11,391.94)(1862.65,391.94)(1862.65,0)
\z(1862.65,0)(1862.65,400.86)(1865.18,400.86)(1865.18,0)
\z(1865.18,0)(1865.18,408.70)(1867.72,408.70)(1867.72,0)
\z(1867.72,0)(1867.72,394.64)(1870.25,394.64)(1870.25,0)
\z(1870.25,0)(1870.25,405.73)(1872.79,405.73)(1872.79,0)
\z(1872.79,0)(1872.79,425.73)(1875.32,425.73)(1875.32,0)
\z(1875.32,0)(1875.32,428.70)(1877.86,428.70)(1877.86,0)
\z(1877.86,0)(1877.86,419.51)(1880.40,419.51)(1880.40,0)
\z(1880.40,0)(1880.40,445.73)(1882.93,445.73)(1882.93,0)
\z(1882.93,0)(1882.93,459.52)(1885.47,459.52)(1885.47,0)
\z(1885.47,0)(1885.47,466.82)(1888.00,466.82)(1888.00,0)
\z(1888.00,0)(1888.00,489.52)(1890.54,489.52)(1890.54,0)
\z(1890.54,0)(1890.54,473.30)(1893.07,473.30)(1893.07,0)
\z(1893.07,0)(1893.07,464.11)(1895.61,464.11)(1895.61,0)
\z(1895.61,0)(1895.61,481.68)(1898.14,481.68)(1898.14,0)
\z(1898.14,0)(1898.14,497.09)(1900.68,497.09)(1900.68,0)
\z(1900.68,0)(1900.68,504.93)(1903.22,504.93)(1903.22,0)
\z(1903.22,0)(1903.22,474.65)(1905.75,474.65)(1905.75,0)
\z(1905.75,0)(1905.75,497.09)(1908.29,497.09)(1908.29,0)
\z(1908.29,0)(1908.29,505.74)(1910.82,505.74)(1910.82,0)
\z(1910.82,0)(1910.82,508.44)(1913.36,508.44)(1913.36,0)
\z(1913.36,0)(1913.36,509.79)(1915.89,509.79)(1915.89,0)
\z(1915.89,0)(1915.89,497.36)(1918.43,497.36)(1918.43,0)
\z(1918.43,0)(1918.43,516.55)(1920.96,516.55)(1920.96,0)
\z(1920.96,0)(1920.96,497.36)(1923.50,497.36)(1923.50,0)
\z(1923.50,0)(1923.50,512.77)(1926.04,512.77)(1926.04,0)
\z(1926.04,0)(1926.04,512.77)(1928.57,512.77)(1928.57,0)
\z(1928.57,0)(1928.57,493.31)(1931.11,493.31)(1931.11,0)
\z(1931.11,0)(1931.11,510.61)(1933.64,510.61)(1933.64,0)
\z(1933.64,0)(1933.64,510.88)(1936.18,510.88)(1936.18,0)
\z(1936.18,0)(1936.18,500.60)(1938.71,500.60)(1938.71,0)
\z(1938.71,0)(1938.71,498.44)(1941.25,498.44)(1941.25,0)
\z(1941.25,0)(1941.25,508.44)(1943.78,508.44)(1943.78,0)
\z(1943.78,0)(1943.78,507.09)(1946.32,507.09)(1946.32,0)
\z(1946.32,0)(1946.32,492.22)(1948.86,492.22)(1948.86,0)
\z(1948.86,0)(1948.86,476.01)(1951.39,476.01)(1951.39,0)
\z(1951.39,0)(1951.39,490.87)(1953.93,490.87)(1953.93,0)
\z(1953.93,0)(1953.93,467.90)(1956.46,467.90)(1956.46,0)
\z(1956.46,0)(1956.46,478.44)(1959.00,478.44)(1959.00,0)
\z(1959.00,0)(1959.00,467.63)(1961.53,467.63)(1961.53,0)
\z(1961.53,0)(1961.53,476.01)(1964.07,476.01)(1964.07,0)
\z(1964.07,0)(1964.07,459.25)(1966.60,459.25)(1966.60,0)
\z(1966.60,0)(1966.60,464.92)(1969.14,464.92)(1969.14,0)
\z(1969.14,0)(1969.14,444.11)(1971.68,444.11)(1971.68,0)
\z(1971.68,0)(1971.68,419.24)(1974.21,419.24)(1974.21,0)
\z(1974.21,0)(1974.21,413.30)(1976.75,413.30)(1976.75,0)
\z(1976.75,0)(1976.75,415.46)(1979.28,415.46)(1979.28,0)
\z(1979.28,0)(1979.28,404.11)(1981.82,404.11)(1981.82,0)
\z(1981.82,0)(1981.82,388.16)(1984.35,388.16)(1984.35,0)
\z(1984.35,0)(1984.35,382.21)(1986.89,382.21)(1986.89,0)
\z(1986.89,0)(1986.89,383.02)(1989.42,383.02)(1989.42,0)
\z(1989.42,0)(1989.42,357.61)(1991.96,357.61)(1991.96,0)
\z(1991.96,0)(1991.96,368.43)(1994.49,368.43)(1994.49,0)
\z(1994.49,0)(1994.49,361.94)(1997.03,361.94)(1997.03,0)
\z(1997.03,0)(1997.03,332.47)(1999.57,332.47)(1999.57,0)
\z(1999.57,0)(1999.57,343.29)(2002.10,343.29)(2002.10,0)
\z(2002.10,0)(2002.10,327.34)(2004.64,327.34)(2004.64,0)
\z(2004.64,0)(2004.64,315.99)(2007.17,315.99)(2007.17,0)
\z(2007.17,0)(2007.17,313.82)(2009.71,313.82)(2009.71,0)
\z(2009.71,0)(2009.71,281.93)(2012.24,281.93)(2012.24,0)
\z(2012.24,0)(2012.24,306.26)(2014.78,306.26)(2014.78,0)
\z(2014.78,0)(2014.78,278.14)(2017.31,278.14)(2017.31,0)
\z(2017.31,0)(2017.31,270.57)(2019.85,270.57)(2019.85,0)
\z(2019.85,0)(2019.85,266.79)(2022.39,266.79)(2022.39,0)
\z(2022.39,0)(2022.39,264.90)(2024.92,264.90)(2024.92,0)
\z(2024.92,0)(2024.92,243.00)(2027.46,243.00)(2027.46,0)
\z(2027.46,0)(2027.46,220.84)(2029.99,220.84)(2029.99,0)
\z(2029.99,0)(2029.99,214.35)(2032.53,214.35)(2032.53,0)
\z(2032.53,0)(2032.53,202.73)(2035.06,202.73)(2035.06,0)
\z(2035.06,0)(2035.06,200.57)(2037.60,200.57)(2037.60,0)
\z(2037.60,0)(2037.60,187.59)(2040.13,187.59)(2040.13,0)
\z(2040.13,0)(2040.13,188.40)(2042.67,188.40)(2042.67,0)
\z(2042.67,0)(2042.67,176.78)(2045.21,176.78)(2045.21,0)
\z(2045.21,0)(2045.21,168.13)(2047.74,168.13)(2047.74,0)
\z(2047.74,0)(2047.74,166.51)(2050.28,166.51)(2050.28,0)
\z(2050.28,0)(2050.28,162.99)(2052.81,162.99)(2052.81,0)
\z(2052.81,0)(2052.81,140.56)(2055.35,140.56)(2055.35,0)
\z(2055.35,0)(2055.35,145.15)(2057.88,145.15)(2057.88,0)
\z(2057.88,0)(2057.88,137.31)(2060.42,137.31)(2060.42,0)
\z(2060.42,0)(2060.42,133.26)(2062.95,133.26)(2062.95,0)
\z(2062.95,0)(2062.95,120.56)(2065.49,120.56)(2065.49,0)
\z(2065.49,0)(2065.49,108.39)(2068.03,108.39)(2068.03,0)
\z(2068.03,0)(2068.03,104.34)(2070.56,104.34)(2070.56,0)
\z(2070.56,0)(2070.56,105.69)(2073.10,105.69)(2073.10,0)
\z(2073.10,0)(2073.10,98.66)(2075.63,98.66)(2075.63,0)
\z(2075.63,0)(2075.63,99.47)(2078.17,99.47)(2078.17,0)
\z(2078.17,0)(2078.17,93.26)(2080.70,93.26)(2080.70,0)
\z(2080.70,0)(2080.70,80.28)(2083.24,80.28)(2083.24,0)
\z(2083.24,0)(2083.24,82.17)(2085.77,82.17)(2085.77,0)
\z(2085.77,0)(2085.77,76.23)(2088.31,76.23)(2088.31,0)
\z(2088.31,0)(2088.31,71.09)(2090.85,71.09)(2090.85,0)
\z(2090.85,0)(2090.85,59.74)(2093.38,59.74)(2093.38,0)
\z(2093.38,0)(2093.38,54.60)(2095.92,54.60)(2095.92,0)
\z(2095.92,0)(2095.92,58.39)(2098.45,58.39)(2098.45,0)
\z(2098.45,0)(2098.45,45.41)(2100.99,45.41)(2100.99,0)
\z(2100.99,0)(2100.99,53.79)(2103.52,53.79)(2103.52,0)
\z(2103.52,0)(2103.52,49.47)(2106.06,49.47)(2106.06,0)
\z(2106.06,0)(2106.06,45.95)(2108.59,45.95)(2108.59,0)
\z(2108.59,0)(2108.59,39.19)(2111.13,39.19)(2111.13,0)
\z(2111.13,0)(2111.13,37.57)(2113.66,37.57)(2113.66,0)
\z(2113.66,0)(2113.66,41.36)(2116.20,41.36)(2116.20,0)
\z(2116.20,0)(2116.20,35.14)(2118.74,35.14)(2118.74,0)
\z(2118.74,0)(2118.74,26.22)(2121.27,26.22)(2121.27,0)
\z(2121.27,0)(2121.27,28.11)(2123.81,28.11)(2123.81,0)
\z(2123.81,0)(2123.81,28.11)(2126.34,28.11)(2126.34,0)
\z(2126.34,0)(2126.34,26.22)(2128.88,26.22)(2128.88,0)
\z(2128.88,0)(2128.88,22.16)(2131.41,22.16)(2131.41,0)
\z(2131.41,0)(2131.41,22.16)(2133.95,22.16)(2133.95,0)
\z(2133.95,0)(2133.95,16.76)(2136.48,16.76)(2136.48,0)
\z(2136.48,0)(2136.48,22.44)(2139.02,22.44)(2139.02,0)
\z(2139.02,0)(2139.02,18.38)(2141.56,18.38)(2141.56,0)
\z(2141.56,0)(2141.56,15.95)(2144.09,15.95)(2144.09,0)
\z(2144.09,0)(2144.09,14.60)(2146.63,14.60)(2146.63,0)
\z(2146.63,0)(2146.63,14.87)(2149.16,14.87)(2149.16,0)
\z(2149.16,0)(2149.16,11.62)(2151.70,11.62)(2151.70,0)
\z(2151.70,0)(2151.70,14.87)(2154.23,14.87)(2154.23,0)
\z(2154.23,0)(2154.23,13.79)(2156.77,13.79)(2156.77,0)
\z(2156.77,0)(2156.77,9.73)(2159.30,9.73)(2159.30,0)
\z(2159.30,0)(2159.30,7.03)(2161.84,7.03)(2161.84,0)
\z(2161.84,0)(2161.84,7.30)(2164.38,7.30)(2164.38,0)
\z(2164.38,0)(2164.38,9.46)(2166.91,9.46)(2166.91,0)
\z(2166.91,0)(2166.91,6.22)(2169.45,6.22)(2169.45,0)
\z(2169.45,0)(2169.45,8.11)(2171.98,8.11)(2171.98,0)
\z(2171.98,0)(2171.98,4.32)(2174.52,4.32)(2174.52,0)
\z(2174.52,0)(2174.52,4.32)(2177.05,4.32)(2177.05,0)
\z(2177.05,0)(2177.05,3.51)(2179.59,3.51)(2179.59,0)
\z(2179.59,0)(2179.59,2.70)(2182.12,2.70)(2182.12,0)
\z(2182.12,0)(2182.12,4.32)(2184.66,4.32)(2184.66,0)
\z(2184.66,0)(2184.66,5.14)(2187.20,5.14)(2187.20,0)
\z(2187.20,0)(2187.20,3.51)(2189.73,3.51)(2189.73,0)
\z(2189.73,0)(2189.73,2.16)(2192.27,2.16)(2192.27,0)
\z(2192.27,0)(2192.27,2.97)(2194.80,2.97)(2194.80,0)
\z(2194.80,0)(2194.80,2.97)(2197.34,2.97)(2197.34,0)
\z(2197.34,0)(2197.34,2.16)(2199.87,2.16)(2199.87,0)
\z(2199.87,0)(2199.87,1.35)(2202.41,1.35)(2202.41,0)
\z(2202.41,0)(2202.41,1.62)(2204.94,1.62)(2204.94,0)
\z(2204.94,0)(2204.94,1.89)(2207.48,1.89)(2207.48,0)
\z(2207.48,0)(2207.48,2.16)(2210.02,2.16)(2210.02,0)
\z(2210.02,0)(2210.02,0.54)(2212.55,0.54)(2212.55,0)
\z(2215.09,0)(2215.09,1.35)(2217.62,1.35)(2217.62,0)
\z(2217.62,0)(2217.62,0.81)(2220.16,0.81)(2220.16,0)
\z(2220.16,0)(2220.16,0.81)(2222.69,0.81)(2222.69,0)
\z(2222.69,0)(2222.69,1.35)(2225.23,1.35)(2225.23,0)
\z(2225.23,0)(2225.23,0.54)(2227.76,0.54)(2227.76,0)
\z(2227.76,0)(2227.76,1.08)(2230.30,1.08)(2230.30,0)
\z(2230.30,0)(2230.30,0.27)(2232.84,0.27)(2232.84,0)
\z(2232.84,0)(2232.84,0.54)(2235.37,0.54)(2235.37,0)
\z(2235.37,0)(2235.37,0.27)(2237.91,0.27)(2237.91,0)
\z(2237.91,0)(2237.91,0.27)(2240.44,0.27)(2240.44,0)
\z(2240.44,0)(2240.44,0.27)(2242.98,0.27)(2242.98,0)
\z(2253.12,0)(2253.12,0.27)(2255.65,0.27)(2255.65,0)
\z(2255.65,0)(2255.65,0.27)(2258.19,0.27)(2258.19,0)
\z(2258.19,0)(2258.19,0.27)(2260.73,0.27)(2260.73,0)
\z(2265.80,0)(2265.80,0.27)(2268.33,0.27)(2268.33,0)
\z(2275.94,0)(2275.94,0.27)(2278.47,0.27)(2278.47,0)
\color{black}
\dottedline{12}(1400,106.50)(2400,106.50)
\z(1400,0)(1400,640)(2400,640)(2400,0)(1400,0)
\end{picture}

%% file: g1_distrib.tex
\input{g1_a1_s2}

\hspace{48pt}$s_1$\hspace{210pt}$s_2$\\
\vspace{60pt}
\input{g1_s3_s4}

\hspace{48pt}$s_3$\hspace{210pt}$s_4$\\
\vspace{36pt}

%% file: g2_distrib.tex
\input{g2_a1_s2}

\hspace{48pt}$s_1$\hspace{210pt}$s_2$\\
\vspace{36pt}
\input{g2_s3_s4}

\hspace{48pt}$s_3$\hspace{210pt}$s_4$\\
\vspace{36pt}
\input{g2_s5_s6}

\hspace{48pt}$s_5$\hspace{210pt}$s_6$\\
\vspace{36pt}

%% file: g3_distrib.tex
\input{g3_a1_s2}

\hspace{48pt}$s_1$\hspace{210pt}$s_2$\\
\vspace{15pt}
\input{g3_s3_s4}

\hspace{48pt}$s_3$\hspace{210pt}$s_4$\\
\vspace{15pt}
\input{g3_s5_s6}

\hspace{48pt}$s_5$\hspace{210pt}$s_6$\\
\vspace{15pt}
\input{g3_s7_s8}

\hspace{48pt}$s_7$\hspace{210pt}$s_8$

%% file: agct11.bbl
\providecommand{\bysame}{\leavevmode\hbox to3em{\hrulefill}\thinspace}
\providecommand{\MR}{\relax\ifhmode\unskip\space\fi MR }
\providecommand{\MRhref}[2]{%
  \href{http://www.ams.org/mathscinet-getitem?mr=#1}{#2}
}
\providecommand{\href}[2]{#2}

%% file: agct11.bbl
\begin{thebibliography}{10}

\bibitem{Aigner:Catalan}
Martin Aigner, \emph{Catalan and other numbers: a recurrent theme}, Algebraic
  combinatorics and computer science, a tribute to {G}ian-{C}arlo {R}ota,
  Springer, 2001, pp.~347--390.

\bibitem{Akiyama:SatoTateConvergence}
Shigeki Akiyama and Yoshio Tanigawa, \emph{Calculation of values of
  {$L$}-functions associated to elliptic curves}, Mathematics of Computation
  \textbf{68} (1999), no.~227, 1201--1231.

\bibitem{Artin:ZetaFunction}
Emil Artin, \emph{Quadratische {K}\"orper in gebiete der h\"{o}heren
  {K}ongruenzen. {I}, {II}}, Mathematische Zeitschrift \textbf{19} (1924),
  153--246.

\bibitem{Baier:SatoTateLangTrotter}
Stephan Baier, \emph{A remark on the conjectures of {L}ang-{T}rotter and
  {S}ato-{T}ate on average}, 2007, preprint,
  \url{http://arxiv.org/abs/0708.2535v3}.

\bibitem{Baier:SatoTateSmallAngles}
Stephan Baier and Lianyi Zhao, \emph{The {S}ato-{T}ate conjecture on average
  for small angles}, 2007, to appear,
  \url{http://arxiv.org/abs/math/0608318v4}.

\bibitem{Banks:SatoTateOnAverage}
William~D. Banks and Igor~E. Shparlinski, \emph{{S}ato-{T}ate, cyclicity, and
  divisibility statistics on average for elliptic curves of small height},
  2007, to appear, \url{http://arxiv.org/abs/math/0609144}.

\bibitem{Cassels:Prolegomena}
J.~W.~S. Cassels and E.~V. Flynn, \emph{Prolegomena to a middlebrow arithmetic
  of curves of genus 2}, London Mathematical Society Lecture Note Series, vol.
  230, Cambridge University Press, 1996.

\bibitem{Chen:CrossingsNestings}
William~Y.C. Chen, Eva~Y.P. Deng, Rosena~R.X. Du, Richard~P. Stanley, and
  Catherine~H. Yan, \emph{Crossings and nestings of matchings and partitions},
  Transactions of the American Mathematical Society \textbf{359} (2007),
  1555--1575.

\bibitem{Clozel:SatoTate}
Laurent Clozel, Michael Harris, and Richard Taylor, \emph{Automorphy for some
  $l$-adic lifts of automorphic mod $\ell$ {G}alois representations}, 2006,
  preprint, \url{http://www.math.harvard.edu/~rtaylor/twugnew.pdf}.

\bibitem{Cojocaru:SurjectiveGalois}
Alina~Carmen Cojocaru, \emph{On the surjectivity of the {G}alois
  representations associated to non-{CM} elliptic curves}, Canadian Mathematics
  Bulletin \textbf{48} (2005), no.~1, 16--31, With an appendix by Ernst Kani.

\bibitem{Cremona:Database}
John Cremona, \emph{The elliptic curve database for conductors to 130000},
  Algorithmic Number Theory Symposium--{ANTS VII}, Lecture Notes in Computer
  Science, vol. 4076, 2006, pp.~11--29.

\bibitem{Deuring:ComplexMultiplication}
Max Deuring, \emph{Die {K}lassenk\"{o}rper der komplexen {M}ultiplication},
  Enzyklop\"{a}die der mathematischen Wissenschaften, vol. 12 (Book 10, Part
  {II}), B.G. Teubner Verlagsgesellschaft, Stuttgart, 1958.

\bibitem{Cohen:HECHECC}
Henri Cohen~(Ed.) et~al., \emph{Handbook of elliptic and hyperelliptic curve
  cryptography}, Chapman and Hall, 2006.

\bibitem{Faltings:Theorem}
Gerd Faltings, \emph{Endlichkeitss\"{a}tze f\"{u}r abelsche {V}ariet\"{a}ten
  \"{u}ber {Z}ahlk\"{o}rpern. [{F}initeness theorems for abelian varieties over
  number fields]}, Inventiones Mathematicae \textbf{73} (1983), no.~3,
  349--366. \MR{718935 (85g:11026a) reviewed by James Milne.}

\bibitem{Frey:ModularCurves}
Gerhard Frey and Michael M\"{u}ller, \emph{Arithmetic of modular curves and
  applications}, Algorithmic algebra and number theory (Matzat et~al., ed.),
  Springer-Verlag, 1999, pp.~11--48.

\bibitem{Gessel:SymmetricFunctions}
Ira~M. Gessel, \emph{Symmetric functions and {P}-recursiveness}, Journal of
  Combinatorial Theory A \textbf{53} (1990), 257--285.

\bibitem{Gonzalez:ModularGenus2}
Enrique Gonz\'{a}lez-Jim\'{e}nez and Josep Gonz\'{a}lez, \emph{Modular curves
  of genus 2}, Mathematics of Computation \textbf{72} (2003), no.~241,
  397--418.

\bibitem{Grabiner:WeylChamber}
David~J. Grabiner and Peter Magyar, \emph{Random walks in {W}eyl chambers and
  the decomposition of tensor powers}, Journal of Algebraic Combinatorics
  \textbf{2} (1993), no.~3, 239--260.

\bibitem{Guy:LatticePaths}
Richard~K. Guy, Christian Krattenthaler, and Bruce Sagan, \emph{Lattice paths,
  reflections and dimension-changing bijections}, Ars Combinatorica \textbf{34}
  (1992), 3--15.

\bibitem{Hall:OpenImageTheorem}
Chris Hall, \emph{An open image theorem for a general class of abelian
  varieties}, 2008, preprint, \url{http://arxiv.org/abs/0803.1682v1}.

\bibitem{Harris:SatoTateProof}
Michael Harris, Nick Shepherd-Barron, and Richard Taylor, \emph{A family of
  {C}alabi-{Y}au varieties and potential automorphy}, May 2006, preprint,
  \url{http://www.math.harvard.edu/~rtaylor/}.

\bibitem{Hashimoto:SatoTateGenus2}
{Ki-Ichiro} Hashimoto and Hiroshi Tsunogai, \emph{On the {S}ato-{T}ate
  conjecture for {$QM$}-curves of genus two}, Mathematics of Computation
  \textbf{68} (1999), no.~228, 1649--1662.

\bibitem{Howe:TorsionSubgroupsJacobians}
Everett~W. Howe, Franck Lepr\'{e}vost, and Bjorn Poonen, \emph{Large torsion
  subgroups of split jacobians of curves of genus two or three}, Forum
  Mathematicum \textbf{12} (2000), no.~3, 315--364.

\bibitem{Igusa:Genus2Invariants}
Jun-Ichi Igusa, \emph{Arithmetic variety of moduli for genus two}, Annals of
  Mathematics \textbf{2} (1960), no.~72, 612--649.

\bibitem{Katz:LarsensAlternative}
Nicholas~M. Katz, \emph{Larsen's alternative, moments, and the monodromy of
  {L}efschetz pencils}, Contributions to automorphic forms, geometry, and
  number theory, Johns Hopkins University Press, 2004, pp.~521--560.

\bibitem{Katz:PersonalCommunication}
\bysame, 2008, e-mail correspondence.

\bibitem{Katz:RandomMatrices}
Nicholas~M. Katz and Peter Sarnak, \emph{Random matrices, {F}robenius
  eigenvalues, and monodromy}, American Mathematical Society, 1999.

\bibitem{KedlayaSutherland:HyperellipticLSeries}
Kiran~S. Kedlaya and Andrew~V. Sutherland, \emph{Computing {$L$}-series of
  hyperelliptic curves}, Algorithmic Number Theory Symposium--{ANTS VIII},
  Lecture Notes in Computer Science, vol. 5011, Springer, 2008, pp.~312--326.

\bibitem{Knapp:BasicAlgebra}
Anthony~W. Knapp, \emph{Basic algebra}, Birkh\"{a}user, 2006.

\bibitem{Knuth:Notation}
Donald~E. Knuth, \emph{Two notes on notation}, American Mathematical Monthly
  \textbf{99} (1992), no.~5, 403--422.

\bibitem{Koosis:LogarithmicIntegral}
Paul Koosis, \emph{The logarithmic integral {I}}, Cambridge University Press,
  1998.

\bibitem{Kraus:EffectiveBound}
Alain Kraus, \emph{Une remarque sur les points de torsion des courbes
  elliptiques}, C. R. Acad. Sci. Paris S\'{e}r. I Math. \textbf{321} (1995),
  no.~9, 1143--1146.

\bibitem{Kuwata:SplitJacobians}
Masato Kuwata, \emph{Quadratic twists of an elliptic curve and maps from a
  hyperelliptic curve}, Mathematics Journal of Okayama University \textbf{47}
  (2005), 85--97.

\bibitem{Lang:ModularForms}
Serge Lang, \emph{Introduction to modular forms}, Springer, 1976.

\bibitem{Lang:AlgebraicNumberTheory}
\bysame, \emph{Algebraic number theory}, second ed., Springer, 1994.

\bibitem{Lang:FrobeniusDistributions}
Serge Lang and Hale Trotter, \emph{{F}robenius distributions in {$GL_2$}
  extensions}, Lecture Notes in Mathematics, vol. 504, Springer-Verlag, 1976.

\bibitem{Larsen:SatoTate}
Michael Larsen, \emph{The normal distribution as a limit of generalized
  {S}ato-{T}ate measures}, early 1990s, preprint,
  \url{http://mlarsen.math.indiana.edu/~larsen/papers/gauss.pdf}.

\bibitem{Mazur:SatoTate}
Barry Mazur, \emph{Finding meaning in error terms}, Bulletin of the American
  Mathematical Society \textbf{45} (2008), no.~2, 185--228.

\bibitem{Moonen:HodgeClasses}
B.J.J. Moonen and Yu.~G. Zarhin, \emph{Hodge classes on abelian varieties of
  low dimension}, Mathematische Annalen \textbf{315} (1999), 711--733.

\bibitem{Mumford:Genus4Example}
David Mumford, \emph{A note of {S}himura's paper ``{D}iscontinuous groups and
  abelian varities"}, Mathematische Annalen \textbf{181} (1969), 345--351.

\bibitem{Rains:Powers}
Eric~M. Rains, \emph{High powers of random elements of compact {L}ie groups},
  Probability Theory and Related Fields \textbf{107} (1997), 219--241.

\bibitem{Rains:EigenvalueDistribution}
\bysame, \emph{Images of eigenvalue distributions under power maps},
  Probability Theory and Related Fields \textbf{125} (2003), no.~4, 522--538.

\bibitem{Rodriguez-Villegas:SplitCM}
Fernando Rodriguez-Villegas, \emph{Explicit models of genus 2 curves with split
  {CM}}, Algorithmic Number Theory Symposium--{ANTS IV}, Lecture Notes in
  Computer Science, vol. 1838, Springer-Verlag, 2000, pp.~505--514.

\bibitem{Serre:CourseNotes1986}
Jean-Pierre Serre, \emph{Groupes lin\'{e}aires modulo $p$ et points d'ordre
  fini des vari\'{e}t\'{e}s ab\'{e}liennes}, 1986, Coll\`{e}ge de {F}rance
  course notes by Eva Bayer-Fluckiger,
  \url{http://alg-geo.epfl.ch/~bayer/files/Serre-cours.pdf}.

\bibitem{Serre:MTGSp26odd}
\bysame, \emph{Propri\'{e}t\'{e}s conjecturales des groupes de {G}alois
  motiviques et des repr\'{e}sentations l-adiques}, Motives (Seattle, WA,
  1991), Proceedings of Symposia in Pure Mathematics, vol.~55, American
  Mathematical Society, 1994, pp.~377--400. \MR{1265537 (95m:11059)}

\bibitem{Serre:GSp26odd}
\bysame, \emph{Lettre a {M}arie-{F}rance {V}igneras}, Oeuvres, Collected
  Papers, Springer, 2000, pp.~38--55.

\bibitem{Sloane:OEIS}
N.~J.~A. Sloane, \emph{The on-line encyclopedia of integer sequences}, 2007,
  \url{www.research.att.com/~njas/sequences/}.

\bibitem{Smart:Genus2}
N.P. Smart, \emph{{$S$}-unit equations, binary forms and curves of genus 2},
  Proceedings of the London Mathematical Society \textbf{75} (1997), no.~2,
  271--307.

\bibitem{Spanier:AtlasOfFunctions}
Jerome Spanier and Keith~B. Oldham, \emph{An atlas of functions}, Hemisphere,
  1987.

\bibitem{SteinWatkins:Database}
William~A. Stein and Mark Watkins, \emph{A database of elliptic curves--first
  report}, Algorithmic Number Theory Symposium--{ANTS V}, Lecture Notes in
  Computer Science, vol. 2369, Springer, 2002, pp.~267--275.

\bibitem{Sutherland:HasseWitt}
Andrew~V. Sutherland, \emph{Notes on the {H}asse-{W}itt matrix and the density
  of zero {F}robenius traces}, 2008, in preparation.

\bibitem{Tate:SatoTate}
J.~T. Tate, \emph{Algebraic cycles and poles of zeta functions}, Arithmetical
  Algebraic Geometry, Harper and Row, 1965.

\bibitem{Taylor:SatoTate}
Richard Taylor, \emph{Automorphy for some $\ell$-adic lifts of automorphic mod
  $\ell$ {G}alois representations {II}}, 2006, preprint,
  \url{http://www.math.harvard.edu/~rtaylor/twugk6.pdf}.

\bibitem{Wamelen:Genus2CM}
Paul~van Wamelen, \emph{Examples of genus two {CM} curves defined over the
  rationals}, Mathematics of Computation \textbf{68} (1999), no.~225, 307--320.

\bibitem{Wamelen:ProvingCM}
\bysame, \emph{Proving that a genus 2 curve has complex multiplication},
  Mathematics of Computation \textbf{68} (1999), no.~228, 1663--1677.

\bibitem{Weil:ZetaFunction}
Andr\'{e} Weil, \emph{Numbers of solutions of equations in finite fields},
  Bulletin of the American Mathematical Society \textbf{55} (1949), 497--508.

\bibitem{Weyl:ClassicalGroups}
Hermann Weyl, \emph{Classical groups}, second ed., Princeton University Press,
  1946.

\bibitem{Zarhin:HyperellipticNoCM}
Yuri~G. Zarhin, \emph{Hyperelliptic jacobians with complex multiplication},
  Mathematical Research Letters \textbf{7} (2000), no.~1, 123--132.

\end{thebibliography}
